\input amstex
\magnification 1200
\documentstyle{amsppt}
\NoRunningHeads
\NoBlackBoxes
\document

\topmatter
\title 
Group schemes of period $p>2$
\endtitle

\author Victor Abrashkin 
\footnote 
{Partially supported by EPSRC, GR/S72252/01}
\endauthor

\email victor.abrashkin@durham.ac.uk 
\endemail 
\address Maths Dept., Durham University, Sci. Laboratories, 
South Rd., Durham, DH1 3LE, U.K. 
\endaddress 
\abstract For a prime number $p>2$, we give a direct proof 
of Breuil's classification of  
finite flat group schemes killed by $p$ over the valuation ring  
of a $p$-adic field with perfect residue field. As application we
establish a correspondence between finite flat group schemes 
and Faltings's strict modules which respects associated Galois modules
via the Fontaine-Wintenberger field-of-norms functor
\endabstract 

\endtopmatter 

\def\G{\operatorname{G}}

\def\Map{\operatorname{Map}}

\def\Tr{\operatorname{Tr}}
\def\pr{\operatorname{pr}}
\def\Ker{\operatorname{Ker}}

\def\Gal{\operatorname{Gal}}
\def\sep{\operatorname{sep}}
\def\Hom{\operatorname{Hom}}
\def\id{\operatorname{id}}

\def\max{\operatorname{max}}
\def\rk{\operatorname{rk}}

\def\Aut{\operatorname{Aut}}
\def\MG{\operatorname{MG}}
\def\FM{\operatorname{\Cal M\Cal F}}

\def\MF{\operatorname{MF}}
\def\GL{\operatorname{GL}}
\def\Aug{\operatorname{Aug}}

\def\m{\operatorname{m}}
\def\cl{\operatorname{cl}}
\def\Spec{\operatorname{Spec}}
\def\Gr{\operatorname{Gr}}

\def\deg{\operatorname{deg}}

\def\Tr{\operatorname{Tr}}

\def\Ext{\operatorname{Ext}}
\def\Frac{\operatorname{Frac}}
\def\LT{\operatorname{LT}}
\def\Alg{\operatorname{Alg}}
\def\Br{\operatorname{Br}}
\def\MG{\operatorname{M\Gamma }}

\document 

\subhead 0. Introduction 
\endsubhead 
\medskip

Let $p$ be a prime number   
and let $k$ be a perfect field of characteristic $p$. Denote by  
$K_{00}$ the fraction field of the ring of Witt vectors
with coefficients in $k$. Let $K_0$ be a totally ramified 
field extension of $K_{00}$ of degree $e$ and let $O_0$ be the 
valuation ring of $K_0$. Let $\Gr '_{O_0}$ be the category 
of finite flat commutative $p$-group schemes $G$ (i.e. the order $|G|$ 
of $G$ is an integral power of $p$) over $O_0$. We shall use the
notation $\Gr _{O_0}$ for the full subcategory in $\Gr '_{O_0}$
consisting of group schemes killed by $p$ (i.e. such that $p\id
_G=0$). All further notation from this introduction will be carefully
reminded in due course in the main body of the paper. 
\medskip 

\subsubhead{\rm 0.1} Motivation 
\endsubsubhead 

There were various approaches to the problem of description of the
category $\Gr '_{O_0}$; especially should be mentioned 
[Fo2] in the case $e=1$, 
[Co] in the case $e<p-1$, [Br1] and [Ki1-3] in the case of 
arbitrary $e$ (everywhere results are not complete 
if $p=2$). In all these cases the classification 
of group schemes appears in
terms of categories of filtered modules 
and was deduced from the corresponding classification of
$p$-divisible groups. More precisely, the case $e<p-1$ was treated in
[Co] via
Fontaine's results about $p$-divisible groups [Fo1] and the case of
arbitrary $e$ uses essentially either in [Br1] and [Ki1] 
the crystalline Dieudonne theory from [BBM], or in [Ki3] the Fontaine-Messing
theory, or in [Ki2]  
the Zink theory of ``displays and
windows'' [Zi]. 

On the other hand, there is an alternative approach resulted in an
explicit description of algebras of group schemes together with the
corresponding coalgebra structures. On the first place one should
mention two classical papers [TO] and [Ra]. 
They give (in the case of the basic ring $O_0$) 
an explicit description of all simple objects of the category $\Gr
'_{O_0}$ ($\equiv $ simple objects of $\Gr _{O_0}$). 
For small $e$, the author disseminated these results to the whole 
category $\Gr _{O_0}$, cf. [Ab2] for the case $e=1$ and [Ab3] for 
the case
$e\leqslant p-1$. It is worth mentioning that:
a) the case $p=2$ is studied completely in [Ab2]; 
b) under the assumption $e<p-1$ all constructions in [Ab3] become
extremely simple and require, as a matter of fact, only the knowledge
of the classical Dieudonne theory of group schemes over perfect field 
of characteristic $p$.

 Notice that the problem of alternative 
and direct description of objects of $\Gr '_{O_0}$ 
(and especialy of $\Gr _{O_0}$) in the case of
arbitrary $e$ was considered in [Br2]. The reason is that 
Breuil's description of $\Gr '_{O_0}$ in [Br1] appears in a very elegant and
natural way in terms of crystalline sheaves but in the very end all 
crystalline concepts can be successfully eliminated. 
Such simplified interpretation of the 
classification of $p$-divisible groups and finite flat group schemes 
was suggested in [Br2] and achieved in different ways in 
[Ki1-3]. Notice also that it is rather
easy to construct and to prove the full faithfulness of the functor from 
an appropriate category of filtered modules to the category of finite flat
group schemes. The main problem appears when proving that this functor
is essentially surjective. This is where the crystalline (resp. the usual)
Dieudonne theory 
plays a crucial role in [Br1] (resp. [Ab2,3]). Despite of the 
beauty and conceptuality of the crystalline Dieudonne theory 
this looks like a very long way around 
and it would be very interesting to understand what
are the properties of finite flat group schemes we do need to 
establish this surjectivity. (These properties should be implicitly hidden 
in the crystalline Dieudonne theory.)  
Mention also that the surjectivity on the level of group schemes 
killed by $p$ implies immediately the surjectivity for the whole
category of $p$-group schemes and for many applications,  
e.g. [BCDT], we do need 
the knowledge of a complete classification of group schemes only on the level
of objects killed by $p$.

In this paper we extend the approach from [Ab2,3] to the whole
category $\Gr _{O_0}$ with no restrictions on $e$. The basic idea 
can be explained as follows. Suppose $G_0=\Spec A_0\in\Gr _{O_0}$. 
Then one can use the methods from [Ab2,3] (and to some extent from [Ab1])
if there are sufficiently many functions $a\in
A_0$ such that for any 
$g_1,g_2\in G_0(\bar K)$, 
$$a(g_1+g_2)\equiv a(g_1)+a(g_2)
\operatorname{mod }p\bar O.\tag{0.1.1}$$
(Here $\bar K$ is an algebraic closure of $K_0$ and 
$\bar O$ is the valuation ring of $\bar K$.)
Notice that if 
$e\leqslant p-1$ then these functions appear just from  
the classical Dieudonne module $M(\bar G_0)$ of $\bar G_0=G_0\otimes k=\Spec
\bar A_0$. (In this case the elements of $M(\bar G_0)$ appear as
covectors $(a_{-n})_{n\geqslant 0}\in\operatorname{CW}(\bar A_0)$ and
their zero components $a_0\in\bar A_0$ give rise to such functions.) 
In the case of arbitrary $e$ 
such functions generally do
not exist (the example comes easily from extensions of the etale constant
group scheme of order $p$ via the multiplicative constant group scheme
of order $p$) but they do appear if we pass to the extension of scalars 
$G_0\otimes _{O_0}O$, where $O$ is the valuation ring 
of $K=K_0(\pi )$ and $\pi ^p=\pi _0$ is a uniformising element of
$K_0$. This explains why our approach depends on a choice of
uniformising element $\pi _0$ of $O_0$. 
This situation corresponds to the fact that the crystalline Dieudonne theory 
for $G_0\otimes (O_0/p)$ 
provides not just a Dieudonne module but a sheaf of Dieudonne
modules in the $fppf$ (or rather syntomic) topology and this sheaf is not generated by its
global sections. In particular, Breuil's classification of group
schemes over $O_0$ requires the sections over $O_0[\{\pi _0^{p^{-n}}\ |\
n\geqslant 0\}]$ and also depends on a choice of $\pi _0$. 
Vice versa, we start with a suitable Breuil's category of filtered
modules $\MF ^e_S$, cf. the definition below, and apply the methods
from [Ab2,3] to construct the functor $\Cal G_O:\MF
^e_S\longrightarrow\Gr _{O}$. If $G=\Spec A\in\operatorname{Im}\Cal
G_O$ then the $O$-algebra $A$ contains sufficiently many functions
satisfying  above condition (0.1.1). This part works perfectly well
for any valuation ring $O$ but the problem is that the image of $\Cal
G_O$ is considerably less than $\Gr _{O}$. Fortunately, if $O$ appears in the
form $O_0[\root p\of {\pi _0}]$ then (the author studied this from
[Br1]) such group scheme $G$ appears as the extension of scalars of 
a unique $G_0\in\Gr _{O_0}$. Finally, it remains 
to prove that for any $G_0\in\Gr _{O_0}$, the corresponding 
$G=G_0\otimes _{O_0}O$ belongs to the image of the functor $\Cal
G_O$. This is the most difficult part of the paper, where 
the methods from [Ab1] were very helpful. This part can be considered
as a replacement of the crystalline Dieudonne theory in the context of 
group schemes over $O_0$ which are killed by $p$. 

The interrelation between $\Gr _{O_0}$ and the image
$\operatorname{Im}\Cal G_O\subset\Gr _O$ can be illustrated by the
following example. Consider the group 
$\Ext _{\Gr _{O_0}}((\Bbb Z/p)_{O_0},\mu _{p,O_0})$ 
of extensions of the constant etale 
group scheme 
$(\Bbb Z/p)_{O_0}$ via the constant multiplicative group scheme 
$\mu _{p,O_0}$ in $\Gr _{O_0}$. This group of extensions is naturally isomorphic 
to $O_0^*/O_0^{*p}$. The image of $\Cal G_O$ 
gives only the subgroup 
$$(1+pO)^{\times }/(O^{*p}\cap (1+pO)^{\times })\subset  
O^*/O^{*p}\simeq \Ext _{\Gr _O}((\Bbb Z/p)_O,\mu
_{p,O}).$$ 
But the embedding $O_0\subset O$ induces the group isomorphism 
$$O_0^*/O_0^{*p}\simeq (1+pO)^{\times }/(O^{*p}\cap (1+pO)^{\times
}).$$
In addition, the fact that the multiplicative structure on $(1+pO)^{\times }$ can be  
transformed into the additive structure on $pO$ 
via the $p$-adic logarithm explains why the additive filtered modules 
are very helpful to describe  
the structure of $G_0\otimes _{O_0}O$ but can't be used directly for 
$G_0\in\Gr _{O_0}$.

We must notice that our strategy should work for all prime
numbers $p$,  
but in this paper we consider only the case $p>2$. 
The case $p=2$ requires much more careful calculations and the author has
not yet checked all details. But having in mind the results from [Ab2] 
one can expect in the case $p=2$ the classification will be obtained 
 for a slightly
different category $\Gr ^*_{O_0}$ under the additional assumption that $k$
is algebraicly closed. The category $\Gr ^*_{O_0}$ contains the same objects as
$\Gr _{O_0}$ but the morphisms of this category come from 
the morphisms $G_1\longrightarrow
G_2$ in $\Gr _{O_0}$ modulo those which factor through 
the canonical projection to the maximal etale 
quotient $G_1\longrightarrow G_1^{et}$ and the 
embedding of the maximal multiplicative subobject 
$G_2^{m}\longrightarrow G_2$. As a matter of fact, the idea to use the
category $\Gr ^*_O$ allows one to modify the constructions from
Sections 2 and 3 for $p=2$. Then one can adjust the content of 
Sections 1,4 and 5 for arbitrary $p$. The main problem appears with
calculations in Sections 6 and their interpretations in Section 7. 
There are some technical complications, e.g. one has more complicated 
formula for opposite elements in the corresponding Lubin-Tate groups. 
Much more serious problems arise because we use
systematically Lubin-Tate logarithms and modulo $p$ calculations 
should be replaced by calculations modulo 4 if $p=2$. 
Anyway, our approach requires an essential restructuring and the case 
$p=2$ deserves a separate study. 
\medskip 

\subsubhead{\rm 0.2} The main statement 
\endsubsubhead 

As earlier, $O=O_K$, where 
$K=K_0(\pi )$, $\pi ^p=\pi _0$ is a uniformising
element in $K_0$. Suppose $S=k[[t]]$ and   
$\sigma :S\longrightarrow S$ is such that $\sigma (s)=s^p$ for
any $s\in S$. Fix a ring identification 
$\kappa _{SO}:S/t^{ep}S\longrightarrow
O/pO$ such that $\kappa _{SO}|k=\id $. 

Let $\FM _S$ be the category of triples $(M^0,M^1,\varphi _1)$,
where $M^0$ is an $S$-module, $M^1$ is its
submodule and $\varphi _1:M^1\longrightarrow M^0$ 
is a $\sigma $-linear morphism. 
The morphisms
in this category are morphisms of filtered $S$-modules which commute
with the corresponding $\varphi _1$'s. 
By $\MF ^e_S$ we denote 
the full subcategory in $\FM _S$ consisting of $(M^0,M^1,\varphi _1)$
such that $M^0$ is a free $S$-module of finite rank, $M^1\supset
t^eM^0$ and $\varphi _1(M^1)S=M^0$.  

Our main result is the following 

\proclaim{Theorem} There is an antiequivalence of categories 
$\Cal G^O_{O_0}:\MF ^e_S\longrightarrow \Gr _{O_0}$. 
\endproclaim 

Notice that this antiequivalence essentially depends on the choice of 
the field extension $K$ of $K_0$ (the ring 
identification $\kappa _{SO}$ is introduced just for technical reasons). 
It would be very interesting to understand how 
$\Cal G^O_{O_0}$ depends on the choice of $K$. It is worth mentioning
that $\Cal G_{O_0}^O$ coincides with the restriction of Breuil's
antiequivalence to the category $\Gr _{O_0}$, cf. Section 8, but 
in our approach the construction  
of this antiequivalence is quite direct and explicit and
all proofs are given entirely in the limits of the theory of finite
flat group schemes. 

The above Theorem is proved in first 7 sections. 

In Section 1 we introduce the category $\MF ^e_S$ and explain that it
is equivalent to the category of filtered $S_0$-modules with slope
$\leqslant e$, where $S_0=k[[t^p]]\subset S$. 
Then we introduce the concept of
a $\varphi _1$-nilpotent lift of objects of $\FM _S$ to $\MF ^e_S$. 
It is related to the situation, where the knowledge of a quotient 
$\Cal N\in\FM _S$ of $\Cal M\in\MF ^e_S$ is sufficient for a unique 
recovering of 
$\Cal M$ from $\Cal N$. As a matter of
fact, this is the only idea from crystalline cohomology which survives
in our setting. 

In Section 2 we construct the functor 
$\Cal G_O:\MF ^e_S\longrightarrow\Gr _O$ by applying directly the
ideas from [Ab2,3]. For each $\Cal M\in\MF ^e_S$ we construct a family
of explicitly given $O$-algebras $\Cal A(\Cal M)$, for any
$A\in\Cal A(\Cal M)$ provide $G=\Spec A$ with a structure of an object
of the category $\Gr _O$ and prove that all these $G$'s are naturally
isomorphic. A special case of algebras from $\Cal A(\Cal M)$  
plays a very important role in Sections 3 and 6 (it appears also in
Breuil's paper [Br1, Section 3.1]), but in Sections 4 and 7 we do need
more general algebras from $\Cal A(\Cal M)$. 

In Section 3 we prove that $\Cal G_O$ is
fully faithful. Namely, suppose $\Cal M\in\MF ^e_S$ 
and $G=\Cal G_O(\Cal M)=\Spec A$. Then $\Cal M$ can be recovered
uniquely as a $\varphi _1$-nilpotent lift of the following object
$\Cal N$ of the category $\FM _S$. Suppose $e_G:A\longrightarrow O$ 
and $\Delta _G:A\longrightarrow A\otimes A$ are the counit and the
comultiplication of $G$. Set 
$I_A=\Ker e_G$, $I_{A\otimes A}=\Ker e_{G\times G}$ 
and $I_{A\otimes A}(p)=\{a\in I_{A\otimes A}\ |\ a^p\in pA\otimes
A\}$.  Introduce the ideal $I_A^{DP}$ as the maximal ideal in $I_A$
with the structure of nilpotent divided powers. 
(This means $a_0\in I^{DP}_A\Leftrightarrow$ if for all 
$i\geqslant 0$, $a_{i+1}=-a_i^p/p$ then $a_i\to 0$ as $i\to\infty $.) 
Then $\Cal N=(N^0, N^1,\varphi _1)$, where 
(compare with (0.1.1))
$$N^0=\{a\in I_A\ |\ \Delta _G(a)\equiv a\otimes 1+1\otimes
a\operatorname{mod}I_{A\otimes A}(p)^p\}\operatorname{mod}I_A^{DP},\tag{0.2.1}$$
$N^1=(I_A(p)/I_A^{DP})\bigcap N^0$ and $\varphi _1$ is induced by the
correspondence $a\mapsto -a^p/p$ with $a\in I_A(p)$. 
As we have noticed earlier, we introduce here a special way to construct 
the $O$-algebras of group schemes $\Cal G_O(\Cal M)$, $\Cal M\in\MF
^e_S$,  
and define for all $\alpha\in O$, the special ideals 
$I_A(\alpha )$ in $A$ and $I_{A\otimes A}(\alpha )$ in $A\otimes A$. 
These technical notions will play an important role
in Section 6. 

In Section 4 we prove that one can study the image of
$\Cal G_O$ by replacing if necessary the ring $O$ by the valuation
ring of any tamely ramified extension of its fraction field $K$. 
In particular, this will allow us later to treat any $G\in\Gr _O$ as a
result of successive extensions via group schemes of order $p$. 
We
also prove in this section 
that any $G\in\Cal G_O(\Cal M)$ comes via extension of 
scalars from a unique $G_0\in Gr _{O_0}$. 

In section 5 we describe torsors of group schemes of
order $p$ over flat $O$-algebras. As a matter of fact, this is a 
detailed revision of the corresponding result from [Ab1]. 
Then we use this description (again following the strategy from 
[Ab1]) to describe
the group of extensions of $H\in\Gr _O$ via group schemes of order $p$ 
in the category $\Gr _O$. Notice that in this section we do not assume that 
$O$ is obtained in the form $O_0[\pi ]$, where $\pi ^p=\pi _0$.

Section 6 contains the proof of the Main Lemma. This lemma is quite
technical, the calculations are done in the spirit of [Ab1] but are
based on a different background. As we have noticed earlier, this
technical lemma provides us with the fact that $A(G_0)_{O_0}\otimes O$ 
contains sufficiently many functions satisfying condition (0.1.1), or
equivalently, that the crystalline Dieudonne sheaf associated with
$G_0$ contains sufficiently many sections over $O$. 

Finally, in Section 7 we 
apply the results of previous Sections 4-6 to deduce that 
for any $G_0\in Gr _{O_0}$, its extension of scalars $G=G_0\otimes
_{O_0}O$ belongs to the image of the functor $\Cal G_O$. This gives the
existence of a functor $\Cal G_{O_0}^O:\MF ^e_S\longrightarrow\Gr
_{O_0}$ such that for any $\Cal M$, $\Cal G_{O_0}^O(\Cal M)=G_0$, 
and the full faithfulness of $\Cal G_O$ implies that 
$\Cal G_{O_0}^O$ is an antiequivalence of categories. 

In section 8 we give several applications of our methods. 
In particular, we prove that for the same choice of the uniformising
element $\pi _0\in O_0$, our antiequivalence coincides with 
the Breuil antiequivalence restricted to the category $\Gr _{O_0}$. 
 Luckily, when proving 
this compatibility it was possible to
use a technical result about sections of Dieudonne
crystalline sheaves from [Br1] to avoid diving into crystalline aspects of 
Breuil's theory. 
We also establish a criterion for $\Bbb F_p[\Gamma _{K_0}]$-modules to
appear in the form $G_0(\bar K)$, $G_0\in\Gr _{O_0}$, and apply it to relate these 
modules to Galois modules of kernels of isogenies of Drinfeld
modules via the field-of-norms functor, cf. Subsection 0.3 below 
for more commentaries.  Finally, we 
establish the interpretation of the Cartier duality in $\Gr _{O_0}$
in terms of filtered modules from $\MF ^e_S$. 
The Cartier duality is described in terms of special 
algebras for $H=H_0\otimes _{O_0}O$ and 
$\widetilde{H}=\widetilde{H}_0\otimes _{O_0}O$, where 
$H_0, \widetilde{H}_0\in\Gr _{O_0}$ and $\widetilde{H_0}$ is the Cartier dual 
to $H_0$, via an explicit construction of the corresponding non-degenerate 
bilinear pairing of group functors $H\times \widetilde{H}\longrightarrow\mu _{p,O}$. 

\medskip 

\subsubhead{\rm 0.3} Relation to Faltings's strict modules  
\endsubsubhead

Let $V$ be a finite $\Bbb F_p[\Gamma _K]$-module, where 
$\Gamma _K=\Gal (\bar K/K)$. Introduce the object 
$T(V)=(T(V)^0,T(V)^1,\varphi _1)$ of the category $\FM _S$ such that 
$T(V)^0=\Hom ^{\Gamma _K}(V,\bar O\operatorname{mod}p)$, 
$T(V)^1=\{a\in T(V)^0\ |\ a^p=0\}$ and $\varphi
_1:T(V)^1\longrightarrow T(V)^0$ is induced by the correspondences
$o\mapsto -o^p/p$ where $o\in \bar O$. 

In Section 8.2  we prove the following criterion: 

\proclaim{(0.3.1)} Suppose $\Cal M\in\MF ^e_S$, $G=\Cal G_O(\Cal
M)$ and $|G(\bar K)|=|V|$. Then the $\Bbb F_p[\Gamma _K]$-modules 
$V$ and $G(\bar K)$ are isomorphic if and only if in the
category $\FM _S$ there is a $\varphi _1$-nilpotent morphism 
$\Cal M\longrightarrow T(V)$.
\endproclaim 

This criterion makes precise the fundamental role of functions satisfying 
condition (0.1.1). It also allows us to study the following problem. Remind that
by Raynaud's theorem any finite flat commutative group scheme 
over $O_0$ 
arises as the kernel of an isogeny in the category
$\operatorname{Ab}_{O_0}$ of abelian schemes over $O_0$. The
characteristic $p$ analogue of $\operatorname{Ab}_{O_0}$ is the
category of Drinfeld modules $\operatorname{Dr}(S_{00})_{S_0}$ 
over $S_0$, where $S_{00}=\Bbb F_p[\tau _{00}]\subset S_0$ is such
that 
$\Cal K_{00}=\Frac \Bbb F_p[[\tau _{00}]]$ is a closed subfield in 
$\Cal K_0=\Frac S_0$ and the ramification index of $\Cal K_0$ over
$\Cal K_{00}$ is $e$. The kernels of isogenies in 
$\operatorname{Dr}(S_{00})_{S_0}$ are analogs of classical group
schemes. They can be introduced and studied directly as 
finite flat group 
schemes with strict action of $S_{00}$, cf. [Fa], [Ab4]. 
In this setting the characteristic $p$ analog of $\Gr _{O_0}$ is
the category $\Gr (S_{00})_{S_0}$ of finite flat commutative group
schemes over $S_0$ with strict action of $S_{00}$ 
which are killed by the action
of $\tau _{00}$. (In [Ab4] this category was denoted by 
$\operatorname{DGr}_1^{\prime *}(S_{00})_{S_0}$). 

As a special case of the classification of strict modules from [Ab4] we 
have the antiequivalence $\Cal G^S_{S_0}:
\MF ^e_S\longrightarrow \Gr (S_{00})_{S_0}$. (The category 
$\MF ^e_S$ was denoted in [Ab4] by $\operatorname{BR}_1(S_{00})_{S_0}$.) 

Suppose $\Cal M\in\MF ^e_S$, $H_0=\Cal G_{O_0}^O(\Cal M)\in\Gr _{O_0}$ and 
$\Cal H_0=\Cal G^S_{S_0}(\Cal M)\in\Gr (S_{00})_{S_0}$. 

Let $\bar\Cal K$ be an algebraic closure of $\Cal K_0$ and 
$\Gamma _{\Cal K_0}=\Aut _{\Cal K_0} (\bar\Cal K)$. As 
earlier, $\bar K$ is an algebraic closure of $K$ and $\Gamma _{K_0}=
\Gal (\bar K/K_0)$. Consider the 
$\Gamma _{K_0}$-module $V_0=H_0(\bar K)$ and the 
$\Gamma _{\Cal K_0}$-module $\Cal V_0=\Cal H_0(\bar\Cal K)$. 
(Notice that by [Ab4], $\Cal H_0$ has an etale generic fibre.) 

As an application of the above criterion (0.3.1) we show that the Galois modules 
$V_0$ and $\Cal V_0$ can be identified via the Fontaine-Wintenberger
functor field-of-norms. More precisely, consider the arithmetically
profinite extension 
$$K_{\infty }=K_0(\{\pi _n\ |\ \pi _1=\pi, \pi
_{n+1}^p=\pi _n\}).$$
 Then the field-of-norms functor gives an 
identification of $\Gamma _{K_{\infty }}=\Gal (\bar K/K_{\infty })$ 
and $\Gamma _{\Cal K_0}$ and we have the following property:

\proclaim{(0.3.2)} With the above identification $\Gamma
_{\Cal K_0}=\Gamma _{K_{\infty }}$, it holds $\Cal V_0\simeq V_0|_{\Gamma
_{K_{\infty }}}$. 
\endproclaim 

As a matter of fact, we can say more. Suppose 
$\Gamma _{K_0}(V_0)=\{\tau\in\Gamma _{K_0}\ |\ \tau |_{V_0}=\id \}$ and 
$\Gamma _{\Cal K_0}(\Cal V_0)=
\{\tau\in\Gamma _{\Cal K_0}\ |\ \tau |_{\Cal V_0}=\id \}$. Then the embedding 
$\Gamma _{\Cal K_0}=\Gamma _{K_{\infty }}\longrightarrow \Gamma _{K_0}$
induces a group isomorphism 
$\Gamma _{\Cal K_0}/\Gamma _{\Cal K_0}(\Cal V_0)\simeq\Gamma
_{K_0}/\Gamma _{K_0}(V_0)$. Therefore, the Galois modules 
$\Cal V_0$ and $V_0$ can be uniquely recovered one from another. 
There are another situations where there is definitely similar
relation between the kernels of isogenies of Drinfeld modules and 
the kernels of isogenies of abelian schemes. The study of this problem
should be useful when studying  
the image of the functor $V\mapsto V|_{\Gamma _{\Cal K_0}}$, 
where $V$ is a ``geometrically interesting'' (e.g. crystalline,
semistable) representation 
of $\Gamma _{K_0}$.

The above approach can be applied to the study of the functor from the
category of finite flat $\Bbb Z_p[\Gamma _{K_0}]$-modules 
(i.e. the Galois modules of the form $G(\bar K)$, where $G\in\Gr
'_{O_0}$) to the category of $\Gamma _{K_{\infty }}$-modules 
given by the restriction of Galois action to $\Gamma _{K_{\infty
}}\subset\Gamma _{K_0}$. 
Our method allows to obtain Breuil's result, [Br4, Theorem 3.4.3], 
about full faithfulness of this functor just from Fontaine's 
ramification estimate:  for all $v>ep/(p-1)-1$, 
the ramification subgroups $\Gamma
_{K_0}^{(v)}$ act trivially on $G_0(\bar K)$, where $G_0\in\Gr
_{O_0}$. This idea also works 
in the context of finite subquotients of     
crystalline representations 
over unramified base with Hodge-Tate weights of length $<p$ 
by using the ramification estimate from [Ab5]. (The case of 
Hodge-Tate weights $\leqslant p-2$ was considered in [Br3] and can be
retrieved even with Fontaine's ramification estimate from [Fo5].) 
We must mention here that recently Kisin [Ki3, Theorem 02] 
proved the full faithfulness of the functor 
$V\mapsto V|_{\Gamma _{K_{\infty }}}$ in the context of all crystalline 
$\Bbb Q_p[\Gamma _{K_0}]$-modules. (Everywhere in this paragraph 
$p$ is any prime number.)

\ \

We shalll use without special reference the following notation:

{\it Basic Notation.} 
\medskip 

$\bullet $ \ \ If $A,B,C$ are sets and $f:A\longrightarrow B$, $g:B\longrightarrow C$
then their composition will be denoted by $fg$ or, sometimes, by
$f\circ g$. 
\medskip 

$\bullet $\ \ $k$ is a perfect field of characteristic $p>2$, $K_{00}$ is the
fraction field of the ring of Witt vectors $W(k)$, $K_0$ is a totally
ramified extension of $K_{00}$ of degree $e$ with fixed
uniformising element $\pi _0$ and the valuation ring $O_0$, 
$K=K_0(\pi )$ and $O=O_0[\pi ]$, where $\pi ^p=\pi _0$; 
$\bar K$ is a fixed algebraic closure of $K$, $\bar O$ is the
valuation ring of $\bar K$, and for any 
field extension $E$ of $K_0$ in $\bar K$, $\Gamma
_E=\operatorname{Gal}(\bar K/E)$;
\medskip 

$\bullet $\ \ $\Gr '_{O_0}$, resp., $\Gr '_{O}$, is the category 
of all finite flat $p$-group
schemes over $O_0$, resp., $O$; $\Gr _{O_0}$ and $\Gr _O$ are the
corresponding full subcategories consisting of objects killed by $p$; for any
finite flat group scheme $H$ we denote by $A(H)$, $\Delta _H$ and
$e_H$, resp., the affine algebra, the comultiplicatiopn and the counut
of $H$;
\medskip 

$\bullet $ \ \ In the category $\Aug _O$ of augmented $O$-algebras,
$I_A$ is always the augmentation ideal of $A\in\Aug _O$; 
\medskip 

$\bullet $ \ \ $\phi (X,Y)=(X^p+Y^p-(X+Y)^p)/p\in\Bbb Z[X,Y]$ is the
first Witt polynomial and $\phi (X)$ is just an abbreviation for 
$\phi (X\otimes 1, 1\otimes X)$;
\medskip 

$\bullet $ \ \ for an indeterminate $t$, we set $S=k[[t]]$, $\sigma
:S\longrightarrow S$ is the morphism of $p$-th power,
$S_0=k[[t_0]]\subset S$ with $t_0=t^p$, 
the corresponding fraction fields are $\Cal K=\Frac S$ 
and $\Cal K_0=\Frac S_0$; $\kappa _{SO}:S/t^{ep}\longrightarrow O/pO$ is
a ring isomorphism such that $\kappa _{SO}(t)=\pi\operatorname{mod}p$
and for any $\alpha\in k$, $\kappa _{SO}(\alpha )=[\alpha
]\operatorname{mod}p$, where $[\alpha ]$ is the Teichm{\"u}ller 
representative of $\alpha $. 
\medskip 

{\it Acnowledgements.} The author expresses deep gratitude to the
referee for very interesting and substantial remarks and advices. 
\medskip 
\medskip

\subhead 1. A category of filtered modules 
\endsubhead 
\medskip

\medskip

\subsubhead{\rm 1.1.} The categories $\MF _S^e$ and $\FM _S$ 
\endsubsubhead

Suppose $e\in\Bbb N$. 
Denote by $\MF _{S}^e$ the category of triples $(M^0,M^1,\varphi _1)$,
where $M^0$ is a free $S$-module of finite rank, $M^1\subset M^0$ is
a submodule such that $M^1\supset t^eM^0$, and 
$\varphi _1:M^1\longrightarrow M^0$ is a $\sigma $-linear morphism 
of $S$-modules such that the set $\varphi _1(M^1)$ generates $M^0$
over $S$. Notice that $M^1$ is a free $S$ module, $\rk
_{S}M^0=\rk _{S}M^1$, and $\varphi _1$ maps any
$S$-basis of $M^1$ to an $S$-basis of $M^0$. 
The morphisms $(M^0,M^1,\varphi _1)\longrightarrow (M_1^0,M_1^1,\varphi
_1)$ in $\MF _{S}^e$ are given by $S$-linear morphisms $f:M^0\longrightarrow M_1^0$
such that $f(M^1)\subset M_1^1$ and $f\varphi _1=\varphi _1f$. 

Let $S_0=\sigma (S)$. Then $S_0=k[[t_0]]$, where $t_0=t^p$. 
Consider the category 
$\MF _{e,S_0}$ of $S_0$-modules with slope
$\leqslant e$. Its objects are couples $(M,\varphi )$, where $M$ is a
free $S_0$-module of finite rank, 
$\varphi :M^{(\sigma )}\longrightarrow M$, where 
$M^{(\sigma )}=M\otimes _{(S_0,\sigma )}S_0\longrightarrow M$, is 
an $S_0$-linear morphism such that its image contains $t_0^eM$. 
The morphisms in $\MF _{e,S_0}$ are  
morphisms of the corresponding modules which commute with $\varphi $. 
The category $\MF _{e,S_0}$ is equivalent to $\MF ^e_{S}$. This equivalence
can be given by the identification of rings 
$S_0\otimes _{(S_0,\sigma )}S_0=S$ (where $t_0\otimes 1\mapsto t$), 
the correspondence $(M,\varphi )\mapsto (M^{(\sigma )},
M^1,\varphi _1)$, where $M^1=\varphi ^{-1}(t_0^eM)$ and $\varphi
_1=t_0^{-e}\eta \varphi |_{M^1}$ with any fixed $\eta\in S^*$.  
\medskip

Introduce the category $\FM _{S}$ as the category of triples
$(M^0,M^1,\varphi _1)$, where $M^0$ is an 
$S$-module, $M^1$ is a 
submodule in $M^0$, and $\varphi _1$ is a $\sigma $-linear morphism  
from $M^1$ to $M^0$. The morphisms $f:(N^0,N^1,\varphi _1)\longrightarrow 
(M^0,M^1,\varphi _1)$ in $\FM _S$ are given by linear morphisms $f:N^0\longrightarrow
M^0$, such that $f(N^1)\subset M^1$ and $f\varphi _1=\varphi _1f$. 

The category $\FM _S$ is additive, in particular, 
$f$ is epimorphic iff $f(N^0)=M^0$. We shall call $f$ {\it strictly epimorphic} if
in addition $f(N^1)=M^1$. 

Notice that $\MF _S^e$ is a full subcategory in $\FM _S$.
\medskip

\subsubhead{\rm 1.2} $\varphi _1$-nilpotent lifts 
\endsubsubhead  
 
Suppose $\Cal N=(N^0,N^1, \varphi _1)$ and 
$\Cal M=(M^0,M^1, \varphi _1)$ are objects of $\FM _S$ and 
$\theta\in \Hom _{\FM _S}(\Cal N, \Cal M)$.

\definition{Definition} a) $\theta $ will be called 
 $\varphi _1${\it -nilpotent} if  
for $T=\ker \theta\subset N^0$, it holds 
$\Ker (\theta |_{N^1})=T$, $\varphi _1(T)\subset T$ and 
$\varphi _1|_T$ is topologically nilpotent. (This means that if
$T=T^{(0)}$ and for $i\geqslant 0$, $T^{(i+1)}=\varphi
_1(T^{(i)})S$, then $\dsize\bigcap \Sb i\geqslant 0\endSb
T^{(i)}=0$.) 
\newline 
b) if $\Cal N\in\MF _S^e$ and $\theta $ is strictly epimorphic  
and $\varphi _1$-nilpotent  then we say 
that $\theta $ is {\it a $\varphi _1$-nilpotent lift} of $\Cal M$ to $\MF _S^e$. 
\enddefinition 

Notice that the composition of $\varphi _1$-nilpotent lifts is again
$\varphi _1$-nilpotent.

\proclaim{Proposition 1.2.1} Suppose $\Cal M,\Cal M_1\in\FM _S$ and $\theta
_1:\Cal N_1\longrightarrow \Cal M_1$, $\theta :\Cal N\longrightarrow
\Cal M$ are $\varphi _1$-nilpotent lifts to $\MF _S^e$. Then for any 
$f\in\Hom _{\FM _S}(\Cal M_1,\Cal M)$ there is a unique 
$\hat f\in\Hom _{\MF _S^e}(\Cal N_1,\Cal N)$ such that 
$\hat f\theta =\theta _1f$. 
\endproclaim 

{\demo{Proof} The proof can be easily obtained from the following Lemma. 
\enddemo

\proclaim{Lemma 1.2.2}  
Suppose $\Cal N,\Cal M\in\FM _S$ and  
$\theta \in\Hom _{\FM _S}(\Cal N,\Cal M)$ is 
$\varphi _1$-nilpotent and strictly epimorphic. 
Then for any $\Cal N_1\in\MF _S^e$, the map 
$\theta _*:\Hom _{\FM _S}(\Cal N_1,\Cal N)\longrightarrow\Hom _{\Cal\FM
_S}(\Cal N_1,\Cal M)$ is bijective. 
\endproclaim

\demo{Proof of Lemma} Let $\Cal N_1=(N_1^0,N_1^1,\varphi _1)$. Choose an $S$-basis 
$n_1^1,\dots ,n_u^1$ in $N_1^1$ and set $\bar n_1^1=(n_1^1,\dots ,n_u^1)$. 
Let $\bar n_1=\varphi _1(\bar n_1^1):=(n_1,\dots ,n_u)$, where for all
$1\leqslant i\leqslant u$, 
$n_i=\varphi _1(n_i^1)$. Then $n_1,\dots ,n_u$ is an $S$-basis of
$N_1^0$ and 
there is a matrix $U\in M_u(S)$ such that $\bar n_1^1=\bar n_1U$.
We shall use below a similar vector notation. 

Suppose $\Cal N=(N^0,N^1,\varphi _1)$, $\Cal M=(M^0,M^1,\varphi
_1)$ and $f\in\Hom _{\FM _S}(\Cal N_1,\Cal M)$. Then $f$ is 
uniquely given by two vectors $f(\bar n_1)=\bar m$ and 
$f(\bar n_1^1)=\bar m^1$ with coordinates in $M^0$ and  
$M^1$, resp.,  
such that $\varphi _1(\bar m^1)=\bar m$
and $\bar m^1=\bar mU$. 

Choose a vector $\bar n^{(0)1}$ with coordinates in $N^1$ such
that $\theta (\bar n^{(0)1})=\bar m^1$ and set 
$\bar n^{(0)}=\varphi _1(\bar n^{(0)1})$. For $i\geqslant 0$,
define by induction on $i$ the vectors $\bar n^{(i)}$ and
$\bar n^{(i)1}$ as follows: $\bar n^{(i+1)1}=\bar
n^{(i)}U$ and $\varphi _1(\bar n^{(i+1)1})=\bar n^{(i+1)}$. 
\medskip 

Then the sequences $\{\bar n^{(i)1}\}_{i\geqslant 0}$
and $\{\bar n^{(i)}\}_{i\geqslant 0}$ converge in $N^1$ 
and, resp., $N^0$. 
\medskip 

 Indeed, use the submodules $T^{(i)}$, $i\geqslant 0$, in
$N^1$ from the above definition of $\varphi _1$-nilpotent morphism 
$\pi $. 
Then $\bar n^{(0)1}-\bar n^{(1)1}$ has coordinates in $T=T^{(0)}$, 
its $\varphi _1$-image $\bar
n^{(0)}-\bar n^{(1)}$ has coordinates in 
$\varphi _1(T^{(0)})$ and $\bar n^{(1)1}-\bar n^{(2)1}$ 
has coordinates in 
$\varphi _1(T^{(0)})S_1=T^{(1)}$. Similarly, for any $i\geqslant 0$, $\bar
n^{(i)}-\bar n^{(i+1)}$ and $\bar n^{(i)1}-\bar n^{(i+1)1}$ have
coordinates 
in $T^{(i)}$. So, the condition 
$\dsize\bigcap\Sb i\geqslant 0\endSb T^{(i)}=0$ implies that the both
sequences are Cauchy and, therefore, converge. 
\medskip

Now let $\bar n^1\in N^1$ and $\bar n\in N^0$ be limits of the 
above sequences  $\{\bar n^{(i)1}\}_{i\geqslant 0}$
and $\{\bar n^{(i)}\}_{i\geqslant 0}$, resp. 
Then $\varphi _1(\bar n^1)=\bar n$ and $\bar n^1=\bar nU$. 
So, the  correspondences 
$\bar n_1^1\mapsto \bar n^1$ and $\bar n_1\mapsto \bar n$ define $g\in\Hom
_{\FM _S}(\Cal N_1,\Cal N)$. But $\theta (\bar n)=\bar m$ and $\pi (\bar
n^1)=\bar m^1$ because for all $i\geqslant 0$, $\pi (\bar
n^{(i)})=\bar m$ and $\pi (\bar n^{(i)1})=\bar m^1$. So, 
$f=\theta _*(g)$ and $\theta _*$ is surjective.

Suppose $g'\in\Hom _{\FM _S}(\Cal N_1,\Cal N)$ is such that $\theta
_*(g')=f$. 
Then $g'(\bar n_1)=\bar n'$ and, resp., $g'(\bar n_1^1)=
\bar n^{\prime 1}$ have
coordinates in $N^0$ and $N^1$, resp.  
Notice also that $\varphi _1(\bar n^{\prime 1})=\bar n'$, 
$\bar n^{\prime 1}=\bar n'U$, $\theta
(\bar n')=\bar m$ and $\theta (\bar n^{\prime 1})=\bar m^1$. 
Therefore, 
$\bar
n^{\prime 1}-\bar n^{(0)1}$ has coordinates in $T^{(0)}$, $\bar n'-\bar
n^{(0)}$ has coordinates in $\varphi _1(T^{(0)})$, $\bar n^{\prime 1}-\bar
n^{(1)1}$ has coordinates in $\varphi _1(T^{(0)})S_1=T^{(1)}$
and so on. In other words, for any $i\geqslant 0$, $\bar n^{\prime 1}-\bar
n^{(i+1)1}$ and $\bar n'-\bar n^{(i)}$ have coordinates in
$T^{(i)}$. Taking limits we obtain that $\bar n^{\prime 1}=\bar n^1$ and $\bar
n'=\bar n$, i.e. $g=g'$. 

The Lemma is proved. $\square $
\enddemo 

1.2.3. With notation from above Proposition 1.2.1, 
$\hat f$ will be called {\it a $\varphi _1$-nilpotent lift} of $f$. 
Notice also the following properties: 

a) the correspondence $f\mapsto \hat f$ induces an injective
homomorphism 
$$\Hom _{\FM _S}(\Cal M_1,\Cal M)\longrightarrow \Hom
_{\MF ^e_S}(\Cal N_1,\Cal N);$$

b) a $\varphi _1$-nilpotent lift is unique up to a unique isomorphism
in $\MF ^e_S$; in particular, a $\varphi _1$-nilpotent lift of 
$\Cal M\in\MF ^e_S$ to $\MF ^e_S$ is an isomorphism; 

c) if $\Cal M\in\FM _S$, $\Cal N=(N^0,N^1,\varphi _1)\in\MF ^e_S$ 
and $\theta :\Cal N\longrightarrow \Cal M$ is a $\varphi _1$-nilpotent
lift of $\Cal M$ then $\Ker \theta\subset tN^1$. (Indeed, 
if $n\in\Ker\theta $ and $n\notin tN^1$ then for any $i\geqslant 1$, 
$\varphi _1^i(n)\in N^1\setminus tN^0$ and, therefore, does not
converge to 0.)  
\medskip

\subsubhead{\rm 1.3} Simplest objects and their extensions 
\endsubsubhead 

Let $\Cal M=(M^0,M^1,\varphi _1)$, $\Cal M_1=(M^0_1,M^1_1,\varphi
_1)$, 
$\Cal M_2=(M_2^0,M_2^1,\varphi _1)$ be objects of $\FM _S$. By
definition, the sequence 
$0\longrightarrow \Cal M_1\longrightarrow\Cal M
\longrightarrow\Cal M_2\longrightarrow 0$ 
is {\it short exact} in $\FM _S$ if the corresponding sequence 
of $S$-modules 
$0\longrightarrow M^0_1\longrightarrow M^0
\longrightarrow M^0_2\longrightarrow 0$ and the induced sequence of 
maps of their submodules 
$0\longrightarrow M^1_1\longrightarrow M^1
\longrightarrow M^1_2\longrightarrow 0$ 
are short exact. 
Then for given $\Cal M_1,\Cal M_2$, one can define, as usually, the
set of   
classes of equivalence of short exact sequences   
$\Ext _{\FM _S}(\Cal M_2, \Cal M_1)$, and this set has a natural structure of
abelian group.

\definition{Definition} If $\tilde s\in S$ is such that $\tilde s|t^e$ 
define the object $\Cal M_{\tilde s}$ of $\MF _S^e$ as 
$(Sm,Sm^1,\varphi _1)$ such that 
$\varphi _1(m^1)=m$ and $m^1=\tilde sm$.   
Such objects $\Cal M_{\tilde s}$ will be called {\it simplest}. 
\enddefinition 

\remark{Remark {\rm 1.3.1}}
Notice that, if $\tilde s,\tilde s'\in S$ are divisors of $t^e$
then $\Cal M_{\tilde s}\simeq\Cal M_{\tilde s'}$
iff $\tilde s'=\tilde su^{p-1}$, where $u\in S^*$. In particular, 
by enlarging if necessary the residue field $k$ we can 
always assume that $\tilde s$ is just an integral power of $t$.
\endremark 
\medskip

Let  
$\Cal K^{\prime }$ be a finite field extension of
$\Cal K=\Frac S$ and let $S^{\prime }$ be the valuation ring of $\Cal
K^{\prime }$. 
If 
$e_0$ is the ramification index of the field extension 
$\Cal K^{\prime }/\Cal K$ and $e'=ee_0$ 
then there is a functor from $\MF ^e_{S}$ to $\MF _{S'}^{e'}$ 
given by the extension of scalars $\Cal M=(M^0,M^1,\varphi _1)
\mapsto\Cal M\otimes _{S}S':=(M^0\otimes _{S}S', 
M^1\otimes _{S}S',\varphi
_1\otimes _{S}S')$.  

\proclaim{Proposition 1.3.2} If 
$\Cal M\in\MF _{S}^e$ then there is a tamely ramified extension 
$\Cal K^{\prime }/\Cal K$ such that $\Cal M\otimes _{S}S'$ can be obtained 
by a sequence of successive extensions via simplest objects of the
category $\MF _{S'}^{e'}$. 
\endproclaim 

\demo{Proof }
Let $\Cal M=(M^0,M^1,\varphi _1)$.  
The embedding $M^1\subset M^0$ induces the identification of 
$\Cal K$-vector spaces 
$V:=M^0\otimes _{S}\Cal K=M^1\otimes _{S}\Cal K$ and $\varphi _1$ induces a $\sigma $-linear
morphism $\varphi _1:V\longrightarrow V$ such that $\varphi _1(V)\Cal
K=V$. Therefore, $V$ is an etale $\varphi _1$-module, cf. [Fo4], 
and we can apply the antiequivalence of the category of etale $\varphi
_1$-modules and the category of continuous $\Cal K[\Gamma _{\Cal
K}]$-modules $H$, where $\Gamma _{\Cal K}=\Gal (\Cal K_{\sep }/\Cal
K)$, cf. [loc cit]. This antiequivalence is given by the
correspondence 
$$V\mapsto H:=\{f\in\Hom _{\Cal K}(V,\Cal K_{\sep })\ |\ 
\forall v\in V, f(\varphi _1(v))=f(v)^p\}.$$
(Here $\Gamma _{\Cal K}$ acts on $H$ via its natural action on $\Cal
K_{\sep }$.) Notice that the inverse functor is induced by the
correspondence 
$H\mapsto V=\Hom _{\Cal K[\Gamma _{\Cal K}]}(H,\Cal K_{\sep
})$. 

Then use that the action of the wild inertia subgroup of $\Gamma _{\Cal
K}$ on $H$ is unipotent. This implies the existence of a finite
tamely ramified field extension $\Cal K^{\prime }$ of $\Cal K$ such that  
$H\otimes _{\Cal K}\Cal K'$ has a
decreasing filtration by its $\Cal K'[\Gamma _{\Cal
K'}]$-submodules such that the corresponding quotients 
are 1-dimensional $\Cal K'$-vector spaces with the trivial action of
$\Gamma _{\Cal K'}$. Therefore, $V\otimes _{\Cal K}\Cal K'$ has a
$\Cal K'$-basis $v_1,\dots ,v_u$ such that for all $1\leqslant
i\leqslant u$, 
$$\varphi _1(v_i)=v_i+\sum\Sb j>i\endSb v_j\alpha _{ji},\ \alpha
_{ji}\in\Cal K'.$$
Therefore, $M^1\otimes _{S}S'$ has an $S'$-basis $m_1^1,\dots
,m_u^1$ such that for $1\leqslant i\leqslant u$,
$$\varphi _1(m_i^1)=\sum\Sb j\geqslant i\endSb m^1_ju_{ji},\ u_{ji}\in
S'.$$
It remains to notice that for $1\leqslant i\leqslant u$, 
$m_i=\varphi _1(m_i^1)$ is an $S'$-basis of 
$M^0_{S'}=M^0\otimes _{S}S'$, and the
condition $M^1\supset t^eM^0$ implies that 
all $\tilde s_i:=u_{ii}$ divide $t^e$. The proposition is proved.
$\square $
\enddemo

\proclaim{Proposition 1.3.3} 
Suppose $\tilde s\in S$, $ \tilde s|t^e$  
and $\Cal N=(N^0,N^1,\varphi _1)\in\MF _S^e$. Then there is
a natural isomorphism of the group 
$\Ext _{\MF _S^e}(\Cal M_{\tilde s},\Cal N)$ onto 
$Z_{\tilde s}(\Cal N)/B_{\tilde s}(\Cal N)$, where 
$Z_{\tilde s}(\Cal N)=\{n\in N^0\ |\ t^e\tilde s^{-1}n\in N^1\}$ and  
$B_{\tilde s}(\Cal N)=\{v^1-\tilde s\varphi _1(v^1)\ |\ v^1\in N^1\}$.  
\endproclaim 

\demo{Proof} By definition,    
$\Cal M_{\tilde s}=(Sm,Sm^1,\varphi _1)$, where 
$\varphi _1(m^1)=m$ and $m^1=\tilde sm$. 
Suppose $\Cal M=(M^0,M^1,\varphi _1)
\in\Ext _{\MF _S^e}(\Cal M_{\tilde s},\Cal N)$. 
Then $\Cal M$ can be described as follows: 
$M=N^0\oplus S\hat m$,  
$M^1=N^1+S\hat m^1$, where 
$\hat m^1=\tilde s\hat m+n(\Cal M)$ with $n(\Cal M)\in N^0$ 
and $\varphi _1(\hat m^1)=\hat m$. 
Notice that $M^1\supset t^eM^0$ holds if and only if 
$n(\Cal M)\in Z_{\tilde s}(\Cal N)$ and the
morphism $\varphi _1$ is uniquely defined. 
Any equivalent to $\Cal M$ extension
$\Cal M'$ 
can be decribed by another lifts $\hat m'=\hat m+v$, 
$\hat m^{\prime 1}=\hat m^1+v^1$ with $v\in N^0$ and $v^1\in N^1$ such that 
$\varphi _1(v^1)=v$. Then the corresponding element 
$n(\Cal M')$ equals $n(\Cal M)+v^1-\tilde s\varphi _1(v^1)$, 
i.e. $n(\Cal M)\equiv n(\Cal M')
\operatorname{mod}B_{\tilde s}(\Cal N)$.  Finally, a straighforward 
verification shows that the correspondence 
$\Cal M\mapsto n(\Cal M)\operatorname{mod}B_{\tilde s}(\Cal N)$ 
gives the required isomorphism. 
$\square $
\enddemo 

\remark{Remark {\rm 1.3.4}} If $\tilde s\in S^*$ then we can always 
choose $n(\Cal M)\in N^1+tN^0$ 
(use that $\varphi _1(N^1)$ generates $N^0$).  
\endremark 
\medskip 

Let $S'=S[t']$ where $t^{\prime p}=t$. Consider the extension of
scalars 
\linebreak 
$\Cal M\mapsto\Cal M\otimes _{S}S'\in\MF _{S'}^{ep}$, 
where $\Cal M\in\MF
^e_{S}$. Consider the induced group homomorphism 
$\pi _{SS'}:\Ext _{\MF _{S}^e}(\Cal M_{\tilde s},\Cal N)
\longrightarrow\Ext _{\MF _{S'}^{ep}}(\Cal M_{\tilde s}\otimes _{S}S',\Cal
N\otimes _{S}S')$.

Choose a basis $n_1^1,\dots ,n_u^1$ of $N^1$ such that for $1\leqslant
i\leqslant u$, there are $\tilde s_i\in S$ such that the elements 
$\tilde s_i^{-1}n_i^1=n_i$ form a basis of $N^0$. With this notation 
$$Z_{\tilde s}(\Cal N)=\left\{\sum\Sb i\endSb \alpha _in_i\ |\ \text{
all }\alpha _i\in S 
\text{ and }t^e\tilde s^{-1}\alpha _i\equiv 0\operatorname{mod}\tilde s_i\right\}.$$
The module $Z_{\tilde s}(\Cal N\otimes _{S}S')$ is given similarly with the
only difference that  
all coefficients $\alpha _i$ should belong to $S'$. 

\proclaim{Proposition 1.3.5} With the above notation suppose 
$z=\sum\Sb i\endSb \alpha _in_i\in Z_{\tilde s}(\Cal N\otimes
_{S}S')$. Then $z\operatorname{mod}B_{\tilde s}(\Cal
N\otimes _{S}S')$ 
belongs to the image of $\pi _{SS'}$ if and only if 
for all $i$, $\alpha _i
\in S\operatorname{mod}\tilde s_i$. 
\endproclaim 

\demo{Proof} Suppose $z\operatorname{mod}B_{\tilde s}(\Cal N\otimes
_SS')
=\pi _{SS'}(y)$, where $y=\sum\Sb i\endSb
\beta _in_i\operatorname{mod}B_{\tilde s}(\Cal N)$ with 
all $\beta _i\in S$. This means that 
$\sum\Sb i\endSb \alpha _in_i=\sum\Sb i\endSb \beta _in_i+v^1-\tilde s\varphi
_1(v^1)$, 
where $v^1=\sum\Sb i\endSb \gamma _in_i^1=\sum\Sb i\endSb \gamma
_i\tilde s_in_i\in N^1\otimes _{S}S'$. 
Then $\varphi _1(v^1)=\sum\Sb i\endSb \gamma _i^p\varphi _1(n^1_i)
=\sum\Sb i\endSb \delta _in_i\in N^0\subset N^0\otimes _SS'$, 
all $\alpha _i=\beta _i+\gamma _i\tilde s_i-\tilde s\delta _i$ 
and $\alpha _i\in (\beta _i-\tilde s\delta
_i)\operatorname{mod}\tilde s_i\in
S\operatorname{mod}\tilde s_i$. 

Conversely, suppose for all $i$, $\alpha _i=\alpha _i^0+\tilde
s_i\alpha _i'$, where $\alpha _i^0\in S$ and $\alpha _i'\in S'$. 
Then $\sum\Sb i\endSb \tilde s_i\alpha _i'n_i\in N^1\otimes _SS'$, 
$$\sum\Sb i\endSb \alpha _in_i\equiv \sum\Sb i\endSb \alpha _i^0n_i
+\tilde s\varphi _1\left (\sum\Sb i\endSb \tilde s_i\alpha
_i'n_i\right )\operatorname{mod}B_{\tilde s}(N\otimes _SS')$$
and the right-hand side is defined over $S$.The proposition is
proved. 
$\square $
\enddemo 
\medskip 

\subsubhead{\rm 1.4.} Special bases 
\endsubsubhead 

Let $\Cal M=(M^0,M^1,\varphi _1)\in\MF _S^e$.

\definition{Definition} An $S$-basis $m_1^1,\dots ,m_u^1$ of $M^1$ 
will be called {\it special} if the non-zero images of $m^1_i$, $1\leqslant
i\leqslant u$, in $M^0\operatorname{mod}tM^0$ are linearly 
independent over $k$. 
\enddefinition 

Suppose $m_1^1,\dots ,m_u^1$ is a special basis of $M^1$, 
 $m_i=\varphi _1(m_i^1)$ if $1\leqslant i\leqslant u$,  
and $U\in M_u(S_1)$ is such that 
$(m_1^1,\dots ,m_u^1)=(m_1,\dots ,m_u)U$. 
Notice that the condition $M^1\supset t^eM^0$ implies that 
$U$ divides the scalar matrix $t^eE$ in $M_u(S)$ 
(where $E$ is the unit matrix of order $u$), i.e. 
there is an $V\in M_u(S)$ such that $UV=t^eE$. Let 
$U=(u_{ij})$, $V=(v_{ij})$ where all entries $u_{ij},v_{ij}\in S$. 

\proclaim{Proposition 1.4.1} With the above notation, 
if $1\leqslant i,j,r\leqslant u$, then   
$$u_{ij}v_{jr}\equiv 0\operatorname{mod}t.$$ 
\endproclaim 

\demo{Proof} Because the basis $m_1^1,\dots ,m_u^1$ is special 
we can assume that there is an index $i_0$, such that 
$m_1^1,\dots ,m_{i_0}^1\in tM^0$ and $m^1_{i_0+1},\dots ,m^1_u$
are linearly independent modulo $tM^0$. Consider the image of the equality 
$\dsize t^em_r=\sum \Sb j\endSb m_j^1v_{jr}$ in
$M^0\operatorname{mod}t$, where $1\leqslant r\leqslant u$.  This gives 
$\dsize 0=\sum\Sb i_0<j\leqslant u\endSb
(m_j^1\operatorname{mod}t)v_{jr}$ and, therefore, 
$v_{jr}\in tS$ if
$i_0<j\leqslant u$. On the other hand, 
if $1\leqslant j\leqslant i_0$ then 
$\dsize m_j^1=\sum\Sb i\endSb m_iu_{ij}\in tM^0$ and 
for any $1\leqslant i\leqslant u$,  $u_{ij}\in
tS$. The proposition is proved.
$\square $ 
\enddemo 
\medskip

\subhead 2. Construction of the functor 
$\Cal G _O:\MF ^e_{S}\longrightarrow\Gr _{O}$ 
\endsubhead 
\medskip   

\subsubhead {\rm 2.1.} The category $\Aug _{O}$ and the 
functor $\iota :\Aug _{O}\longrightarrow\FM _S$
\endsubsubhead 

The objects of the category $\Aug _{O}$ 
are flat $O$-algebras $A$ 
of finite rank over $O$ with a given augmentation ideal $I_A$. 
The morphisms are morphisms of augmented algebras. 

\definition{Definition} If $A\in \Aug _{O}$ then: 

a) $I_A(p):=\{a\in I_A\ |\ a^p\in pA\}$;

b) $I_A^{DP}$ is the maximal ideal of $A$ with nilpotent
divided powers or, equivalently, such that if $a_1=a\in I^{DP}_A$ 
and for any $i\in\Bbb N$, 
$a_{i+1}=a^p_{i}/p$, then all $a_i\in I_A^{DP}$ and $\dsize \lim\Sb
i\to\infty \endSb a_i=0$. 
\enddefinition

Notice that $I_A/I_A^{DP}$ is killed
by $p$ (remind that $p>2$) and we can use the identification 
$\kappa _{SO}$  
to provide $I_A/I_A^{DP}$ with an
$S$-module structure. Then the triple 
$\iota ^{DP}(A):=(I_A/I^{DP}_A, I_A(p)/I^{DP}_A,\varphi _1)$, where $\varphi _1$ is
induced by the correspondence 
\linebreak 
$a\mapsto -a^p/p$ with $a\in I_A(p)$, 
is an object of the category 
$\FM _S$. The correspondence $A\mapsto\iota ^{DP}(A)$ gives rise to 
the functor $\iota ^{DP}$ from $\Aug _{O}$ to $\FM _S$. 

\proclaim{Proposition 2.1.1} Suppose $A\in\Aug _{O}$, $u\geqslant 1$, 
$b_1,\dots ,b_u\in I_A/I_A^{DP}$ and elements 
$b_1^1,\dots ,b^1_u\in I_A(p)/I_A^{DP}$ are such that for $1\leqslant
i\leqslant u$, 
$\varphi _1(b_i^1)=b_i$. Suppose $\hat U\in M_u(O)$ is such that 
$(b_1^1,\dots ,b_u^1)=(b_1,\dots ,b_u)\hat U$. Then for $1\leqslant
i\leqslant u$, there are unique  
$\hat b_i\in I_A$, $\hat b_i^1\in I_A(p)$ such that 
$\hat b_i\operatorname{mod} I_A^{DP}=b_i$, 
$\hat b_i^1\operatorname{mod}I_A^{DP}=b_i^1$, 
$(\varphi _1(\hat b_1^1), \dots ,\varphi _1(b_u^1))=(\hat b_1, \dots
,\hat b_u)$  and 
$(\hat b_1^1, \dots ,\hat b_u^1)=(\hat b_1, \dots ,\hat b_u)\hat U$. 
\endproclaim 

\demo{Proof} The proof is very similar to the proof of proposition 1. 

Use the vector notation, e.g. $\bar b=(b_1,\dots ,b_u)$, $\bar b^1=(b_1^1,\dots
,b_u^1)$. Choose $\bar b^{(0)1}$ with coordinates in $I_A(p)$
such that $\bar b^{(0)1}\operatorname{mod}I_A^{DP}=\bar b^1$. 
Then define for $i\geqslant 1$, $\bar b^{(i)}$ and 
$\bar b^{(i)1}$ via the relations 
$\bar b^{(i+1)}=\varphi _1(\bar b^{(i)1})$ and 
$\bar b^{(i)1}=b^{(i-1)}\hat U$. Consider the sequence 
of ideals $J_i$, $i\geqslant
0$, such that $J_0=I_A^{DP}$ and 
$J_{i+1}=I_A(p)J_i+\varphi _1(J_i)$, where $\varphi _1(J_i)$ is the
ideal generated by all elements $\varphi _1(a)$, $a\in J_i$. 
Notice that for all $i\geqslant 1$, 
$\bar b^{(i)1}\equiv \bar b^{(i-1)1}
\operatorname{mod}J_{i-1}$ and 
$\bar b^{(i)}\equiv \bar b^{(i-1)}\operatorname{mod}J_{i-1}$. This 
proves our proposition because 
$\dsize\bigcap \Sb i\geqslant 0\endSb J_i=0$. 
$\square $
\enddemo 
\medskip 

\subsubhead{\rm 2.2.} The family of augmented $O$-algebras $\Cal A(\Cal M)$,
$\Cal M\in\MF _S^e$ 
\endsubsubhead

Suppose $\Cal M=(M^0,M^1,\varphi _1)$ and the coordinates of the vector 
$\bar m^1=(m_1^1,\dots ,m_u^1)$ form a special basis in $M^1$. As
earlier, the coordinates $m_1,\dots ,m_u$ of 
$\varphi _1(\bar m^1)=\bar m$  
form an
$S$-basis of $M^0$ and there is an $U\in M_u(S)$ such that $\bar m^1=\bar
mU$. 

Choose $\hat U\in M_u(O)$ such that $\hat U\operatorname{mod}p=\kappa _{SO} 
(U\operatorname{mod}t^{ep})$. 
Introduce the augmented $O$-algebra $A$ as a quotient of $O[Y_1,\dots ,Y_u]$
by the ideal 
$$J_A:=J_{A,K}\bigcap O[Y_1, \dots ,Y_u],$$ 
where $J_{A,K}$ is the ideal in $K[Y_1,\dots ,Y_u]$ generated by
the coordinates $F_1,\dots ,F_u$ of the vector 
$\bar F=(\bar Y\hat U)^{(p)}+p\bar Y$. By definition the augmentation
ideal $I_A$ of $A$ is generated by $Y_1\operatorname{mod}{J_{A}},\dots
,Y_u\operatorname{mod}J_A$. 
Here and everywhere below we use the vector notation 
$\bar Y=(Y_1,\dots ,Y_u)$ and for any
matrix $C=(c_{ij})$, $C^{(p)}:=(c_{ij}^p)$. 
So, if $\hat U=(\hat
u_{ij})$ then for $1\leqslant i\leqslant u$, 
$F_i=(\sum\Sb j\endSb Y_j\hat u_{ji})^p+pY_i$. If there is no
risk of confusion we shall use just the notation $Y_1,\dots ,Y_u$ for the elements 
$Y_1\operatorname{mod}J_A,\dots ,Y_u\operatorname{mod}J_A$ of $A$.  

\proclaim{Proposition 2.2.1} $A$ is a flat $O$-algebra of rank $p^u$. 
\endproclaim 

\demo{Proof} First, we need the following property. 

\proclaim{Lemma 2.2.2} {\rm a)} $\hat U^{(p)}=(\hat u_{ij}^p)$ divides the scalar
matrix $pE$ in $M_u(O)$; 

b) if $V^0=(v_{ij}^0)\in M_u(O)$ is such that $\hat U^{(p)}V^0=pE$
then for any $1\leqslant i,r,j\leqslant u$, 
$$\hat u_{ir}v^0_{rj}\equiv 0\operatorname{mod}\pi .$$  
\endproclaim 

\demo{Proof of lemma} Let $V=(v_{ij})\in M_{u}(S)$ be such that 
$UV=t^eE$. Choose $\hat v_{ij}\in O$ such that 
for all $1\leqslant i,j\leqslant u$, 
$\hat v_{ij}\operatorname{mod}p=\kappa _{SO}(v_{ij}
\operatorname{mod}t^{ep})$. 
Then the equality $UV=t^eE$ implies that 
$\dsize\sum\Sb r\endSb \hat u_{ir}\hat v_{rj}\equiv 
\pi ^e\delta _{ij}\operatorname{mod}p$, where $\delta $ is the
Kronecker symbol.  
Now Proposition 1.4.1 implies that all products $\hat u_{ir}\hat
v_{rj}\equiv 0\operatorname{mod}\pi $ and, therefore, 
$$\sum\Sb r \endSb \hat u_{ir}^p\hat v_{rj}^p\equiv\pi ^{ep}\delta
_{ij}\operatorname{mod}p\pi .$$
This gives the existence of $v_{ij}'\in O$ such that 
$v'_{ij}\equiv \hat v_{ij}^p\operatorname{mod}\pi $ and 
$$\sum\Sb r \endSb \hat u_{ir}^pv'_{rj}=\pi ^{ep}\delta
_{ij}.$$
Therefore, we can take the matrix 
$V^0=(v^0_{ij})=(v'_{ij}p\pi ^{-ep})$ 
to satisfy the requirement $\hat U^{(p)}V^0=pE$. Clearly, the
condition b) follows from Proposition 1.4.1. 

The lemma is proved. 
$\square $
\enddemo

Continue the proof of proposition 2.2.1. Let 
$(F_1', \dots ,F_u')=((\bar Y\hat U)^{(p)}+p\bar Y)\hat
U^{(p)^{-1}}$. 
Let $J'_A$ be the ideal in $O[Y_1,\dots ,Y_u]$ generated by 
$F_1',\dots ,F_u'$. Clearly, $J'_A\otimes _OK=J_A\otimes _OK$. 
By above Lemma 2.2.2 all $F_i'\in O[Y_1,\dots ,Y_u]$ and, therefore, 
$J'_A\subset J_A$.

\definition{Definition} $O^{<p}[Y_1,\dots ,Y_u]$ will denote the
$O$-submodule in $O[Y_1,\dots ,Y_u]$ generated by all monomials 
$Y^{\underline{i}}:=Y_1^{i_1}\dots Y_u^{i_u}$, 
where $\underline i=(i_1,\dots ,i_u)$ is a 
multi-index such that $0\leqslant i_1,\dots ,i_u<p$. 
\enddefinition 

\proclaim{Lemma 2.2.3} With the above notation 
$$O[Y_1,\dots ,Y_u]=\oplus\Sb k_1,\dots ,k_u\geqslant 0\endSb 
O^{<p}[Y_1,\dots ,Y_u]F_1^{\prime k_1}\dots F_u^{\prime k_u}.$$
\endproclaim 
\demo{Proof of lemma} 
First, prove that for all $1\leqslant i\leqslant u$, 
$F_i'=Y_i^p+G_i'$, where 
$G_i'\operatorname{mod}\pi \in k[Y_1,\dots ,Y_s]$ are linear
polynomials. 
Indeed, the non-linear terms of the polynomial 
$F'_i-\sum\Sb j\endSb Y_j^p\hat
u_{ji}^p$ have coefficients divisible by elements of the form 
$p\hat u_{j_1i}\dots \hat
u_{j_pi}$. By above Lemma 2.2.2, ${U^{(p)}}^{-1}=(v^0_{ij}/p)$. Therefore, 
the coefficients of non-linear terms of $F_{i}'$ are linear
combinations of $p\hat u_{j_1i}\dots \hat u_{j_pi}v^0_{ij}/p\equiv
0\operatorname{mod}\pi $ because $\hat u_{j_1i}v^0_{ij}\equiv
0\operatorname{mod}\pi $. 

Now the division algorythm in each variable $Y_1,\dots
,Y_u$ gives the required decomposition modulo $\pi $. This
immediately implies the required decomposition on 
the level of $O$-modules. The lemma is proved.  
$\square $
\enddemo

Lemma 2.2.3 implies that the projection $\pr
_{\bar 0}$ of 
$O[Y_1,\dots ,Y_u]$ onto the $(0,\dots ,0)$-component 
$O^{<p}[Y_1,\dots ,Y_u]$ of the corresponding decomposition has the kernel $J'_A$ and 
it identifies $O[Y_1,\dots ,Y_u]/J'_A$ with the flat $O$-module 
$O^{<p}[Y_1,\dots ,Y_u]$.  

The embedding $J'_A\subset J_A$ induces an epimorphic map 
of $O$-modules 
$$\alpha :O^{<p}[Y_1,\dots ,Y_u]\longrightarrow A.$$ 
But $J'_A\otimes _OK=J_A\otimes _OK$ implies that 
$\Ker\alpha\otimes _OK=0$ (because $O^{<p}[Y_1,\dots ,Y_u]$ has no
$O$-torsion). Therefore, $\Ker\alpha =0$, $J_A=J'_A$ and the
proposition is  proved. 
$\square $ 
\enddemo

\remark{Remark} 
Notice that for any $1\leqslant i\leqslant u$, one has  
$dF_i=p(1+H_i)dY_i$, where all $H_i$ belong to the maximal ideal of  
the ring of formal power series $O[[Y_1, \dots ,Y_u]]$. Therefore, 
$dF_i$,  $1\leqslant i\leqslant u$, form an $K[[Y_1,\dots ,Y_u]]$-basis of 
$\Omega ^1_{K[[Y_1,\dots ,Y_u]]/K}$ and 
$A_{K}=A\otimes _{O}K$ is etale over $K$. 
\endremark 

\definition{Definition} For a given $\Cal M=(M^0,M^1,
\varphi _1)\in\MF _S^e$ denote by $\Cal A(\Cal M)$ the family of
$O$-algebras obtained by the above procedure for all 
choices of a special basis in $M^1$ and the corresponding 
lift $\hat U\in M_u(O)$ of the matrix $\kappa _{SO}(U\operatorname{mod}t
^{ep})\in M_u(O/pO)$. 
\enddefinition 

\subsubhead{\rm 2.3.} $\varphi _1$-nilpotent lifts $\theta ^{DP}_A$  
\endsubsubhead 

Suppose $\Cal M=(M^0,M^1,\varphi _1)\in\MF _S^e$ and 
$A\in\Cal A(\Cal M)$ is given in the notation from n.2.2. 
Consider $\iota ^{DP}(A)=(I_A/I_A^{DP}, I_A(p)/I_A^{DP},\varphi
_1)\in\FM _S$. Define the $S$-linear morphism 
$\theta _A^0:M^0\longrightarrow I_A/I_A^{DP}$ by 
the correspondences $m_i=\varphi _1(m_i^1)\mapsto
Y_i\operatorname{mod}I_A^{DP}$, $1\leqslant i\leqslant u$. 
Then $\theta _A^0$ induces $\theta _A^1:M^1\longrightarrow
I_A(p)/I_A^{DP}$, which 
is also uniquely determined by the correspondences 
$m_i^1\mapsto Z_i\operatorname{mod}I_A^{DP}$, $1\leqslant i\leqslant
u$. So,

$$\theta _A^0(M^0)=N^0=\left\{\sum\Sb i\endSb o_iY_i\operatorname{mod}I_A^{DP}\ |\ o_1,\dots
,o_u\in O\right\},$$ 
$$\theta _A^1(M^1)=N^1=\left\{\sum\Sb i\endSb o_iZ_i\operatorname{mod}I_A^{DP}\ |\ o_1,\dots
,o_u\in O\ \right\},$$
where for all $1\leqslant i\leqslant u$, 
$$Z_i=\sum\Sb j\endSb Y_j\hat u_{ji}. \tag{2.3.1}$$

Clearly, $\varphi _1:I_A(p)/I_A^{DP}\longrightarrow I_A/I_A^{DP}$ 
induces $\varphi _1:N^1\longrightarrow N^0$ and we obtain 
$\Cal N=(N^0,N^1,\varphi _1)\in\FM _S$ together with the natural
embedding $\Cal N\longrightarrow\iota ^{DP}(A)$ in the category $\FM
_S$. On the other hand, it is not obvious that 
$\theta _A^{DP}:=(\theta _A^0,\theta _A^1)$ gives a morphism from
$\Cal M$ to $\Cal N$ in the category $\FM _S$: we must verify the
compatibility of $\theta _A^{DP}$ with $\varphi _1$'s in $\Cal M$ and
$\Cal N$. As a matter of fact, we have more. 

\proclaim{Proposition 2.3.2} $\theta ^{DP}_A$ is a $\varphi
_1$-nilpotent 
morphism in the category $\FM _S$. 
\endproclaim 

\demo{Proof} Consider the map $\tilde{\iota }_A:M^0\longrightarrow
A\otimes O/pO$ given for $1\leqslant i\leqslant u$, by the 
correspondences $m_i\mapsto\widetilde Y_i:=Y_i\operatorname{mod}p$. If
$\widetilde{\Cal M}:=\tilde{\iota }_A(\Cal
M)=(\widetilde{M}^0,\widetilde{M}^1,\varphi _1)$, then 

$\bullet $\ $\widetilde{M}^0$ is a free $O/p=\kappa _{SO}(S/t^{ep})$-module
with the basis
$\widetilde{Y}_1,\dots ,\widetilde{Y}_u$; 

$\bullet $\  $\widetilde{M}^1$ is
generated over $O/p$ by
$\widetilde{Z}_i:=Z_i\operatorname{mod}pA$, where $i=1,\dots ,u$ and
$Z_1,\dots ,Z_u$ are given by above relations (2.3.1); 

$\bullet $\ $\varphi _1:\widetilde{M}^1\longrightarrow \widetilde{M}$
is a unique $\sigma $-linear map such that 
$\widetilde{Z}_i\mapsto\widetilde{Y}_i$, $1\leqslant i\leqslant u$. 

Clearly, $\tilde{\iota }_A:\Cal M\longrightarrow \tilde{\iota }_A(\Cal M)$ is 
$\varphi _1$-nilpotent (use that $p>2$ and $\varphi
_1(t^{ep}M^0)\subset\varphi _1(t^{e(p-1)}M^1)\subset t^{ep(p-1)}M^0$). 
So, if $h:\widetilde{M}\longrightarrow \theta ^{DP}_A(\Cal M)$ is the natural
projection and $T=\Ker h:\widetilde{M}^0\longrightarrow N^0$ then it 
will be sufficient to prove that $\Ker h|_{\widetilde M^1}=T$, $\varphi _1(T)\subset T$
and $\varphi _1|_T$ is nilpotent. 

\proclaim{Lemma 2.3.3} If for $o_1,\dots ,o_u\in O$, 
$\dsize\sum\Sb i\endSb o_iY_i\in I_A(p)$ then 
$\dsize\sum \Sb i\endSb o_i\widetilde{Y}_i\in\widetilde{M}^1$.
\endproclaim 

\demo{Proof} If $\sum\Sb i\endSb o_iY_i\in I_A(p)$ then $\sum\Sb
i\endSb o_i^pY_i^p\in pI_A$. Consider the generators
$F_i'=Y_i^p+G_i'$, $1\leqslant i\leqslant u$, of the ideal $J_A$ 
from the proof of proposition 2.2.1. Then 
$$\sum\Sb i\endSb o_i^pG_i'\in 
pO[Y_1,\dots ,Y_u].\tag{2.3.4}$$ 
Indeed, 
$\sum\Sb i\endSb o_i^pY_i^p\equiv -\sum\Sb i\endSb
o_i^pG_i'\operatorname{mod}J_A$ 
and the polynomial from the right hand side is a canonical
presentation of the element from the left hand side as a polynomial 
from $O^{<p}[Y_1,\dots ,Y_u]$. 

Notice now that the linear terms of the coordinates of the vector 
$(G_1',\dots ,G_u')$ are equal to $p\bar Y(\hat U^{(p)})^{-1}$. 
Therefore, above condition (2.3.4) implies that 
$$\hat U^{(p)^{-1}}\left (\matrix o_1^p\\\dots \\ o_u^p\endmatrix
\right )
=\left (\matrix \alpha _1\\ \dots \\ \alpha _n\endmatrix \right )$$
with all $\alpha _1,\dots ,\alpha _u\in O$. 

Clearly, for $1\leqslant i\leqslant u$, there are 
$\alpha _i'\in O$ such that $\alpha _i^{\prime p}
\equiv\alpha _i\operatorname{mod}p$. Then we obtain 
$$\left (\matrix o_1 \\ \dots \\ o_u\endmatrix\right )\equiv \hat
U\left (\matrix \alpha _1' \\ \dots \\ \alpha _u'\endmatrix\right )
\operatorname{mod}\pi ^e$$
and, therefore, 
$$\sum\Sb i\endSb o_i\widetilde{Y}_i=\widetilde{\bar Y}\left (\matrix
o_1 \\ \dots \\ o_u\endmatrix\right )\equiv 
\widetilde{\bar Y}\hat U\left (\matrix\alpha _1' \\ \dots \\ \alpha
_u'\endmatrix\right )
\equiv \widetilde{\bar Z}\left (\matrix \alpha _1' \\ \dots \\ \alpha
_u' \endmatrix \right )\operatorname{mod}\pi ^e.$$

In other words, 
 $\sum\Sb i\endSb o_i\widetilde{Y}_i$ is an 
$O$- linear combination of
$\widetilde{Z}_1,\dots ,\widetilde{Z}_u$ modulo $\pi
^e\widetilde{M}^0$ and it remains to notice that $\pi
^e\widetilde{M}^0\subset\widetilde{M}^1$. 

The lemma is proved.
$\square $
\enddemo 

$\bullet $\ \ Prove that $T=\Ker h|_{\widetilde{M}^1}$. 
\medskip 

Suppose $o_1,\dots ,o_u\in O$ are such that 
$\sum\Sb i\endSb o_i\widetilde{Y}_i\in\Ker h$. Then 
$\sum\Sb i\endSb o_iY_i\in I_A^{DP}$. In particular, 
$\sum\Sb i\endSb o_iY_i\in I_A(p)$ and $\sum\Sb i\endSb
o_i\widetilde{Y}_i\in\widetilde{M}^1$ by the above lemma. 
So, $T\subset \widetilde{M}^1$ and, therefore, 
$T=\Ker h|_{\widetilde{M}^1}$. 
\medskip 

$\bullet $\ \ Prove that $\varphi _1(T)\subset T$. 
\medskip 

Let $J^*$ be the ideal in $A$ generated by $p$ and all products 
$Z_{i_1}\dots Z_{i_p}$, where 
$1\leqslant i_1,\dots ,i_p\leqslant u$. Because $p>2$ 
and all $Z_i\in I_A(p)$, it holds $J^*\subset
I_A^{DP}$. 

Suppose $o_1,\dots ,o_u\in O$ and 
$\tilde m= \sum\Sb i\endSb o_i\widetilde{Z}_i\in T$. Then $\sum\Sb
i\endSb o_iZ_i\in I_A^{DP}$ and 
$\varphi _1(\tilde m)$ is the image in $I_A/pI_A$ of 
$\sum\Sb i\endSb o_i^pY_i=-p^{-1}\left (\sum\Sb i\endSb
o_iZ_i\right )^p+j^*$, 
where $j^*\in J^*$ (use that $Z_i^p+pY_i=0$). 
This element belongs to $I^{DP}_A$ and, therefore, 
$\varphi _1(\tilde m)\in I_A^{DP}\operatorname{mod}pI_A$,
i.e. $\varphi _1(\tilde m)\in T$. 
\medskip 

It remains to prove that $\varphi _1|_T$ is nilpotent. 

First, introduce the $O$-subalgebra $A'$ of $A$ generated by
$Z_1,\dots ,Z_u$. It can be described as the quotient of 
the polynomial ring $O[Z_1,\dots ,Z_u]$ by the ideal 
generated by all 
$$Z_i^p+pY_i=Z_i^p+p\pi ^{-e}\sum\Sb j\endSb Z_j\tilde v_{ji},\ \
1\leqslant i\leqslant u, $$
where $\tilde V=(\tilde v_{ij})\in M_u(O)$ is such that 
$\hat U\tilde V=\pi ^eE$. (The existence of $\tilde V$ follows from the
existence of $V\in M_u(S)$ such that $UV=t^eE$.) This implies that 
any element $b$ of $A'$ can be written uniquely as 
$b=\sum\Sb \underline i\endSb o_{\underline i}Z^{\underline i}\in
O^{<p}[Z_1,\dots ,Z_u]$ (as earlier, here 
$\underline i=(i_1,\dots ,i_u)$, $0\leqslant i_1,\dots
,i_u<p$, all $o_{\underline i}\in O$ and $Z^{\underline
i}=Z_1^{i_1}\dots Z_u^{i_u}$). 

Let $I_{A'}=(Z_1,\dots ,Z_u)$ be the augmentation ideal of $A'$.  

\proclaim{Lemma 2.3.5} If $\alpha _1,\dots ,\alpha _u\in O$ and 
$\alpha _1Z_1+\ldots +\alpha _uZ_u\in I_{A'}^p+pA'$ then all $\alpha _i\in
p\pi ^{-e}O$. 
\endproclaim 

\demo{Proof} We have 
$\alpha _1Z_1+\dots +\alpha _uZ_u+pb=a\in I_{A'}^p$, where  we can assume
that  
$b=\sum\Sb \underline i\endSb o_{\underline i}Z^{\underline i}\in
O^{<p}[Z_1,\dots ,Z_u]$. 

Notice that $a$ is an $O$-linear combination of the terms 
$Z_1^{j_1}\dots Z_u^{j_u}$ with 
\linebreak 
$j_1+\ldots +j_u\geqslant p$. If all 
$j_i<p$ then such a term can contribute only to the coefficient 
$o_{\underline i}$ from the above decomposition of $b$ with 
$\underline i=(j_1,\dots ,j_u)$ and does not affect $\alpha
_1,\dots ,\alpha _u$. If for some index $i$, 
$j_i\geqslant p$ then $Z_i^p$ must be replaced by 
$-p\pi ^{-e}\sum\Sb j\endSb Z_j\tilde v_{ji}$ and a possible
contribution to $\alpha _1,\dots ,\alpha _u$ is zero modulo 
$p\pi ^{-e}$. 

The lemma is proved. 
$\square $
\enddemo 

$\bullet $\ \ Prove that $\varphi _1|T$ is nilpotent. 
\medskip  

Suppose $m_0\in T$. For all $i\geqslant 0$, set $m_{i+1}=\varphi
_1(m_i)$. 

Choose $\hat m_0=\sum\Sb k\endSb o_{k0}Z_k$ such that 
$\hat m_0\operatorname{mod}pI_A=m_0$. 

Then $\hat m_0\in I_A^{DP}$ and, therefore, if $u_0=\hat m_0$ and for 
$i\geqslant 0$, $u_{i+1}=-u_i^p/p$, then $\dsize\lim _{i\to\infty
}u_i=0$. 
Notice that all $u_i$ belong to the augmentation ideal 
$I_{A'}$ of the above defined $O$-algebra $A'=O[Z_1,\dots ,Z_u]$. 

For $i>0$, define $\hat m_i=\sum\Sb k\endSb
o_{ki}Z_k$ with $o_{ki}\in O$, by the relation 
$$\hat m_{i}=\sum\Sb k\endSb o_{k,i-1}^pY_k=\sum\Sb k\endSb
o_{ki}Z_k.$$
Clearly, $m_i=\hat m_i\operatorname{mod}pI_A$ and there are 
$j_i\in I_{A'}^p$ such that 
$$\hat m_{i+1}=-\frac {\hat m_i^p}{p}+j_i.$$
This means that for any $i\geqslant 0$, 
$\hat m_i\equiv u_i\operatorname{mod}I_{A'}^p$. Therefore, there is an
$i_1\geqslant 0$ such that for any $i\geqslant i_1$, 
$\hat m_i\in I_{A'}^p+pI_{A'}$. Now Lemma 2.3.5 implies that 
$m_{i_1}\in\pi ^{e(p-1)}\widetilde{M}^1$ and  
$\varphi _1 |T$ is nilpotent because $\varphi _1(\pi
^e\widetilde{M}^1)\subset\pi ^{ep}\widetilde{M}
\subset \pi ^{e(p-1)}\widetilde{M}^1$ 
and $p>2$. 

The proposition is proved. 
$\square $

\enddemo

\subsubhead{\rm 2.4.} The functor 
$\Cal G_{O}:\MF _S^e\longrightarrow\Gr _{O}$ 
\endsubsubhead 

Consider $\Cal M=(M^0,M^1,\varphi _1)\in\MF _S^e$, $A\in\Cal A(\Cal M)$
and $B\in\Aug _{O}$. 

\proclaim{Lemma 2.4.1} The correspondence 
$f\mapsto \theta ^{DP}_A\circ\iota ^{DP}(f)$
induces a bijective map from  
$\Hom _{\Aug _{O}}(A,B)$ to $\Hom _{\FM _S
}(\Cal M,\iota ^{DP}(B))$. 
\endproclaim 

\demo{Proof} Suppose $A=O[Y_1,\dots ,Y_u]$ is given in notation of
n.2.2. 
 Then by considering the images $\bar c$ of the 
vector $\bar Y=(Y_1,\dots ,Y_u)$ in $(I_B)^u$   we obtain 
$$\Hom _{\Aug _{O}}(A,B)=\{\bar c\in (I_B)^u, \bar c'\in
I_B(p)^u\ | -\bar c^{\prime (p)}/p=\bar c,\ \bar c'=\bar c\hat U\ \}.$$

Similarly, 
$$\Hom_{\FM _S}(\Cal M,\iota ^{DP}(B))=\{\bar b\in
(I_B/I_B^{DP})^u, \bar b'\in
(I_B(p)/I_B^{DP})^u
\ |\ \varphi _1(\bar b')=\bar b, \bar b'=\bar b\hat U\}.$$ 

The correspondence $f\mapsto \theta ^{DP}_A\circ\iota ^{DP}(f)$ is given by the projections  
$\bar c\mapsto \bar c\operatorname{mod}I_B^{DP}$ and 
$\bar c'\mapsto\bar c'\operatorname{mod}I_B^{DP}$. Therefore, 
our lemma follows from proposition 2.1.1. 
$\square $
\enddemo 

\proclaim{Proposition 2.4.2} Suppose $\Cal M_1,\Cal M_2\in\MF _S^e$,  
$A_1\in\Cal A(\Cal M_1)$ and $A_2\in\Cal A(\Cal M_2)$. Then 
\newline 
{\rm a)} for any $g\in\Hom _{\MF _S^e}(\Cal M_1,\Cal M_2)$, there is 
a unique $f\in\Hom _{\Aug _{O}}(A_1,A_2)$ such that 
$\theta ^{DP}_{A_1}\circ \iota ^{DP}(f)=g\circ \theta ^{DP}_{A_2}$;
\newline 
{\rm b)} with the above notation, the correspondence $g\mapsto f$ 
gives an embedding 
$$p_{A_1A_2}:\Hom _{\MF _S^e}(\Cal M_1,\Cal M_2)\longrightarrow\Hom
_{\Aug _{O}}(A_1,A_2).$$
\endproclaim

\demo{Proof} a) follows from Lemma 2.4.1 applied to 
$g\circ \theta ^{DP}_{A_2}\in\Hom _{\FM _S}(\Cal M_1,\iota ^{DP}(A_2))$. In order to prove b)
notice that 
$\iota ^{DP}(f)(\theta ^{DP}_{A_1}(\Cal M_1))=(g\circ \theta ^{DP}_{A_2})(\Cal M_1)\subset
\theta ^{DP}_{A_2}(\Cal M_2)$. 
Therefore, $g$ appears as a unique $\varphi _1$-nilpotent lift of 
$$\theta ^{DP}_{A_1}\circ \iota ^{DP}(f)\in\Hom _{\FM _S}(\Cal M_1,\theta
^{DP}_{A_2}(\Cal M_2)).$$ 

$\square $
\enddemo

Let $\Cal M\in\MF _S^e$ and $A\in\Cal A(\Cal M)$. Notice that 
$A\otimes _{O}A\in\Cal A(\Cal M\oplus\Cal M)$. Indeed, suppose $A$ 
is constructed via 
the special basis $m_1^1,\dots ,m_u^1$ of $M^1$ and the corresponding
matrix $\hat U\in M_u(O)$. Then $A\otimes _{O}A$ will appear 
in $\Cal A(\Cal M\oplus\Cal M)$ via the special basis 
$(m_1^1,0),\dots ,(m_u^1,0),(0,m_1^1),\dots ,(0,m_u^1)$ of $M^1\oplus M^1$
and the
corresponding matrix $\dsize\left (\matrix \hat U & 0 \\ 0 & \hat U
\endmatrix\right )$. 

Then proposition 2.4.2 immediately implies that 
$\Spec A$ can be provided with a structure of the object of the
category $\Gr _{O}$ by taking:

$\bullet $\ \ $p_{A,A\otimes A}(\nabla ):A\longrightarrow A\otimes A$
as a comultiplication, where 
$\nabla $ is the diagonal embedding of $\Cal M$ into $\Cal M\oplus\Cal
M$; 

$\bullet $\ \ the natural projection $A\longrightarrow A/I_A=O$ as
a counit; notice that this projection also appears in the form $p_{AO}$, where 
$O$ is considered as an element of $\Cal A(0)$ and $0=(0,0,\varphi _1)$ is the zero
object in $\MF _S^e$; 

$\bullet $ \ \ $p_{AA}(-\id _{\Cal M})$ as a coinversion. 
\medskip 

Now we can introduce the {\it functor} $\Cal G_O$. For any $\Cal M\in\MF _S^e$ choose 
$A=A(\Cal M)\in\Cal A(\Cal M)$ and set $\Cal G_O(\Cal M)=\Spec A(\Cal
M)$ with the above defined structure of the object of the category
$\Gr _{O}$. If $\Cal M_1,\Cal M_2\in\MF _S^e$ and 
$f\in\Hom _{\MF _S^e}(\Cal M_1,\Cal M_2)$ then set 
$\Cal G_O(f)=p_{A(\Cal M_1)A(\Cal M_2)}(f)$. 

We can also use proposition 2.4.2 to prove that under any other
choice 
of representatives $A'(\Cal M)\in\Cal A(\Cal M)$, the corresponding
functor $\Cal G'_O$ will be naturally equivalent to $\Cal G_O$. 
\medskip 
\medskip 
\medskip 
\medskip

\subhead 3. Full faithfulness of $\Cal G_O$ 
\endsubhead 
\medskip 

\subsubhead{\rm 3.1.} Special case of $B\in\Cal A(\Cal N)$, 
$\Cal N\in\MF _S^e$ 
\endsubsubhead 
\medskip 

Let $\Cal N=(N^0,N^1,\varphi _1)\in\MF _S^e$. Choose an $S$-basis 
$n_1^1,\dots ,n_u^1$ of $N^1$ such that for $1\leqslant i\leqslant u$, 
there are $\tilde s_i\in S$ such that 
$n_i=n_i^1\tilde s_i^{-1}$ form a basis of $N^0$.  
All $\tilde s_i$ divide $t^e$ because $N^1\supset t^eN^0$. Notice that 
$n_1^1,\dots ,n_u^1$ is a special basis, cf. n.1.4. Indeed, for any $1\leqslant
i\leqslant u$, 
$n_i^1\in tN^0$ iff $\tilde s_i\equiv
0\operatorname{mod}t$. Therefore, the non-zero elements of the set 
$\{n_i^1\operatorname{mod}tN^0\ |\ \tilde s_i\in S^*\}$ are linearly
independent over $k$ 
because the corresponding elements $\tilde s_i^{-1}n_i$ form a part of the $S$-basis 
$n_1,\dots ,n_u$ of $N^0$. 

For $1\leqslant i\leqslant u$, let $n'_i=\varphi _1(n_i^1)$ and let 
$U\in M_u(S)$ be such that $\bar n^1=\bar n'U$. Here  
$\bar n^1=(n_1^1,\dots ,n_u^1)$ and $\bar n'=(n'_1,\dots ,n'_u)$. 
Notice that $U=U_0U_1$, where $U_0\in\GL _u(S)$ is such that 
$\bar n'U_0=\bar n:=(n_1,\dots ,n_u)$ and 
$U_1=(\tilde s_i\delta _{ij})\in M_u(S)$ is diagonal. Let $D\in\GL _u(O)$
be such that $D\operatorname{mod}p
=\kappa _{SO}(U_0\operatorname{mod}t^{ep})$ and let for $1\leqslant
i\leqslant u$, $\tilde\eta _i'\in O$ be such that 
$\kappa _{SO}(\tilde s_i\operatorname{mod}t^{ep})=\tilde\eta
_i'\operatorname{mod}p$.

Denote by $\widetilde{\Omega }$ the diagonal matrix $(\tilde\eta
_i'\delta _{ij})$. 
By the results of section 2  
the coordinates of the vector 
$$((\bar YD\widetilde{\Omega })^{(p)}+p\bar Y)
(D\widetilde{\Omega })^{(p)^{-1}}\tag{3.1.1}$$
give the equations of the algebra $B=O[\bar Y]\in\Cal A(\Cal N)$. 
(This $B$ corresponds to
the above choice of basis $\bar n ^1$ and the structural matrix 
$\hat U=D\widetilde{\Omega }\in M_u(O)$.) 

Introduce the new variables $\bar X=\bar YD$ and notice that 
$(D\widetilde{\Omega
})^{(p)}=D^{(p)}\widetilde{\Omega }^{(p)}$. 
For $1\leqslant i\leqslant u$, let 
$\eta _i=-p/\tilde\eta _i^{\prime p}$. With this notation, the vector 
(3.1.1) can be rewritten as 
$(\dots ,X_i^p-\eta _iY_i ,\dots )D^{(p)^{-1}}$. 
Therefore, the algebra $B$ appears as the quotient of $O[X_1,\dots ,X_s]$ by the
ideal generated by the elements 
\linebreak 
$X_i^p-\eta _i\dsize\sum\Sb j\endSb X_jc_{ji}$, 
where $1\leqslant i\leqslant s$ and $C=(c_{ij})=D^{-1}\in\GL
_s(O)$.

\subsubhead{\rm 3.2} The description of comultiplication 
\endsubsubhead 

As we have just obtained, if $H=\Cal G_O(\Cal N)$ with  
$\Cal N=(N^0,N^1,\varphi _1)\in\MF ^e$, then  
$A(H)=B=O[X_1,\dots ,X_u]$, with the equations 
$$X_i^p-\eta _i\sum\Sb r\endSb X_rc_{ri}=0,\ \ 1\leqslant i\leqslant u,$$
where 
$u=\rk _{S}N^0=\rk _{S}N^1$, $C=(c_{ri})\in\GL _u(O)$ and  
all $\eta _i\in O$, $\eta _i|p$. 
\medskip 

\remark{Remark {\rm 3.2.1}} 
We can assume that all  
$\eta _i'$ are just integral powers of $\pi $. Indeed, 
the elements $\tilde s_i\in S$ from n.3.1 can be chosen as 
integral powers of $t$ and this will allow us to choose all $\tilde \eta _i'$
as integral powers of $\pi $.  We shall also use the notation 
$\tilde\eta _i=\tilde\eta _i^{\prime p}$, i.e. $\tilde\eta _i\eta
_i=-p$. 
\endremark 

By the definition from n.2.4,  
the comultiplication $\Delta :B\longrightarrow B\otimes _{O}B$ 
can be recovered uniquely from the conditions  
$\Delta (X_i)=
X_i\otimes 1+1\otimes X_i+j_i$, where all $j_i\in I_{B\otimes
B}^{DP}$. Using the above equations of $B$ we obtain the following 
recursive relation to recover these elements
$j_i$: 

$$\sum\Sb r \endSb j_rc_{ri}=\tilde\eta _i(\phi (X_i)+\phi (X_i\otimes
1+1\otimes X_i,j_i))+j_i^p/\eta _i,\ \ 1\leqslant i\leqslant u. \tag
{3.2.2} $$

Here $\phi (X,Y)=p^{-1}(X^p+Y^p-(X+Y)^p)$ is the first Witt
polynomial and for all $1\leqslant i\leqslant u$,  
$\phi (X_i\otimes 1, 1\otimes X_i)$ is denoted just by $\phi (X_i)$, 
cf. the basic notation in the end of the Introduction. 

Let $\Cal J_B$ be the ideal in $B\otimes _OB$ generated by the
elements $\tilde\eta _iX_i^r\otimes X_i^{p-r}$, where $1\leqslant
i\leqslant u$ and $1\leqslant r<p$. Notice that all 
$\eta _i\phi (X_i)\in \Cal J_B$. The definition of the ideal $\Cal
J_B$ depends on the chosen construction of the $O$-algebra $B$. But
because $\tilde\eta _i'X_i\in I_B(p)$, all 
$\tilde\eta _iX_i^r\otimes X_i^{p-r}\in I_{B\otimes B}(p)^p$. 
This will allow us (if necessary) to replace $\Cal J_B$ by 
the bigger invariant ideal $I_{B\otimes B}(p)^p$. 
\medskip

\proclaim{Proposition 3.2.3} For all $1\leqslant i\leqslant u$, 
$j_i\in\Cal J_B$.
\endproclaim 

\demo{Proof}  As we above have noticed, for all  
$i$, $\tilde\eta _i\phi (X_i)\in \Cal J_B\subset I^{DP}_{B\otimes B}$. 

Consider the sequence $J_n$, $n\geqslant 0$, of ideals in $B\otimes B$ 
such that 
$J_0=I_{B\otimes B}^{DP}$ and for all $n\geqslant 0$, 
$J_{n+1}=I_{B\otimes B}(p)J_n+\varphi _1(J_n)$, where the ideal 
$\varphi _1(J_n)$ is generated by the elements $\{j^p/p\ |\
j\in J_n\}$. Then the recursive relation (3.2.2) implies that 
for all
$n\geqslant 0$ and $1\leqslant i\leqslant u$, 
$j_i\in \Cal J_B+J_n$. The proposition follows then
from the fact that the intersection of all $J_n$ is the zero ideal. 
$\square $
\enddemo 

\subsubhead {\rm 3.3.} The ideals $I_B(\alpha )$, $\alpha\in O$ 
\endsubsubhead

Notice that by Lemma 2.2.3 any element of $B$ can be uniquely written
as a linear combination 
$\dsize\sum\Sb \underline i\endSb o_{\underline i}
X^{\underline i}\in O^{<p}[X_1,\dots ,X_u]$, 
where as earlier, $\underline i=(i_1,\dots ,i_u)$, $0\leqslant
i_1,\dots ,i_u<p$, $X^{\underline i}=X_1^{i_1}\dots
X_u^{i_u}$ and all coefficients $o_{\underline i}$ belong to $O$.  
 
Our system of generators $X_1,\dots ,X_u$ has a very
special property: if for all $i$, $X_i'=X_i^p/\eta _i$ then 
$X_1',\dots ,X_u'$ is obtained from $X_1,\dots ,X_u$ by a non-degenerate 
linear transformation (given by the matrix $C\in\GL _u(O)$). Then 
any element of $A$ can be written uniquely 
as a (similar to just described) linear combination 
$\dsize\sum\Sb \underline i\endSb o_{\underline i}
X^{\prime\underline i}\in O^{<p}[X_1', \dots ,X_n']$, 
where $X^{\prime\underline i}=
X_1^{\prime i_1}\dots X_{i_u}^{\prime i_u}$. This follows from
the fact that $B$ is the quotient of $O[X_1',\dots ,X_u']$ by 
the ideal generated by the elements 
$$X_i^{\prime p}-\sum\Sb j\endSb \eta _jX_j'c_{ji}^p+pH_i,$$
where $1\leqslant i\leqslant u$ and all $H_i$ are polynomials in
$X_1',\dots ,X_u'$ of total degree $\leqslant p$.

\definition{Definition} 
For $\alpha\in O$, set  
$$I_B(\alpha )=\left\{\sum\Sb \underline i\endSb o_{\underline
i}X^{\underline i}\in O^{<p}[X_1,\dots ,X_u]\ |\ o_{\underline i}\in O,\ o_{\underline i}^p
X^{\underline ip}\in \alpha I_B\right\}.$$ 
\enddefinition 

Notice that,  
$\dsize\sum\Sb \underline i\endSb o_{\underline i}
X^{\underline i}\in I_B(\alpha )$ if and only if for all multi-indices 
$\underline i=(i_1,\dots ,i_u)$ 
(where always $0\leqslant i_1,\dots ,i_u<p$), it holds 
$o_{\underline i}\in O$,  
$o^p_{\underline i}\in\alpha\eta _1^{-i_1}\dots\eta _u^{-i_u}O$ and
$o_{(0,\dots ,0)}=0$. 

The sets $I_B(\alpha )$ depend generally on the choice of generators
$X_1,\dots ,X_u$, cf. n.3.1. But if $\alpha |p$ then $I_B(\alpha )=
\{a\in I_B\ |\ a^p\in\alpha I_B\}$ does not depend on such a
choice. Notice that we have used already in 2.1 a special case of the 
notation $I_B(\alpha )$ when $\alpha =p$.

\proclaim{Lemma 3.3.1} Suppose $o\in O$, $\alpha\in O$ and  
$oX^{\underline i}\in I_B(\alpha )$, where 
$\underline i=(i_1,\dots ,i_u)$ with $0\leqslant i_1,\dots ,i_u<p$. 
 Then for any $1\leqslant j\leqslant u$, 
$oX^{\underline i}X_j\in I_B(\alpha \eta _j)$. 
\endproclaim 

\demo{Proof} We can assume that $j=1$. Notice that 
the statement is obviously true if
$i_1<p-1$. So, assume $i_1=p-1$. 

 Use induction on the number $r$ 
of elements of the set $\{j\ |\ i_j\ne 0\}$. 

Let $r=1$. Then $oX^{\underline i}=oX_1^{p-1}\in
I_B(\alpha )$ implies that 
$o^p\eta _1^{p-1}\equiv 0\operatorname{mod}\alpha $. Then 
$$oX^{\underline i}X _1=o X_1^p\in o\eta _1I_B\subset
I_B(o^p\eta _1^p)\subset I_B(\alpha \eta _1).$$

Suppose $r>1$ and the lemma is proved for all $r'<r$. Then 
$oX_1^{p-1}X_2^{i_2}\dots X_u^{i_u}\in I_B(\alpha )$ means 
$o^p\eta _1^{p-1}\eta _2^{i_2}\dots \eta _u^{i_u}\equiv
0\operatorname{mod}\alpha $. Consider the equality 
$$oX_1^{p-1}X_2^{i_2}\dots X_u^{i_u}=\sum\Sb j\endSb o\eta
_1c_{1j}X_2^{i_2}\dots X_u^{i_u}X_j, \tag{3.3.2}$$
Then for any index $j$, 
$o\eta _1c_{1j}X_2^{i_2}\dots X_u^{i_u}\in I_B(o^p\eta
_1^p\eta _2^{i_2}\dots \eta _u^{i_u})\subset I_B(\alpha\eta
_1)$, 
and clearly $X_j\in I_B(1)$. Therefore, by the inductive assumption each
term of the sum (3.3.2) belongs to $I_B(\alpha\eta _1)$. 

The lemma is proved.
$\square $
\enddemo 

\proclaim{Corollary 3.3.3} {\rm a)} If $\alpha _1,\alpha _2\in O$ then
$I_B(\alpha _1)I_B(\alpha _2)\subset I_B(\alpha _1\alpha _2)$; 
\newline 
{\rm b)} $\forall\alpha\in O$, $I_B(\alpha )$ is an ideal in
$B$. $\square $
\endproclaim 
 
\medskip

\subsubhead{\rm 3.4.} Recovering $\Cal N\in\MF ^e_S$
\endsubsubhead 

Suppose $\Cal N\in\MF _S^e$, $H=\Cal G_O(\Cal N)$ and $B=A(H)$ is the
algebra of $H$ given in the notation and assumptions of n.3.2. 
Then use for $1\leqslant i\leqslant n$, 
the generators $X_i$ of $B$ and the generators 
$X_i\otimes 1$ and $1\otimes X_i$ of  
$B\otimes B$ to introduce the ideals $I_B(\alpha )$ and 
$I_{B\otimes B}(\alpha )$, where $\alpha\in O$. 

For any $a\in I_B$, let $\delta ^+a=\Delta (a)-a\otimes 1-1\otimes a\in
I_{B\otimes B}$. Then by Proposition 3.2.3 for $1\leqslant i\leqslant u$, 
$\delta ^+X_i\in I_{B\otimes B}(p)^p\subset I_{B\otimes B}(p^p)$. 

\proclaim{Lemma 3.4.1} Suppose $a\in I_B$ is such that $\delta ^+a\in
I_{B\otimes B}(p^p)$. Then there are $o_1,\dots ,o_u\in O$ 
such that $\dsize a\equiv\sum \Sb i\endSb
o_iX_i\operatorname{mod}I_B(p^p)$. 
\endproclaim 

\demo{Proof} We can assume that 
$\dsize a=\sum\Sb r(\underline
i)\geqslant 2\endSb o_{\underline i}X^{\underline i}\in
O^{<p}[X_1,\dots ,X_u]$,
where as earlier $\underline i=(i_1,\dots ,i_u)$, 
$0\leqslant i_1,\dots ,i_u<p$ and $r(\underline i)=
i_1+\dots +i_u$. 
Then 
$$\Delta (X^{\underline i})\equiv 
(X_1\otimes 1+1\otimes X_1)^{i_1}\dots (X_u\otimes 1+1\otimes
X_u)^{i_u}\operatorname{mod}I_{B\otimes B}(p^p).$$ 
 It is easy to see that:
\medskip 

a)  $\delta ^+X^{\underline i}$ is a linear 
combination of the terms $X^{\underline j_1}\otimes 
X^{\underline j_2}$ where $\underline j_1$ and $\underline j_2$ are
multi-indices such that $r(\underline j_1), r(\underline j_2)>0$ 
and $\underline j_1+\underline j_2=\underline i$; in addition, all 
such terms 
$X^{\underline j_1}\otimes  X^{\underline j_2}$ 
appear with coefficients from $\Bbb Z_p^*\subset
O^*$; 
\medskip 

b) any term $X^{\underline
j_1}\otimes X^{\underline j_2}$ from the above n.a) 
does not appear with a non-zero
coefficient in the
decomposition of any $\delta ^+ X^{\underline i_1}$ with 
$\underline i_1\ne \underline i$. 
\medskip 

The above two facts a) and b) imply that for any $\alpha\in O$ such
that $\alpha |p^p$, $\delta
^+a\in I_{B\otimes B}(\alpha )$ if and only if $a\in I_B(\alpha )$. 
In particular, if $\alpha =p^p$ we obtain the statement of our
lemma. $\square $
\enddemo

Let $\theta ^{DP}_B:\Cal N\longrightarrow \iota ^{DP}(B)$ 
be the morphism from n.2.3  
and $\theta ^{DP}_B(\Cal N)=(N_1^0,N_1^1,\varphi _1)\in\FM _S$. Then the above
lemma implies that  
$$\widetilde{N}^0=\{ a\in I_B\ |\ \delta ^+a\in I_{B\otimes
B}(p)^p\}/I_B(p^p)$$   
is mapped onto $N_1^0$ by the projection 
$I_B/I_B(p^p)\longrightarrow
I_B/I_B^{DP}$ induced by the embedding 
$I_B(p^p)\subset I_B^{DP}$. Therefore, we have the following
statement. 

\proclaim{Proposition 3.4.2} If $H=\Spec B=\Cal G_O(\Cal N)$ with  
$\Cal N\in\MF ^e$, then $\Cal N$ can be uniquely recovered as a  
$\varphi _1$-nilpotent lift of 
$(N_1^0/I_B^{DP}, N_1^1/I_B^{DP},
\varphi _1)$, where 
\newline 
$N_1^0=\{a\in I_B\ |\ \delta ^+(a)\in
I_{B\otimes B}(p)^p\}$, 
$N_1^1=I_B(p)\cap N_1^0$
and $\varphi _1$ is
induced by the map 
$a\mapsto-a^p/p$.$\square $
\endproclaim

\proclaim{Corollary 3.4.3} The functor $\Cal G_O$ is fully faithful. 
$\square $
\endproclaim

\remark{Remark} We could not prove that the elements of 
$N_1^0/I_B^{DP}$ come from $a\in I_B$ such that 
$\delta ^+(a)\in I_{B\otimes B}^{DP}$. This is why we use more strong
condition $\delta ^+(a)\in I_{B\otimes B}(p)^p$. As a matter of fact, 
we could use either the stronger condition $\delta ^+(a)\in \Cal
J_{B}$ or the weaker one $\delta ^+(a)\in I_{B\otimes
B}(p^p)$, but they both depend on a choice of special generators 
of $B$ and, therefore, are not functorial.
\endremark 
\medskip

\subsubhead{\rm 3.5.} A property of comultiplication 
\endsubsubhead 

Suppose $\Cal N\in\MF _S^e$ and $H=\Cal G_O(\Cal N)$ is given in
notation of n.3.2. 
Then for $1\leqslant i\leqslant
u$, $\Delta (X_i)=X_i\otimes 1+1\otimes X_i+j_i$ with $j_i\in
I_{B\otimes B}(p^p)$, cf. n.3.2. 
Remind that any element of $B$ can be uniquely written as a polynomial 
from $O^{<p}[X_1,\dots ,X_u]$ and  
any element from $B\otimes _{O}B$ can be uniquely written as 
a polynomial from 
\linebreak 
$O^{<p}[X_1\otimes 1,\dots ,X_u\otimes 1,1\otimes
X_1,\dots ,1\otimes X_u]$. We shall use the following property later in 
subsection 6.7.  

\proclaim{Proposition 3.5.1} For $1\leqslant i,r\leqslant u$, $j_i$ 
as an element of $O^{<p}[X_1\otimes 1, \dots ,1\otimes X_u]$ contains
$\phi (X_r)$ with the coefficient $\tilde\eta
_rd_{ri}$ modulo $p\tilde\eta _r$, where $(d_{ri})=C^{-1}$.  
\endproclaim 

\demo{Proof} We can proceed with $j_i$ taken modulo 
$I_{B\otimes B}(p^{2p})$ because if $\alpha\in O$ is such that 
$\alpha\phi (X_r)\in I_{B\otimes B}(p^{2p})$ then 
$\alpha ^p\eta _r^p\equiv 0\operatorname{mod}p^{2p}$ and, therefore, 
$\alpha\equiv 0\operatorname{mod}p\tilde\eta _r$. 

For any $1\leqslant i\leqslant u$, 
$j_i\in I_{B\otimes B}(p^p)$ and, therefore, 
$j_i^p/p\in I_{B\otimes B}(p^{p^2-p})\subset I_{B\otimes B}(p^{2p})$
because $p\geqslant 3$. 
In addition, 
$$\tilde\eta _i\phi (X_i\otimes 1+1\otimes X_i, j_i)
\in I_{B\otimes B}(p^{2p-1}).$$
Therefore, relation (3.2.2) implies that 
for $1\leqslant i\leqslant u$, 
$$j_i\equiv \sum\Sb r\endSb \tilde\eta _r
\phi (X_r)d_{ri}\operatorname{mod}I_{B\otimes B}(p^{2p-1}),$$
and we obtain that 
$$\sum\Sb r \endSb j_rc_{ri}-\tilde\eta _i\phi (X_i)\equiv 
-\tilde\eta _i(X_i\otimes 1+1\otimes X_i)^{p-1}\sum\Sb r\endSb 
\tilde\eta _r\phi (X_r)d_{ri}\operatorname{mod}I_{B\otimes B}(p^{2p}).$$

Notice that for $i\ne r$, the term 
$(X_i\otimes 1+1\otimes X_i)^{p-1}\phi (X_r)$ does not contribute to
the coefficient for $\phi (X_r)$. But if $i=r$ then 
$$\tilde\eta _r(X_r\otimes 1+1\otimes X_r)^{p-1}\tilde\eta _r\phi
(X_r)\in p\tilde\eta _rI_{B\otimes B},$$
 because $X_r^p\in\eta _rI_B$ and
$\tilde\eta _r\eta _r=-p$. In other words, there is no contribution 
to the coefficient for $\phi (X_i)$ modulo $p\tilde\eta _r$.  
The proposition is proved. 
$\square $
\enddemo 
\medskip 
\medskip

\subhead 4. Construction of the functor $\Cal G_{O_0}^O:\MF
^e_S\longrightarrow\Gr _{O_0}$ 
\endsubhead 

In this section we use the basic notation $K,O,\pi, K_0, O_0, \pi _0$. 
We prove that the existence of the subfield $K_0$ in $K$ implies 
that any $G=\Cal G_O(\Cal M)\in\Gr _O$ 
comes from a unique $G_0\in\Gr _{O_0}$ by the 
extension of scalars
$G_0\mapsto G=G_0\otimes _{O_0}O$. Then the correspondence 
$\Cal M\mapsto G_0$ will give rise to a fully faithful 
functor $\Cal G_{O_0}^O$ from $\MF ^e_S$ to $\Gr _{O_0}$. 
\medskip

\subsubhead{\rm 4.1.} Tamely ramified extension of scalars 
\endsubsubhead 

Suppose $K_0'$ is a tamely ramified extension of $K_0$ with the residue
field $k'\supset k$. Let $O'_0$ be the valuation ring of $K'_0$. 
Denote by $e'$ the absolute ramification index of
$K'_0$ and set $e_0=e'/e$. By replacing (if necessary)
$K'_0$ by a bigger
tamely ramified extension we can always assume that $K'_0/K_0$ 
is Galois and there
is a uniformizer 
$\pi _0'\in O'_0$ such 
that ${\pi _0^{\prime}}^{e_0}=\pi _0$. Then we can 
introduce $K'=K(\pi ')$ with ${\pi ^{\prime }}^{p}=\pi _0^{\prime }$ such that
$\pi ={\pi ^{\prime}} ^{e_0}$. Set $O'=O_{K'}$ and 
$\Gamma =\Gal (K'_0/K_0)\simeq\Gal (K'/K)$. 

Notice, first, that there is the following necessary condition for the
existence of a descent of $G'_0\in\Gr _{O'_0}$ to
$O_0$ or 
, resp., of $G'\in\Gr _{O'}$ to $O$ (it holds without
the assumption that the ramification of $K'_0$ over $K_0$ is tame): 
\medskip 

$(\alpha )$\ for all $\tau\in\Gamma $ there is a $\tau $-linear
bialgebra automorphism $f_{\tau }$ of $A(G'_0)$ (or, resp., 
$A(G')$) such that 
$\forall\tau _1,\tau _2\in\Gamma $, $f_{\tau _1\tau _2}=
f_{\tau _1}f_{\tau _2}$. 
\medskip 

If $K'_0$ is unramified over $K_0$ then the above condition $(\alpha )$ is
also sufficient for the existence of a such descent. If
$K'_0$ is
totally ramified over $K_0$ of degree $n$ and $K_0$ contains 
a primitive $n$-th root of unity, then the corresponding $\Gamma $-action is
semi-simple and one can easily see that $G'_0$ admits a
descent to $O_0$, resp., $G'$ admits a descent to $O$, 
if and only if the condition $(\alpha )$ holds and
$\forall\tau\in\Gamma $, $f_{\tau }$ induces the identity maps on 
special fibres. 

With the relation to the category of filtered modules introduce the 
appropriate characteristic $p$ object $S^{\prime }=k'[[t']]$, 
where $t=t^{\prime e_0}$. Then the identification 
$\kappa _{S'O'}:S'/t^{ep}S'
\longrightarrow O'/pO'$ 
is induced by the correspondence 
$t'\mapsto \pi '$. This identification allows us to identify the
Galois group $\Gamma $ with the Galois
group of $\Frac S'$ over $\Frac S$. 
In this situation we have the obvious functor of extension of scalars 
 $\otimes _{S}S':\MF _{S}^e\longrightarrow\MF ^{e'}_{S'}$  
and, clearly, if $\Cal M\in\MF ^e_{S}$ then 
$\Cal G_{O'}(\Cal M\otimes _{S}S')=\Cal G_O(\Cal M)\otimes
_{O}O'$. 
\medskip 

4.2. The following proposition allows us to pass to tamely ramified
extensions when studying the image of the functor $\Cal G_O$. 

\proclaim{Proposition 4.2.1} $G\in\Gr _{O}$ is in the image of the
functor 
$\Cal G_O:\MF ^e_{S}\longrightarrow\Gr _{O}$ if and only if 
$G'=G\otimes _{O}O'$ is in the image of the functor 
$\Cal G_{O'}:\MF ^{e'}_{S'}\longrightarrow\Gr _{O'}$. 
\endproclaim 

\demo{Proof} It will be sufficient to consider separately 
the cases of an unramified extension $K'_0/K_0$ and a 
totally ramified extension $K'_0/K_0$ 
of degree $n$ where $n$ is prime to $p$ and 
$K_0$ contains a primitive $n$-th root of unity. 

Clearly, if $G=\Cal G_O(\Cal M)$, where $\Cal M\in\MF ^e_{S}$, then 
$G'=\Cal G_{O'}(\Cal M\otimes _{S}S')$. 

Now assume that $G'=\Cal G_{O'}(\Cal M')$, $\Cal M'\in\MF
^{e'}_{S'}$. We must prove that 
there is an $\Cal M\in\MF _{S}^e$ such that 
$\Cal M'=\Cal M\otimes _{S}S'$. 

For $\tau\in\Gamma $, consider 
$\tau $-linear bialgebra automorphisms $f_{\tau }$ 
such that $\forall\tau _1,\tau _2\in\Gamma $, $f_{\tau _1\tau
_2}=f_{\tau _1}f_{\tau _2}$ and 
$A(G)=\{a\in A(G')\ |\ \forall\tau\in\Gamma ,\  f_{\tau }
(a)=a\}$. 

By Proposition 3.4.2 there are induced $\tau $-linear automorphisms of 
$\theta ^{DP}_{A(G')}(\Cal M')$ in the category $\FM _{S'}$ and by Proposition 1.2.1 
they give rise to $\tau $-linear automorphisms $g_{\tau }\in\Aut _{\MF
^{e'}_{S'}}(\Cal M')$ such that 
for any $\tau _1,\tau _2\in\Gamma $, $g_{\tau _1\tau _2}=g_{\tau
_1}g_{\tau _2}$.

If $\Cal M'=(M^{\prime 0},M^{\prime 1},\varphi _1)$, then for any 
$\tau\in\Gamma $, $g_{\tau }$ is a $\tau $-linear automorphism of 
the $S'$-module $M^{\prime 0}$,  
$g_{\tau }(M^{\prime 1})=M^{\prime 1}$ and $\varphi _1g_{\tau
}|_{M^{\prime 1}}=g_{\tau }|_{M^{\prime 1}}\varphi _1$. 

If $K'_0/K_0$ is unramified, this already implies that for 
$M^0:=$ 
\linebreak 
$\{m\in M^{\prime 0}\ |\ \forall\tau\in\Gamma , g_{\tau
}(m)=m\}$ and $M^1:=M^0\cap M^{\prime 1}$, it holds 
\linebreak 
$\Cal M=(M^0,M^1,\varphi _1|_{M^1})\in\MF _{S}^e$ and 
$\Cal M\otimes _{S}S'=\Cal M'$. 

If $K^{\prime 0}/K^0$ is totally ramified then the $\Gamma $-action is semi-simple 
and (under our assumptions) there is an  
$S'$-basis of $M^{\prime 1}$ such that for any 
$\tau\in\Gamma $, $g_{\tau }(m_i')=\chi _i(\tau
)m_i'$, where $1\leqslant i\leqslant u$ and $\chi
_i:\Gamma\longrightarrow k^*$ are 1-dimensional characters. Then 
the elements of the
$S'$-basis $\varphi _1(m_1'),\dots ,\varphi _1(m_{u}')$ of 
$M^{\prime 0}$ satisfy the conditions $g_{\tau }(\varphi _1(m_i'))=
\chi _i(\tau )^p\varphi _1(m_i')$, 
\linebreak 
$1\leqslant i\leqslant u$. Therefore, 
$A'=A(G')$ appears in the form $O'[X_1,\dots ,X_{u}]$, 
where for all $\tau\in\Gamma $, there are induced
$\tau $-linear automorphisms $h_{\tau }:A'\longrightarrow A'$ 
such that $h_{\tau }(X_i)=[\chi _i(\tau )]X_i$, 
$1\leqslant i\leqslant u$. 
(Here $[\alpha ]\in O'$ is the Teichm\"uller representative of
$\alpha\in k$.) 

Notice that for all $\tau\in\Gamma $, $h_{\tau }=f_{\tau }$ by the
uniqueness property from Proposition 2.4.2.  It remains to notice that the
elements of $A'$ are presented uniquely as polynomials from 
$O^{\prime <p}[X_1,\dots ,X_u]$. Therefore, $A'$ can be descended to
the $O$-algebra $A(G)$ if and only if all characters $\chi _i$ are
trivial. So, there is an $\Cal M\in\MF ^e_{S}$ such that 
$\Cal M\otimes _{S}S'=\Cal M'$. 

The proposition is proved. 
$\square $
\enddemo 
\medskip 

4.3. Suppose $G_0\in\Gr _{O_0}$ and $G=G_0\otimes _{O_0}O\in\Gr _{O}$. 

\proclaim{Proposition 4.3.1} If $H_0\in\Gr _{O_0}$ is such that 
$H_0\otimes _{O_0}O=G$ then $G_0=H_0$. 
\endproclaim 

\demo{Proof} 
Let $V=G(\bar K)$ be the $\Gamma _{K}$-module of $\bar K$-points
of $G$. Then there is a canonical embedding $A(G)\subset\Map
^{\Gamma _{K}}(V,\bar K)$ given for any $a\in A(G)$, by the
correspondence 
$$a\mapsto\{v\mapsto a(v)\ |\ \forall v\in V\}.$$
The existence of 
$G_0\in\Gr _{O_0}$ such that $G_0\otimes _{O_0}O=G$ implies the existence 
of a $\Gamma _{K_0}$-module $V_0$ such that $V_0|_{\Gamma _{K}}=V$ and 
$A(G_0)=A(G)\cap\Map ^{\Gamma _{K_0}}(V_0,\bar K)$ with respect to the
natural embedding 
$$\Map ^{\Gamma _{K_0}}(V_0,\bar K)\subset \Map ^{\Gamma _{K}}(V_0|_{\Gamma
_{K}},\bar K)=\Map ^{\Gamma _{K}}(V,\bar K).$$
Therefore, it will be sufficient to prove that if $V_0'$ is 
any $\Gamma _{K_0}$-module such
that $V_0'|_{\Gamma _{K}}=V$ and $V'_0=H_0(\bar K)$ with 
$H_0\in \Gr _{O_0}$,  then 
$V_0$ and $V_0'$ coincide as $\Gamma _{K_0}$-modules.

Suppose the group morphisms 
$\xi :\Gamma _{K_0}\longrightarrow \Aut _{\Bbb F_p}V$ and 
$\xi ':\Gamma _{K_0}\longrightarrow\Aut _{\Bbb F_p}V$ give the structures 
of $\Gamma
_{K_0}$-modules $V_0$ and, respectively, $V_0'$ on $V$. Notice that 
$\xi |_{\Gamma _{K}}=\xi '|_{\Gamma _{K}}$.  
By Fontaine's estimates [Fo] for $e^*=ep/(p-1)$ and any 
$v>e^*-1$, $\xi (\Gamma _{K_0}^{(v)})=\xi '(\Gamma
_{K_0}^{(v)})=
\id _{V}$, where $\Gamma _{K_0}^{(v)}$ is the ramification subgroup of 
$\Gamma _{K_0}$ in the upper numbering. Therefore, 
$\xi $ and $\xi '$ factor through the natural projection 
$\Gamma _{K_0}\longrightarrow \Gamma _{K_0}/\Gamma _{K_0}^{(e^*)}$. 
Notice that $\Gamma _{K_0}^{(e^*)}$ acts non-trivially on $K$. 
(More precisely, $\Gamma _{K_0}^{(v)}$ acts trivially on $K$ if and only
if $v>e^*$.) This implies that $\Gamma _{K}\Gamma
_{K_0}^{(e^*)}=\Gamma _{K_0}$ 
(use that $(\Gamma _{K_0}:\Gamma _{K})=p$). So, the natural 
embedding $\Gamma _{K}\subset\Gamma _{K_0}$ induces the group epimorphism 
$\Gamma _{K}\longrightarrow \Gamma _{K_0}/\Gamma
_{K_0}^{(e^*)}$. Therefore, 
the coincidence of $\xi $ and $\xi '$ on $\Gamma _{K}$ implies that
$\xi =\xi '$. 

The proposition is proved.
$\square $ 
\enddemo 
\medskip 

4.4. Suppose $G\in\Gr _{O}$ and $G'=G\otimes _{O}O'\in\Gr
_{O'}$. 

\proclaim{Proposition 4.4.1} $G$ admits a descent to $O_0$ if and only if
$G'$ admits a descent to $O'_0$. 
\endproclaim 

\demo{Proof} It will be sufficient to consider the cases where
$K'_0/K_0$ is unramified and totally tamely ramified of degree $n$,
where $K_0$ contains a primitive $n$-th root of unity. 

Clearly, the existence of $G_0$ such that $G=G_0\otimes _{O_0}O$ 
implies the existence of $G'_0$ such that 
$G'=G'_0\otimes _{O'_0}O'$ in both cases. 

Now suppose that $G'_0$ exists. Because $G'=G\otimes _{O}O'$,
for any $\tau\in\Gamma $, there is a $\tau $-linear automorphism of 
bialgebras $f'_{\tau }:A(G')\longrightarrow A(G')$ such that 
$\forall\tau _1,\tau _2\in\Gamma $, $f'_{\tau _1\tau
_2}=f'_{\tau _1}f'_{\tau _2}$ and in the case of totally ramified
$K'/K$, 
$\forall\tau\in\Gamma $, $f'_{\tau
}\operatorname{mod}\pi '=\id _{A(G')\otimes k'}$. 
\medskip 

By the uniqueness property from Proposition 4.3.1 for any $\tau\in\Gamma
$, $f'_{\tau }|_{A(G'_0)}=f_{\tau }$ are $\tau $-linear automorphisms of
the coalgebra $A(G'_0)$ and they satisfy the similar properties:  
$\forall\tau _1,\tau _2\in\Gamma $, $f_{\tau _1\tau
_2}=f_{\tau _1}f_{\tau _2}$ and  in the case of totally ramified
$K'_0/K_0$,  
$\forall\tau\in\Gamma $, $f_{\tau
}\operatorname{mod}\pi _0'=\id _{A(G'_0)\otimes k'}$ (because the
embedding $A(G'_0)\subset A(G')$ induces the identity map on
reductions). 
Therefore, $G'_0$ admits a descent $G_0$ to $O_0$ and one can easily see
that $G_0\otimes _{O_0}O=G$. 

The proposition is proved. $\square $
\enddemo

\subsubhead{\rm 4.5.} The Lubin-Tate group law 
\endsubsubhead 

For any $p$-adic ring $R$ denote by $\m (R)$ its topological nilradical,
i.e. the ideal of all $r\in R$ such that $\dsize\lim_{n\to\infty
}r^n=0$. Let 
$$l_{\LT }(X)=X+\frac{X^p}{p}+\dots +\frac{X^{p^n}}{p^n}+\dots \in\Bbb
Q_p[[X]]$$
be the Lubin-Tate logarithm. If $R$ has no $p$-torsion then $\m (R)$ can be
provided with the Lubin-Tate structure of abelian group 
such that for any $f,g\in\m (R)$, 
$[f]+[g]=l_{\LT }^{-1}(l_{\LT }(f)+l_{\LT }(g))$, cf. eg. [Ha].  

We need the following simple properties: 
\medskip 

4.5.1) for any $f,g\in\m (R)$, 
$[f]+[g]=[f+g]+[\phi
_1(f,g)]+\dots +[\phi _n(f,g)]+\dots $, 
where for all $n\geqslant 1$, $\phi _n(f,g)\in\Bbb Z_p[X,Y]$ is a
homogeneous polynomial of degree $p^n$ and, in particular, $\phi
_1(X,Y)=\phi (X,Y)=(X^p+Y^p-(X+Y)^p)/p$;
\medskip 

4.5.2) for any $f\in\m (R)$, 
$[p](f)=[pf]+[\alpha _1f^p]+\dots +[\alpha _nf^{p^n}]+\dots $, 
where all $\alpha _n\in\Bbb Z_p$ and, in particular, $\alpha
_1=1-p^{p-1}$ and for $n\geqslant 2$, $\alpha _n\equiv
0\operatorname{mod}(p^{p-1})$.  
\medskip 

4.5.3) if $X\in\m (R)$ then 
$[p](X)\equiv [X^p]+[pX]\operatorname{mod}p^2R$  
(remind that $p>2$); 
\medskip 

4.5.4) the correspondence $a\mapsto [p](a)$ induces a
one-to-one additive map from $pR$ to $p^2R$.
\medskip

4.6. 
 From above nn.4.2-4.4 it follows  that when studying the image of 
the functor $\Cal G_O$ we can make any tamely ramified extension of
scalars. In particular, we can assume that $O_0$ contains a primitive
$p$-th root of unity $\zeta _p$.  
When proving that $G=\Cal G_O(\Cal M)$, where $\Cal M\in\MF
^e_{S}$, is obtained via an extension of scalars from 
$G_0\in\Gr _{O_0}$, it will be convenient to verify this only on the level of
augmented algebras because of the following property.  

\proclaim{Proposition 4.6.1} 
Suppose $\zeta _p\in O_0$, $\Cal M\in\MF _S^e$, $G=\Spec A=\Cal G_O(\Cal M)$ 
and $(A,I_{A})$
is the corresponding augmented $O$-algebra. If there is an 
$(A_0,I_{A_0})\in\Aug _{O_0}$ such that $I_{A_0}\otimes _{O_0}O=I_{A}$ then $\Spec
A_0$ is provided with a unique structure of $G_0\in\Gr _{O_0}$ such
that  $G_0\otimes _{O_0}O=G$. 
\endproclaim 

\demo{Proof} Suppose $b _1,\dots ,b _u$ is an $O_0$-basis of
$I_{A_0}$. Then it can be considered also as an $O$-basis of $I_{A}$. 
Let $\Delta :A\longrightarrow A\otimes _{O}A$ be the
comultiplication. Then for $1\leqslant i\leqslant u$, 

$$\Delta (b _i)=b _i\otimes 1+1\otimes b _i+\sum\Sb j,r
\endSb a^{(i)}_{jr}b _j\otimes b _r,$$
where all coefficients $a^{(i)}_{jr}\in O$. This map $\Delta $ 
will induce a group
structure on $\Spec A_0$ if and only if all coefficients 
$a^{(i)}_{jr}$ belong to $O_0$. 

Suppose $\tau\in\Gamma =\Gal (K/K_0)$. Then 
$$\Delta ^{(\tau )}:b _i\mapsto \sum\Sb j,r \endSb \tau (a^{(i)}_{jr})b
_j\otimes b _r$$
gives a conjugate group scheme $G^{(\tau )}\in\Gr _{O}$. If
$f_{\tau }$ is a $\tau $-linear automorphism of $A$ given by the
action of $\tau $ on $O$ and the trivial action on $A_0$ 
with respect to  the decomposition $A=A_0\otimes _{O_0}O$ 
then $\Delta ^{(\tau )}=
f^{-1}_{\tau }\circ\Delta  \circ (f_{\tau }\otimes f_{\tau })$. 

Suppose $\Cal M=(M^0,M^1,\varphi _1)\in\MF _S^e$ and 
$G$ is constructed via a special basis $m_1^1,\dots ,m^1_u$
of the $S$-module $M^1$, cf.n.2. If $X_1,\dots ,X_u$ are the variables
attached to the elements 
$m_i=\varphi _1(m_i^1)\in M^0$, $1\leqslant i\leqslant u$, 
then $\Delta $ appears as a unique
$O$-algebra morphism $A\longrightarrow A\otimes _{O}A$ such
that for all $i$, $X_i\mapsto X_i\otimes
1+1\otimes X_i\operatorname{mod}I^{DP}_{A\otimes A}$. Notice that
for 
any $\tau\in\Gamma $ and $1\leqslant i\leqslant u$, 
$f_{\tau }(X_i)\equiv X_i\operatorname{mod}I^{DP}_{A\otimes A}$. 
(Use that for any $c\in I_{A_0}$ and $n\geqslant 1$, 
$(\tau (\pi ^n)-\pi ^n)c
\in (\zeta _p-1)\pi I_{A}\subset 
I_{A}^{DP}$.) Therefore, $\Delta ^{(\tau )}$ appears also as a
unique morphism of $O$-algebras $A\longrightarrow A\otimes
_{O}A$ such that for all $i$, 
$X_i\mapsto X_i\otimes 1+1\otimes X_i\operatorname{mod}I_{A\otimes
A}^{DP}$. 
So, $\Delta ^{(\tau )}=\Delta $ and the proposition is proved. 
$\square $
\enddemo 
\medskip 

\proclaim{Proposition 4.6.2} Suppose $G\in\Cal G_O(\Cal M)$, $\Cal M\in\MF
_S^e$. Then there is $G_0\in\Gr _{O_0}$ such that $G_0\otimes _{O_0}O=G$. 
\endproclaim 

\demo{Proof} By proposition 4.6.1 it will be sufficient to show that the
augmented $O$-algebra $A=A(G)$ comes from an augmented
$O_0$-algebra $A_0$ via the extension of scalars from $O_0$ to $O$.

Suppose $|G|=p$. Then by the results of n.3,  
$A(G)=O[X]$, where $X^p-\eta cX=0$ with $c\in O^*$ 
and $\eta\in O_0$ such that $\eta |p$, cf. Remark 3.2.1. Suppose 
$c=[\alpha ]c_1$, where $[\alpha ]$ is the Teichmuller representative 
of $\alpha\in k$ and $c_1\in O^*$, $c_1\equiv 1\operatorname{mod}\pi
$. 
Let $c_2 =c_1^{1/(p-1)}\in O^*$, $c_2\equiv
1\operatorname{mod}\pi $. Then 
$$(c_2 X)^p-\eta [\alpha ](c_2 X)=0$$
and for the augmentation ideal $I_{A_0}=(c_2X)A_0$ of the 
$O_0$-algebra $A_0=O_0[c_2 X]$, we have 
$I_{A_0}\otimes _{O_0}O=I_{A(G)}$. This proves our proposition if $|G|=p$. 

Suppose $|G|>p$.

By nn.4.2-4.4, 
we can replace $K_0$ by its sufficiently
large tamely ramified extension. Therefore, we can assume that:

$\bullet $\ \ there is a simplest object $\Cal M_{\tilde s}\in\MF _S^e$, 
where $\tilde s\in S$ is an integral power of $t$, $\tilde s|t^e$, 
and there is an $\Cal N=(N^0,N^1,\varphi _1)\in\MF _S^e$
such that $\Cal M=(M^0,M^1,\varphi _1)\in\Ext _{\MF _S^e}(\Cal
M_{\tilde s},\Cal N)$; 

$\bullet $ \ \ if $\tilde\eta '\in O$ is such that 
$\kappa _{SO}(\tilde s\operatorname{mod}t^{ep})=\tilde\eta '\operatorname{mod}p$ 
then $\tilde\eta ^{\prime
p}=\tilde\eta \in O_0$;  in addition, there is a 
$\tilde\lambda '\in O$ such that $\tilde\eta '=\tilde\lambda ^{\prime
p-1}$ and $\tilde\lambda =\tilde\lambda ^{\prime p}\in O_0$ (use again
Remark 3.2.1).
\medskip

Describe the structures of $\Cal M$ and $\Cal N$. There is 
an $S$-basis $m^1,n_1^1,\dots ,n_u^1$ of $M^1$ such that
$N^1=\sum\Sb i\endSb Sn_i^1$. There is an $S$-basis 
$m,n_1,\dots ,n_u$ of $M^0$ such that $N^0=\sum\Sb i\endSb Sn_i$. We
can assume that for $1\leqslant i\leqslant u$, there are 
$\tilde s_i\in S$ such that $\tilde s_i|t^e$, 
$n_i^1=\tilde s_in_i$ and $m^1=\tilde sm+\sum\Sb i\endSb\alpha _in_i$,
where $\alpha _1,\dots ,\alpha _u\in S$ and $t^e\tilde
s^{-1}\alpha _i\equiv 0\operatorname{mod}\tilde s_i$ (or,
equivalently, $t^e\tilde s_i^{-1}\alpha _i\equiv
0\operatorname{mod}\tilde s$). 
The structural morphism $\varphi _1$ is given via the relations 
$\varphi _1(m^1)=m$ and for $1\leqslant i\leqslant u$, 
$\varphi _1(n_i^1)=\sum\Sb j\endSb n_ju_{ji}$, where $(u_{ij})\in\GL
_u(S)$. 

Describe the structure of the corresponding $O$-algebra $B\in\Cal
A(\Cal N)$. It equals $B=O[X_1,\dots ,X_u]$, where for $1\leqslant
i\leqslant u$, 
$X_i^p-\eta _i\sum\Sb j\endSb X_jc_{ji}=0$ with 
$\eta _i=-p/\tilde\eta _i^{\prime p}\in O_0$, $\tilde\eta
_i'\operatorname{mod}p=\kappa _{SO}(\tilde s_i\operatorname{mod}t_1^{ep})$ and 
$C=(c_{ij})\in\GL _u(O)$ is such that 
$C\operatorname{mod}p=\kappa _{SO}(u_{ij}\operatorname{mod}t^{ep})$. 

The structure of the corresponding $A\in\Cal A(\Cal M)$ is given by 
$A=B[Y]$ where 
$$\tilde\eta ^{-1}\left ((Y\tilde\eta '+\sum\Sb i\endSb
r_iX_i)^p+pY\right )=\left (Y+\tilde\eta ^{\prime -1}\sum\Sb i\endSb
r_iX_i\right )^p-\eta Y=0, \tag{4.6.3}$$
with $\tilde\eta '\operatorname{mod}p=\tilde
s\operatorname{mod}t^{ep}$, $\tilde\eta ^{\prime p}\in O_0$, 
$\eta =-p/\tilde\eta ^{\prime p}\in O_0$ and for $1\leqslant i\leqslant
u$, $r_i\in O$, $r_i\operatorname{mod}p=\alpha
_i\operatorname{mod}t^{ep}$ 
and $\eta r_i^p\equiv 0\operatorname{mod}(\tilde\eta _i)$ 
(or, equivalently, $\eta _ir_i^p\equiv 0\operatorname{mod}\tilde\eta
$). Notice that if $\tilde\eta\in O^*$ then by 
Remark 1.3.3 we can assume that 
\linebreak 
$\sum\Sb i\endSb \alpha _in_i\in M^1+tM^0$; so, in this case we can assume
that   
$\sum\Sb i\endSb r_iX_i\in \m (B)$ and, because of equation (4.6.3)
this implies that $Y\in\m (A)$. (Rermind that $\m (A)$ and $\m (B)$
are the topological nilradicals of $A$ and $B$, respectively. ) 

Set $h=\sum\Sb i\endSb r_iX_i$. Then $h\in I_{B}$ and the above
congruences imply that $h^p\in\tilde\eta I_{B}$. By the inductive 
assumption there is an augmented $O_0$-algebra $B_0$ such that
$I_{B}=I_{B_0}\otimes _{O_0}O$. Therefore, for $0\leqslant l<p$, there are 
$b_l\in I_{B_0}$ such that $h=\sum\Sb l\endSb \pi ^lb_l$. 

Let $Y'\in A\otimes _{O}K$ be such that (we use the Lubin-Tate group
law) 
$$[\tilde\lambda Y']=[\tilde\lambda Y+\tilde\lambda 'h]-\sum\Sb
l\endSb [\tilde\lambda '\pi ^lb_l]. \tag{4.6.4}$$

$\bullet $\  If $\tilde\lambda\notin O^*$ then (4.6.4) implies that 
$Y'\equiv Y\operatorname{mod}
\pi I_{A}$. 
\medskip 

$\bullet $\ If $\tilde\lambda \in O^*$ then 
$Y\in\m (A)$, $h\in\m (B)$ and 
$Y'\equiv Y\operatorname{mod} (Y\m (A)+\pi \m (B))$. 
\medskip 

(Cf. the definition of $\m (A)$ and $\m (B)$ in the beginning of
n.4.5.) 
So, in both above cases $A=B[Y']$. 

Find the equation for $Y'$. 

Multiplying (4.6.3) by $\tilde\lambda ^p$ we obtain that 
$(\tilde\lambda Y+\tilde\lambda 'h)^p+p\tilde\lambda Y=0$. 
Using properties of the Lubin-Tate group law from n.4.5 we 
can rewrite this relation as 
$$[p](\tilde\lambda Y')=[p](\tilde\lambda Y+\tilde\lambda 'h)-\sum
_l[p](\tilde\lambda '\pi ^lb_l)\equiv $$
$$[(\tilde\lambda Y+\tilde\lambda 'h)^p]
+[p(\tilde\lambda Y+\tilde\lambda 'h)]-\sum _l[\tilde\lambda
^{\prime p}\pi _0^lb_l^p]-\sum _l[p\tilde\lambda '\pi ^lb_l]\equiv $$
$$-\sum\Sb 0\leqslant l<p\endSb
[\tilde\lambda \pi _0^lb^p_l]\operatorname{mod}p^2I_{A}.$$
Therefore, by replacing $Y'$ by $Y_1=Y'-(p/\tilde\lambda )a$ with a
suitable $a\in I_{A}$ we still have $A=B[Y_1]$ and $Y_1$ will 
satisfy the following relation  
$$[p](\tilde\lambda Y_1)=-\sum\Sb 0\leqslant l<p\endSb 
[\tilde\lambda\pi _0^lb_l^p]. \tag{4.6.5}$$ 
Notice that the right hand side of (4.6.5) equals $\tilde\lambda ^pb_0$, 
$b_0\in I_{B^0}$ (use that $\sum _l\pi _0^lb_l^p\equiv
h^p\operatorname{mod}p$ and $h^p\in\tilde\eta I_{B}$). If 
$E=\exp (X+X^p/p+\dots )$ is the Artin-Hasse exponential then 
$E(\tilde\lambda Y_1)=1+\tilde\lambda Y_2$ and 
we still have $A=B[Y_2]$ (this is obvious if 
$\tilde\lambda\notin O_0^*$ and use that $Y_1\in\m (A)$, 
otherwise). 
If $E(\tilde\lambda ^pb_0)=1+\tilde\lambda ^pb_0'$ then 
$b_0'\in I_{B_0}$ and $Y_2$ is a root of the unitary polynomial
$$F=\frac{(1+\tilde\lambda T)^p-1}{\tilde\lambda ^p}-b_0'\in B_0[T].$$
This implies (use that $\rk _{B}A=p$ and $A=B[Y_2]$) that 
$A=B[T]/(F)$. Therefore, for the augmented algebra 
$A_0=B_0[T]/(F)$ we have $I_{A}=I_{A_0}\otimes _{O_0}O$. 

The proposition is proved. 
$\square $
\enddemo

\proclaim{Corollary 4.6.6} There is a fully faithful functor $\Cal G^O_{O_0}:\MF
^e_{S}\longrightarrow\Gr _{O_0}$ such that for any $\Cal M\in\MF
^e_{S}$, $\Cal G^O_{O_0}(\Cal M)\otimes _{O_0}O=\Cal G_O(\Cal
M)$. $\square $ 
\endproclaim

\medskip 
\medskip

\subhead{\rm 5.} Group of classes of short exact sequences in $\Gr _{O}$ 
\endsubhead 
\medskip 

In this section we do not use that the ring $O$ is obtained from the ring $O_0$
by joining a $p$-th root of some uniformizing element of $O_0$. This
will allow us to apply in n.6 the results of this section also to the
category $\Gr _{O_0}$ with $O$ replaced by $O_0$. 

For technical reasons we shall assume here that the residue field $k=\bar k$ is
algebraically closed and there is $\pi ^*\in O$
such that ${\pi ^{*}}^{p-1}=-p$. An element $\eta\in O$ will always be such
that $\eta |p$ and there is an 
$\lambda\in O$ such that $\lambda ^{p-1}=\eta $. We set 
$\tilde\lambda =\pi ^*/\lambda $ and $\tilde\eta =\tilde\lambda
^{p-1}$. In particular, $\eta\tilde\eta =-p$ and $\lambda\tilde\lambda
=\pi ^*$. 

\subsubhead{\rm 5.1.} The different and the trace 
\endsubsubhead 
\medskip

Suppose $B\in\Alg _{O}$, 
i.e. $B$ is a flat finite $O$-algebra, 
and $A$ is a faithfully flat finite $B$-algebra. 

\proclaim{Proposition 5.1.1} Suppose there is $\theta\in A$ such that
$A=B[\theta ]$. Then:  
\medskip  
{\rm a)} there is a unique monic polynomial $F\in B[X]$ such that
$A=B[X]/(F)$ and $\theta =X\operatorname{mod}F$;
\medskip 
{\rm b)} the ideal $(F'(\theta ))$ does not depend on a choice of
$\theta $;
\medskip 
{\rm c)} $A_{K}=A\otimes _{O}K$ is etale over 
$B_{K}$ if and only if $F'(\theta )\in
A^{*}_{K}$. 
\endproclaim 

\demo{Proof} a) follows because for any maximal ideal $\m $ in $B$ 
$\dim _{B/{\m} B}(A/{\m} A)$ does not depend on $\m $.

b) Suppose $A=B[X _1]$ and $\theta _1=X_1\operatorname{mod}F_1$, 
where $F_1(X)\in B[X_1]$ is monic. Let $G(X)\in B[X]$, $H(X_1)\in
B[X_1]$ be such that $\theta _1=G(\theta )$ and $\theta =H(\theta
_1)$. Then 
$H(G(X))\equiv X\operatorname{mod}F(X)$ implies that 
$H'(\theta _1)G'(\theta )\equiv 1\operatorname{mod}F'(\theta )$. 
In addition, $F_1(G(X))\equiv 0\operatorname{mod}F(X)$ implies that 
$F_1'(\theta _1)G'(\theta )
\equiv 0\operatorname{mod}F'(\theta )$. 
Therefore, $F_1'(\theta _1)\in (F'(\theta ))$ and by symmetry 
$F'(\theta )\in (F_1'(\theta ))$. 

c) $A_{K}$ is etale over $B_{K}$ if and only if $\Omega
^1_{A_{K}/B_{K}}=\Omega ^1_{A/B}\otimes _{O}K=0$. 
It remains to notice that $\Omega ^1_{A/B}=
A/(F'(\theta ))dX$. 
$\square $
\enddemo 

\definition{Definition} With the above notation:
\medskip  
a) $\Cal D(A/B)=(F'(\theta ))$ is the {\it different} of $A$ over $B$;
\medskip 
b) if $A_{K}$ is etale over $B_{K}$ then we set $\Cal
D^{-1}(A/B)=F'(\theta )^{-1}A\subset A_{K}$. 
\enddefinition 
\medskip 

\remark{Remark {\rm 5.1.2}} The part a) of the above definition implies that 
the norm $N_{A/B}(F'(\theta ))$ is the discriminant of the $B$-algebra
$A$. 
\endremark 
\medskip 

Let $\Tr _{A/B}:A\longrightarrow B$ be the trace map and let $\Tr _
{A_{K}/B_{K}}=
\Tr _{A/B}\otimes _{O}K:A_{K}\longrightarrow B_{K}$. 
Suppose $F(X)$ splits completely in $A$, i.e. there are 
$\theta _{\alpha }$, $1\leqslant \alpha \leqslant\deg F$ 
such that $F(X)=\prod \Sb\alpha \endSb (X-\theta _{\alpha })$. 
For $1\leqslant \alpha\leqslant \deg F$ introduce the $O$-algebra
morphisms 
$t_{\alpha }:A\longrightarrow A$ 
such that $t_{\alpha }(\theta )=\theta _{\alpha }$ and $t_{\alpha
}|_B=\id $. Clearly, for any $a\in A$, 
$\Tr _{A/B}(a)=\sum\Sb\alpha \endSb t_{\alpha }(a)$. 

\proclaim{Proposition 5.1.3} Suppose $A=B[\theta ]$, $A_{K}$ is etale over
$B_{K}$ and $F(X)$ splits completely over $A$. Then 
there is an $a\in\Cal D^{-1}(A/B)$ such that $\Tr _{A_{K}/B_{K}}(a)=1$. 
\endproclaim 

\demo{Proof} This follows from the case $k=n-1$ of the following
lemma.
\enddemo 

\proclaim{Lemma 5.1.4} If $X_1,\dots ,X_n$ are independent variables
over $\Bbb Q$, $0\leqslant k\leqslant n-1$ 
and $\delta $ is the Kronecker symbol then 

$$\sum _{i=1}^{n}\frac{X_i^k}{\prod\Sb j\ne i
\endSb (X_i-X_j)}=\delta _{k,n-1}.$$
\endproclaim 

\demo{Proof} Consider the decomposition into a sum of simplest
fractions in $L(X_1)$, where $L=\Bbb Q(X_2,\dots ,X_n)$, 
$$\frac{X_1^k}{\prod\Sb j\ne 1\endSb (X_1-X_j)}=\sum _{j=2}^n 
\frac{A_j}{X_1-X_j}+\delta _{k,n-1}.\tag{5.1.5}$$

Then multiplying this identity by 
$\prod _{j\ne 1}(X_1-X_j)$ and substituting 
for $2\leqslant j\leqslant n$, $X_1=X_j$ we obtain  
$$A_j=\frac{X_j^k}{\prod\Sb s\ne 1,j\endSb (X_j-X_s)}.$$
It remains to substitute these formulas to (5.1.5). 
The lemma is proved. 
$\square $
\enddemo 

\remark{Remark {\rm 5.1.6}} The above lemma implies that 
$$\Cal D^{-1}(A/B)=\{a\in A_{K}\ |\  
\Tr _{A_{K}/B_{K}}(a)\in B\}.$$
\endremark

\subsubhead{\rm 5.2.} Group schemes $\G _{\tilde\lambda }$ 
\endsubsubhead 
\medskip 

Suppose $G\in\Gr _{O}$ is of order $p$. Then by [TO, p.14, Remark] 
there are $\eta ,\tilde\eta\in O$ 
such that $\eta\tilde\eta =-p$ and $G\simeq \G_{\tilde\lambda }$, where 
$A(\G_{\tilde\lambda })=O[Y_{\tilde\lambda }]$ with $Y_{\tilde\lambda
}^p-
\eta Y_{\tilde\lambda }=0$. The counit
$e_{\G_{\tilde\lambda }}:A(\G_{\tilde\lambda })\longrightarrow O$ and the comultiplication
$\Delta _{\G_{\tilde\lambda }}
:A(\G_{\tilde\lambda })\longrightarrow A(\G_{\tilde\lambda })
\otimes _OA(\G_{\tilde\lambda })$ are uniquely determined 
by the conditions $e_{\G_{\tilde\lambda }}(Y_{\tilde\lambda })=0$ and 
$\Delta _{\G_{\tilde\lambda }}(Y_{\tilde\lambda })
=Y_{\tilde\lambda }\otimes 1+1\otimes
Y_{\tilde\lambda }+\tilde\eta\phi (Y_{\tilde\lambda })\operatorname{mod}p\tilde\eta $.  
Notice, if $\tilde\eta _1, \tilde\lambda _1\in O$, $\tilde\eta _1|p$ 
and $\tilde\eta _1=\tilde\lambda _1^{p-1}$ then 
$\G_{\tilde\lambda }$ is isomorphic to $\G_{\tilde\lambda _1}$ iff 
$\tilde\eta \tilde\eta _1^{-1}\in O^{*}$
(remind that $k=\bar k$). 

Notice that:
\medskip 

a) if 
$\eta\in O^{*}$ then $\G_{\tilde\lambda }$ is etale; in particular, if $\eta =1$
then $\G_{\tilde\lambda }=\G _{\pi ^*}$ is constant etale; 
\medskip 

b) if
$\tilde\eta\in O^{*}$ then $\G_{\tilde\lambda }$ is multiplicative. In particular, 
$\G_1$ is isomorphic to the constant multiplicative 
group scheme $\mu _p$ 
of order $p$ given by the $O$-algebra $A(\mu _p)=O[T]$, where $T^p=1$, $e(T)=1$ and
$\Delta (T)=T\otimes T$. This implies the existence of a polynomial
$P\in O[X]$ such that $P(X)\equiv X\operatorname{mod}X^2$, 
$(1+P(Y_1))^p=1$ and $\Delta _{\G_1}(1+P(Y_1))=(1+P(Y_1))\otimes
(1+P(Y_1))$. 
\medskip 

c) there is a natural morphism of group schemes $\delta
_{\tilde\lambda }:\G_{\tilde\lambda }\longrightarrow \G_1$ given by the $O$-algebra
morphism 
$\delta ^*_{\tilde\lambda }:O[Y_1]\longrightarrow O[Y_{\tilde\lambda }]$ such that 
$\delta ^*_{\tilde\lambda }(Y_1)=\tilde\lambda Y_{\tilde\lambda }$. If we use the above
identification 
$G_1\simeq \mu _p$ then the corresponding morphism $\delta _{\tilde\lambda }:\G_{\tilde\lambda }
\longrightarrow \mu _p$ is given by the correspondence 
$T\mapsto 1+P(\tilde\lambda Y_{\tilde\lambda })$ with $P\in O[X]$ from above n.b);
\medskip 

d) the set of all geometric points of $\G_{\tilde\lambda }$ equals 
$\G_{\tilde\lambda }(O)=\{g_{\alpha }\ |\ \alpha\in\Bbb F_p\}$, where 
$g_{\alpha }(Y_{\tilde\lambda })=[\alpha ]\lambda $ 
($[\alpha ]$ is the Teichmuller representative of $\alpha \in\Bbb
F_p$). 
Then for any $\alpha\in\Bbb F_p$ and $\delta _{\tilde\lambda }:\G_{\tilde\lambda
}\longrightarrow \G_1$ from above n.c),  $\delta _{\tilde\lambda }(g_{\alpha })=\zeta
(\alpha )$, where $\zeta (\alpha )\in O$ is the $p$-th root of unity
uniquely determined by the congruence 
$\zeta (\alpha )\equiv 1+[\alpha ]\pi ^*\operatorname{mod}\pi ^*\pi
$.

\subsubhead{\rm 5.3.} $\G_{\tilde\lambda }$-torsors 
\endsubsubhead 

Let 
$B\in\Alg _{O}$. Then a 
$\G_{\tilde\lambda }$-torsor over $B$ is a finite faithfully flat
$B$-algebra $A\in\Alg _{O}$ with the action of $\G_{\tilde\lambda }$ given by an 
$O$-algebra morphism 
$\omega :A\longrightarrow A(\G_{\tilde\lambda })\otimes _{O}A$ such that 
\medskip 

$\bullet $\ \ $\omega \circ (\id \otimes\omega )=
\omega \circ (\Delta _{\G_{\tilde\lambda }}\otimes \id )$; 
\medskip 

$\bullet $\ \ $B=A^{\G_{\tilde\lambda }}=\{a\in A\ |\ \omega (a)=1\otimes a\}$; 
\medskip 

$\bullet $\ \ the correspondence $a_1\otimes a_2\mapsto \omega
(a_1)(1\otimes a_2)$ induces an identification of $O$-algebras 
$A\otimes _BA=A(\G_{\tilde\lambda })\otimes _{O}A$. 
\medskip 

Suppose $A_1$ is another $\G_{\tilde\lambda }$-torsor over $B$ with the $\G_{\tilde\lambda
}$-action given by $\omega _1:A_1\longrightarrow A(\G_{\tilde\lambda })\otimes
A_1$. Then $A$ and $A_1$ are equivalent if there is an isomorphism of
$B$-algebras $\nu :A\longrightarrow A_1$ such that 
$\nu\circ \omega _1=\omega \circ (\id\otimes\nu )$. 

The set $E(\G_{\tilde\lambda }, B)$ 
of equivalence classes of $\G_{\tilde\lambda }$-torsors over $B$ has a natural
structure of abelian group given by the Baer composition $*$. Remind
that  
$A*A_1=(A\otimes _{B}A_1)^{\G_{\tilde\lambda }}$ where $\G_{\tilde\lambda }$ acts on 
$A\otimes _BA_1$ via the composition of the 
antidiagonal embedding into $\G_{\tilde\lambda }\times
\G_{\tilde\lambda }$ and the component-wise action $\omega \otimes\omega _1$ 
of $\G_{\tilde\lambda }\times \G_{\tilde\lambda }$
on $A\otimes _BA_1$. 
\medskip 

\subsubhead{\rm 5.4.} Construction of $\G_{\tilde\lambda }$-torsors 
\endsubsubhead 

As earlier, $B\in\Alg _{O}$ and $\m (B)$ is the topological nilradical
of $B$.

Denote by $\hat\Bbb G_{m,\tilde\lambda }$ the formal group functor
such that if $B$ is an $O$-algebra then 
$\hat\Bbb G_{m,\tilde\lambda }(B)=(1+(\tilde\lambda B\cap \m (B))^{\times }$.  
If $\tilde\lambda =1$ we shall use also the usual notation 
$\hat\Bbb G_m$ for $\hat\Bbb G_{m,1}$. 

Suppose $1+\tilde\lambda ^pb\in\hat\Bbb G_{m,\tilde\lambda ^p}(B)$.   
Let $A$ be the quotient of the polynomial ring 
$B[X]$ by the ideal generated by the monic polynomial 
$F_b(X)=\tilde\lambda ^{-p}((1+\tilde\lambda X)^p-1)-b$. Denote the image of
$X$ in 
$A$ by $\theta _b$, then $A=B[\theta _b]$   
is a faithfully flat $B$-algebra. 

For $\alpha\in\Bbb F_p$, there is a unique $O$-algebra isomorphism 
$t_{\alpha }:A\longrightarrow A$ such that $t_{\alpha }|_B=\id $ and 
$$t_{\alpha }:1+\tilde\lambda\theta _b\mapsto (1+\tilde\lambda\theta
_b)\zeta (\alpha ), \tag {5.4.1}$$
where $\zeta (\alpha )$ is the $p$-th root of unity such that 
$\zeta (\alpha )\equiv 1+[\alpha ]\pi ^*\operatorname{mod}\pi ^*\pi
$. 
Indeed, the correspondence (5.4.1) determines a unique $K$-algebra
automorphism $t_{\alpha K}:A_{K}\longrightarrow A_{K}$, where $A_{K}=A\otimes
_{O}K$, and clearly $t_{\alpha }(A)\subset A$ (use that $\tilde\lambda
|\pi ^*$). 

\proclaim{Proposition 5.4.2} {\rm a)}\ There is a unique action of $\G_{\tilde\lambda }$ on $A$
given by the $O$-algebra homomorphism 
$\omega :A\longrightarrow A(\G_{\tilde\lambda })\otimes _{O}A$ 
such that for any $\alpha\in\Bbb F_p$, $\omega \circ (g_{\alpha }\otimes
1)=t_{\alpha }$, where $g_{\alpha }\in \G_{\tilde\lambda }(O)$ 
were defined in 5.2 d), and this action 
determines on $A$ a structure of $\G_{\tilde\lambda }$-torsor over $B$;
\medskip 

{\rm b)} the correspondence $1+\tilde\lambda ^pb\mapsto A=B[\theta
_b]$ 
determines a group epimomorphism 
$\kappa :\hat\Bbb G_{m,\tilde\lambda ^p}(B) 
\longrightarrow E(\G_{\tilde\lambda } ,B)$ 
and $\Ker\kappa =\hat\Bbb G_{m,\tilde\lambda }(B)^p$.  
\endproclaim 

\demo{Proof} Notice that for any $\alpha\in\Bbb F_p$, 
$g_{\alpha }(1+P(\tilde\lambda Y_{\tilde\lambda }))=\zeta (\alpha
)$. Therefore, the only candidate for such action of $\G_{\tilde\lambda }$ must
satisfy the following requirement  
$$\omega :1+\tilde\lambda\theta _b\mapsto (1+P(\tilde\lambda Y_{\tilde\lambda
}))\otimes (1+\tilde\lambda \theta _b).\tag{5.4.3}$$
Clearly, this requirement determines a unique $O$-algebra homomorphism 
$\omega :A\longrightarrow A(\G_{\tilde\lambda })\otimes _{O}A$. 

Let $\G_{\tilde\lambda ,K}=\G_{\tilde\lambda }\otimes _OK$. 
Prove that $\omega _K:=\omega\otimes K$ defines a $\G_{\tilde\lambda ,K}$-torsor
over $B_{K}=B\otimes _{O}K$.

Notice that $\delta _{\tilde\lambda }\otimes K$ determines the identification 
$\G_{\tilde\lambda ,K}=\mu _{p,K}=\mu _p\otimes K$. Then 
$\omega _{K}(1+\tilde\lambda\theta _b)=
T\otimes (1+\tilde\lambda \theta _b)$, 
where $A(\mu _{p,K})=K(T)$ with the comultiplication 
$\Delta (T)=T\otimes T$. Therefore, 
$\omega _{K}\otimes (\id\otimes\omega _{K})=\omega _{K}\otimes 
(\Delta \otimes\id )$, i.e. 
$\G_{\tilde\lambda ,K}$ acts on $A_{K}$. 
Similarly,  $A_{K}^{\mu _{p,K}}=B_{K}$, and the
correspondence $a_1\otimes a_2\mapsto\omega _{K}(a_1)(1\otimes a_2)$ gives an
isomorphism $\xi _{K}$ of $K$-algebras 
$A_{K}\otimes _{B_{K}}A_{K}$ and $A(\mu _{p,K})\otimes _{K}A_{K}$. 

Therefore, $\omega \circ (\id\otimes\omega )=
\omega \circ (\Delta _{\G_{\tilde\lambda
}}\otimes\id )$, because 
$A\subset A_{K}$ amd $\omega _{K}$ maps $A$ into $A(\G_{\tilde\lambda })\otimes A$.
We have also that 
$A^{\G_{\tilde\lambda }}=(A_{K}^{\mu _{p,K}})\cap A=B_{K}\cap A=B$. 
Finally, $\xi _{K}$ induces an embedding of $O$-algebras 
$\xi :A\otimes _BA\longrightarrow A(\G_{\tilde\lambda })\otimes _{O}A$ 
such that $\xi |_{1\otimes A}=\id $. Now notice that 
$$\Cal D(A\otimes _BA/1\otimes A)=\Cal D(A/B)=(\eta )=
\Cal D(A(\G_{\tilde\lambda })/O)=\Cal D(A(\G_{\tilde\lambda })\otimes _{O}A/1\otimes A).$$
(Use that $A$ is a faithfully flat $B$- and $O$-algebra, 
$F'_b(\theta _b)=\eta (1+\tilde\lambda\theta _b)^{p-1}$ and 
$(1+\tilde\lambda \theta _b)^p=1+\tilde\lambda ^pb\in B^*$.) 
Therefore, the discriminants of both $(1\otimes A)$-algebras, 
$A\otimes _BA$ and $A(\G_{\tilde\lambda })\otimes _{O}A$ are equal, 
and the embedding $\xi $ is an isomorphism. The part a) is proved.

One can see easily that the map $\kappa $ from the part b) of
our proposition is a group homomorphism and its kernel 
is $\hat\Bbb G_{m,\tilde\lambda }(B)^p$. It remains to 
prove that $\kappa $ is epimorphic.

Suppose $A\in E(\G_{\tilde\eta },B)$. 

First, prove that we can use the concept of the differente for
the $B$-algebra $A$. We need the following properties:
\medskip 

1) there is an $\theta\in A$ such that $A=B[\theta ]$;
\medskip 

2) $A_{K}$ is etale over $B_{K}$;
\medskip 

3) $\Cal D(A/B)=(\eta )$ and, therefore, is an invertible ideal of $A$
in $A_{K}$.  
\medskip 

Indeed, we know that $B_1=B/\m (B)$ is a product of finitely many
copies of $k$. Therefore, the $B_1$-algebra $A_1=A/\m (B)A$ can be
provided with augmentation (use that $k$ is algebraicly closed). 
This implies that  
$E(\G_{\tilde\lambda }\otimes k,B_1)=0$ and $A_1=A(\G_{\tilde\lambda }\otimes
B_1)=A(\G_{\tilde\lambda })\otimes B_1=B_1[Y_{\tilde\lambda }]$. So, by the Nakayama
Lemma, $A=B[\theta ]$, where $\theta\in A$ is such that
$\theta\operatorname{mod}\m (B)A=Y_{\tilde\lambda }$. 
The identification $A\otimes _BA=A(\G_{\tilde\lambda })\otimes _{O}A$ implies that 
$A_{K}\otimes _{B_{K}}A_{K}$ is etale $1\otimes A_{K}$-algebra (because
$A(\G_{\tilde\lambda })\otimes K$ is etale over $K$) and by faithful flatness
$A_{K}$ is etale over $B_{K}$. 
Finally, 
$\Cal D(A/B)\otimes _BA=\Cal D(A(\G_{\tilde\lambda })/O) \otimes _{O}A=(\eta
A)\otimes _BA$ 
implies by faithful flatness that $\Cal D(A/B)=(\eta )$. 

Next, prove the existence of 
$v\in\hat\Bbb G_{m,\tilde\lambda }(A)$ such
that 
$\omega (v)=T\otimes v$, where 
$T=1+P(\tilde\lambda Y_{\tilde\eta })\in A(\G_{\tilde\lambda })$, cf. 5.2. 
\medskip 

{\it 1st case}. $\tilde\lambda\in O^{*}$.
\medskip 

In this case $A(\G_{\tilde\lambda })\simeq A(\mu _p)=O[T]$, $e(T)=1$ and 
$\Delta (T)=T\otimes T$. We know that 
$A\operatorname{mod}\m (B)A=O[T]\otimes _{O}B_1$, where $B_1=B/\m (B)$. 
Let $\theta\in A$ be such that $\theta\operatorname{mod}\m
(B)A=T\otimes 1$,
then $\theta\equiv 1\operatorname{mod}\m (A)$ and $\omega (\theta )
\equiv T\otimes\theta \operatorname{mod}\m (B)A$. Therefore, if 
$\dsize\omega (\theta )=\sum\Sb 0\leqslant i<p\endSb T^i\otimes a_i\in
A(\mu _p)\otimes A$ 
then $a_1\equiv \theta \operatorname{mod}\m (B)A$ and, therefore, $a_1\equiv
1\operatorname{mod}\m (A)$. On the other hand, 
$$(\omega \circ (\id\otimes\omega ))(\theta )=\sum\Sb i\endSb
T^i\otimes\omega (a_i)=
(\Delta \otimes\id )(\omega (\theta ))=
\sum\Sb i\endSb T^i\otimes T^i\otimes a_i$$
implies that $\omega (a_1)=T\otimes a_1$ and we can take $v=a_1$. 
\medskip 

{\it 2nd case}. $\tilde\lambda\notin O^*$, i.e. $\tilde\eta =-p\eta ^{-1}\notin
O^*$. 
\medskip 

By Proposition 5.1.3 we can choose $\theta\in A$ such that 
$\Tr _{A/B}\theta =\eta $. Clearly, $\theta\notin B$, otherwise, 
$\eta\equiv 0\operatorname{mod}p$ and $\tilde\eta\in O^*$.  Then
there is $1\leqslant m<p$ such that 
$$v_1=\sum\Sb \alpha\in\Bbb F_p\endSb \zeta (\alpha )^{m}t_{\alpha
}(\theta )\ne 0.$$
Indeed, otherwise, for all $1\leqslant m<p$, 
$\sum\Sb\alpha\endSb \zeta (\alpha )^{m}(t_{\alpha }(\theta )-\theta
)=0$ and this implies that for all $\alpha\in\Bbb F_p$, $t_{\alpha
}(\theta )=\theta $, i.e. $\theta\in B$. 

Let $(1+P(\tilde\lambda Y_{\tilde\lambda }))^m=1+\tilde\lambda h$, $h\in
A(\G_{\tilde\lambda })$. Then for all $\alpha\in\Bbb F_p$, 
$$g_{\alpha }(1+\tilde\lambda h)=\zeta (\alpha )^m=e_{\G_{\tilde\lambda
}}(1+\tilde\lambda t_{\alpha }(h)).$$
(In this situation $t_{\alpha }:A(\G_{\tilde\lambda })
\longrightarrow A(\G_{\tilde\lambda })$ is just
the shift by $g_{\alpha }\in \G_{\tilde\lambda }(O)$.)  Then by Proposition 5.1.3
$$v_1-\eta =\tilde\lambda (e_{\G_{\tilde\lambda }}\otimes\id
)\sum\Sb\alpha \endSb t_{\alpha }(h)\otimes t_{\alpha }(\theta )
=\tilde\lambda (e_{\G_{\tilde\lambda }}\otimes\id )\Tr (h\otimes\theta
)\in\tilde\lambda \eta A,$$ 
where 
$\Tr $ is the trace map on $A(\G_{\tilde\lambda })\otimes A$ induced by the 
diagonal action of $\G_{\tilde\lambda }$. So, 
$v_1=\eta (1+\tilde\lambda a_1)$ with $a_1\in A$ and 
for any $\alpha\in\Bbb F_p$, $t_{\alpha }(1+\tilde\lambda a_1)=\zeta
(\alpha )^m(1+\tilde\lambda a_1)$. 

The required element $v$ then can be obtained by taking $m'$-th power
of $1+\tilde\lambda a_1$, where $mm'\equiv 1\operatorname{mod}p$. 
The second case is also considered.
\medskip 

Finally, for the above constructed element $v$, we have 
$v^p=1+\tilde\lambda ^pb\in\hat\Bbb G_{m,\tilde\lambda ^p}(B)$ and 
the corresponding $\G_{\tilde\lambda }$-torsor 
$\kappa (v^p)\in E(\G_{\tilde\lambda },B)$ is
identified with a $B$-subalgebra $A'$ in $A$. But the 
differentes 
$\Cal D(A/B)$ and $\Cal D(A'/B)$ are both equal to $(\eta )$. 
Therefore, the discriminants of $A$ and $A'$ over $B$ are equal and 
$A=A'$. The proposition is proved.
$\square $
\enddemo

5.5. Suppose $H=\Spec B$ is a finite flat commutative group scheme
over $O$ 
with the counut $e$ and the comultiplication $\Delta $. 
Denote by $\Ext (H,\G_{\tilde\lambda })$ the group of equivalence classes 
of short exact sequences
$0\longrightarrow \G_{\tilde\lambda }\longrightarrow G\longrightarrow
H\longrightarrow 0$ 
in the category $\Gr '_O$ of commutative finite flat group schemes over $O$.  
(Notice that we do not assume that $H$ and $G$ belong to $\Gr _O$, 
i.e. are killed by $p$.)

\definition{Definition} {\rm a)} $Z^2(H,\hat\Bbb G_{m,\tilde\lambda
})$ is the group of all symmetric (with respect to the permutation of
 factors in $B\otimes B$) $\varepsilon \in \hat\Bbb G_{m,\tilde\lambda }
(I_{B\otimes B})$ such that 

$$(\Delta \otimes\id _B)(\varepsilon )\cdot (\varepsilon \otimes 1)=(\id
_B\otimes\Delta )(\varepsilon )\cdot (1\otimes\varepsilon ).$$

{\rm b)} $B^2(H,\hat\Bbb G_{m,\tilde\lambda })$ is the
multiplicative group of all elements of the form 
$\delta ^{\times }a :=\Delta (a)(a\otimes a )^{-1}$, 
where $a\in \hat\Bbb G_{m,\tilde\lambda }(I_B)$. 
\enddefinition 

Then $B^2(H,\hat\Bbb G_{m,\tilde\lambda })$ is a subgroup in 
$Z^2(H,\hat\Bbb G_{m,\tilde\lambda })$ and we set  
$H^2(H,\hat\Bbb G_{m,\tilde\lambda })=Z^2(H,\hat\Bbb G_{m,\tilde\lambda
})/B^2(H,\hat\Bbb G_{m,\tilde\lambda })$. 

If $\varepsilon\in Z^2(H,\hat\Bbb G_{m,\tilde\lambda })$ then 
$\varepsilon ^p\in Z^2(H,\hat\Bbb G_{m,\tilde\lambda ^p})$ and the
correspondence $\varepsilon \mapsto\varepsilon ^p$ induces the group 
homomorphism 
$\Cal F:H^2(H,\hat\Bbb G_{m,\tilde\lambda })\longrightarrow 
H^2(H,\hat\Bbb G_{m,\tilde\lambda ^p})$. 

\definition{Definition} $H^2(H,\hat\Bbb
G_{m,\tilde\lambda })_p:=\Ker \Cal F$.
\enddefinition 

\proclaim{Proposition 5.5.1} There is a functorial in $H$ group
isomorphism 
$$\Ext (H,\G_{\tilde\lambda })=H^2(H,\hat\Bbb G_{m,\tilde\lambda })_p.$$
\endproclaim 

\demo{Proof} Suppose $G\in\Ext (H,\G_{\tilde\lambda })$. Then $A=A(G)$ is
provided with a natural action of $\G_{\tilde\lambda }\subset G$ and 
with respect of this action $A$ becomes 
an element of the group of torsors $E(B,\G_{\tilde\lambda })$. 

By Proposition 5.4.2, $A=B[\theta ]$, where for 
$v=1+\tilde\lambda\theta\in \hat\Bbb G_{m,\tilde\lambda }(I_A)$, 
it holds 
$v^p=1+\tilde\lambda ^pb_0\in\hat\Bbb
G_{m,\tilde\lambda ^p}(I_B)$ (use the 
existence of counit $e_G:A\longrightarrow O$). Then 
$\Delta _G|_B=\Delta $ and $\Delta _G(1+\tilde\lambda\theta )=
(1+\tilde\lambda\theta )\otimes (1+\tilde\lambda\theta )\varepsilon
$, where  $\varepsilon\in\hat\Bbb G_{m,\tilde\lambda }(I_{B\otimes B})$ (use
that $\Delta _G$ relates the actions of $\G_{\tilde\lambda }$ and $\G_{\tilde\lambda
}\times \G_{\tilde\lambda }$ via the composition $\G_{\tilde\lambda }\times \G_{\tilde\lambda
}\longrightarrow \G_{\tilde\lambda }$) and, in addition,   
$\varepsilon\in Z^2(H,\hat\Bbb G_{m,\tilde\lambda })$ 
(use the coassociativity and commutavity of $\Delta _G$ and the
existence of counit $e_G$). 

The above $v\in\hat\Bbb G_{m,\tilde\lambda }(I_A)$ is well-defined
modulo the subgroup $\hat\Bbb G_{m,\tilde\lambda }(I_B)$, this implies
that $G$ depends only on the class $\cl (\varepsilon )\in
H^2(H,\hat\Bbb G_{m,\tilde\lambda })$ of $\varepsilon $. (If we change
$v$ by $va$ with 
$a\in \hat\Bbb G_{m,\tilde\lambda }(I_B)$, then $\varepsilon $ will be
changed by $\varepsilon\delta ^{\times }a$.) Clearly, 
$\varepsilon ^p\in B^2(H,\hat\Bbb G_{m,\tilde\lambda ^p})$ and,
therefore, we have the map 
$\Pi :\Ext (H, \G_{\tilde\lambda })\longrightarrow H^2(H,\hat\Bbb
G_{m,\tilde\lambda })_p$. 
Then a straightforward verification shows that $\Pi $ is a group
isomorphism. $\square $
\enddemo 

Now suppose that $H\in\Gr _O$ and denote by $\Ext _{\Gr _O}(H,\G_{\tilde\lambda
})$ 
the subgroup   
of extensions $G\in\Ext (H,\G_{\tilde\lambda })$ such that $G\in\Gr _O$. 
Consider the map 
$$\delta _{\tilde\lambda *}:\Ext (H,\G_{\tilde\lambda })\longrightarrow\Ext (H,\mu
_p),$$ 
where $\delta _{\tilde\lambda }:\G_{\tilde\lambda }\longrightarrow \mu _p$ is the
morphism from 5.2.1 c). Let $\bar O$ be the valuation ring of an algebraic
closure $\bar K$ of $K$. Clearly, for 
any $G\in\Ext (H,\mu _p)$, 
$$G\in\Gr _O\ \ \Leftrightarrow\ \  
(\delta _{\tilde\lambda *}G)\otimes _O\bar O\in\Ext _{\Gr _{\bar O}}(\bar
H,\bar\mu _p).$$
(Here $\bar H=H\otimes _O\bar O$ and $\bar{\mu }_p=\mu _p\otimes
_O\bar O$). 

\proclaim{Lemma 5.5.2} $\Ext _{\Gr _{\bar O}}(\bar H,\bar\mu _p)=0$. 
\endproclaim 

\demo{Proof} By the Cartier duality we must prove
that 
$\Ext _{\Gr _{\bar O}}((\Bbb Z/p)_{\bar O},\bar H^D)=0$, 
where $\bar H^D$ is the Cartier dual for $\bar H$ and $(\Bbb
Z/p)_{\bar O}$ is the constant group scheme of order $p$ over $\bar
O$. Equivalently, we must prove that in the category $\Gr _{\bar O}$ 
any faithfully flat morphism 
$\gamma :\bar G\longrightarrow (\Bbb Z/p)_{\bar O}$ has a
section $\hat\gamma :(\Bbb Z/p)_{\bar O}\longrightarrow\bar G$. 

Consider the $\bar O$-algebra $C=A\left ((\Bbb Z/p )_{\bar O}\right )
=\oplus _{i\in\Bbb Z/p}\bar O\theta _i$, where for any $i\in\Bbb
F_p=(\Bbb Z/p)_{\bar O}(\bar O)$, 
$\theta _i$ is its characteristic function. Notice that  
$\sum\Sb i\endSb\theta _i=1$. In these terms
the corresponding counit $e_C$ and the comultiplication $\Delta _C$ are defined as
follows: $e_C(\theta _0)=1$, $e_C(\theta _i)=0$ if $i\ne 0$, and 
for all $i$, $\Delta _C(\theta _i)=\sum\Sb j_1+j_2=i\endSb \theta
_{j_1}\otimes\theta _{j_2}$. Fix a section $s:\Bbb F_p\longrightarrow
\bar G(\bar O)$ of the map induced by $\gamma $ on geometric points 
of $\bar G$ and $(\Bbb Z/p)_{\bar O}$. (Such section exists because
$\bar G$ is killed by $p$.) Prove that the section $\hat\gamma $ of $\gamma $ 
can be defined by the $\bar O$-linear morphism  
$\hat\gamma ^*:A(\bar G)\longrightarrow C$ such that 
for any $a\in A(\bar G)$, 
$\hat\gamma ^*(a)=\sum\Sb i\endSb a(s(i))\theta _i$. 

Clearly, $\hat\gamma ^*$ is a morphism of $\bar O$-algebras and 
$\hat\gamma ^*|_C=\id _C$.  It remains to verify that $\hat\gamma 
^*$ is compatible with the comultiplications $\Delta $ on $A(\bar H^D)$
and $\Delta _C$ on $C$. Let $a\in A(\bar H^D)$. Then 
$$\Delta _C(\hat\gamma ^*(a))=\Delta _C(\sum\Sb i\endSb a(s(i))\theta
_i)=\sum\Sb j_1,j_2\endSb a(s(j_1+j_2))\theta _{j_1}\otimes\theta
_{j_2}$$
$$=\sum\Sb j_1,j_1\endSb \Delta (a)(s(j_1),s(j_2))\theta
_{j_1}\otimes\theta _{j_2}=(\hat\gamma ^*\otimes\hat\gamma ^*)(\Delta (a)).$$
The lemma is proved. 
$\square $
\enddemo 

So, the elements of $\Ext _{\Gr _O}(H,\G_{\tilde\lambda })$ are
described via the classes $\cl (\varepsilon )\in H^2(H,\hat\Bbb
G_{m,\tilde\lambda })_p$ such that $\delta _{\tilde\lambda *}(\varepsilon
)\otimes _O\bar O\in B^2(\bar H,\hat\Bbb G_m)$. We can state this result
in the following form.

\proclaim{Proposition 5.5.3} Let $H\in\Gr _O$, $B=A(H)$, $\bar B=B\otimes
_O\bar O$. Let $\Cal H$ be the subgroup of 
$\hat\Bbb G_m(I_{\bar B})$ consisting of $\bar f\in 1+\m
(I_{\bar B})$ such that 
\medskip 

$\alpha )$\ \ $\bar f^p\in \hat\Bbb G_{m,\tilde\lambda ^p}(I_{B})$; 
\medskip 

$\beta )$\ \ $\delta ^{\times }\bar f=\Delta (\bar f)(\bar
f\otimes\bar f)^{-1}\in \hat\Bbb G_{m,\tilde\lambda }
(I_{B\otimes B})$. 
\medskip 

Then there is a group epimorphism 
$\Theta :\Cal H\longrightarrow\Ext _{\Gr _O}(H,\G_{\tilde\lambda })$  
such that for $\bar f\in\Cal H$, $\Theta (\bar f)=\Spec A\in\Gr _O$ with the
counit $e$ and the comultiplication $\Delta $ such that:
\medskip 

{\rm 1)}\ \ $A=B[X]$, where $\tilde\lambda ^{-p}\left
((1+\tilde\lambda X)^p-\bar f^p\right )=0$;
\medskip 

{\rm 2)}\ \ $e(X)=0$;
\medskip 

{\rm 3)}\ \ $1+\tilde\lambda\Delta (X)=[(1+\tilde\lambda X)\otimes
(1+\tilde\lambda X)]\cdot \delta ^{\times }\bar f$. $\square $
\endproclaim

Consider $\m (I_{\bar B})=I_{\bar B}\cap\m (\bar B)$ 
with the Lubin-Tate addition, cf. 4.5.  
For any $\bar f\in\m (I_{\bar B})$, set 
$$\delta _{\LT }(\bar f)=[\Delta (\bar f)]-[\bar f\otimes
1]-[1\otimes\bar f]\in\m (I_{\bar B\otimes \bar B}).$$ Then the above
proposition can be stated in the following equivalent form. 

\proclaim{Proposition 5.5.4} Let $\Cal H_{\LT }\subset\m (I_{\bar B})$ be
the subgroup (with respect to the Lubin-Tate addition) of $\bar f\in\m
(I_{\bar B})$ such that 
\medskip 
$\alpha _{\LT })$\ \ $[p](\bar f)\in\tilde\lambda ^pI_B$; 
\medskip 
$\beta _{\LT })$\ \ $\delta _{\LT }(\bar f)\in\tilde\lambda 
I_{B\otimes B}$. 
\medskip

Then there is a group epimorphism 
$\Theta _{\LT }:\Cal H_{\LT }\longrightarrow\Ext _{\Gr _O}(H,\G_{\tilde\lambda })$
such that $\forall\bar f\in\Cal H_{\LT }$, 
$\Theta _{\LT }(\bar f)=\Theta(E(\bar f))$, where 
$E$ is the Artin-Hasse exponential. $\square $
\endproclaim  

Notice that this proposition is obtained from proposition 5.5.3 just
by applying the Lubin-Tate logarithm. As a matter of fact this is a
first (though completely trivial) step towards relating the
multiplicative structures in the description of extensions 
from $\Ext _{\Gr _O}(H,G_{\tilde\lambda })$ in this section with  
additive constructions of the algebra and coalgebra structures of 
group schemes from $\operatorname{Im}\Cal G_O$ in Sections 2 and 3.

\medskip
\medskip 
\medskip

\subhead{\rm 6.} Calculations in the $O$-algebra of
$H\in\operatorname{Im}\Cal G_O$
\endsubhead 

In this section we use earlier notation and assumptions about 
$S,O,O_0,t,\pi,\pi _0$. 

\subsubhead{\rm 6.1.} Reminder and the statement of the Main Lemma 
\endsubsubhead 

Remind earlier notation and agreements:

$\bullet $\ \ $\Cal N\in\MF _S^e$, $H=\Cal G_O(\Cal N)$; $B=A(H)$ is the 
$O$-algebra of $H$ given in notation of n.3.2 as 
$O[X_1,\dots ,X_u]$ with the equations 
$X_i^p-\eta _i\sum\Sb j\endSb X_jc_{ji}=0,\ \ 1\leqslant i\leqslant
u$. 
Here $(c_{ij})\in\GL _u(O)$, all $\eta _i\in O_0$ and $\eta _i|p$. 
The comultiplication of $H$ appears as a unique 
$O$-algebra morphism $\Delta $ from $B$ to $B\otimes _{O}B$
such that 
$\Delta (X_i)=$ 
\linebreak 
$X_1\otimes 1+1\otimes X_i+j_i$, 
where all $j_i\in \Cal J_B\subset I^{DP}_{B\otimes B}$.  
Here the ideal $\Cal J_B$ is generated by
the elements $\tilde\eta _iX_i^r\otimes X_i^{p-r}$, where $1\leqslant
i\leqslant u$, $0<r<p$ and $\tilde\eta _i=-p/\eta _i$.

$\bullet $\ \ $i$ will be always an index such that $1\leqslant i\leqslant u$;
$\underline{i}$ will be always a multi-index $(i_1,\dots ,i_u)$, where 
$0\leqslant i_1,\dots ,i_u<p$; an index $i$ can be considered as a
special case of the multi-index $(\delta _{i1},\dots ,\delta _{iu})$,
where $\delta $ is the Kronecker symbol;  
$r(\underline{i})=i_1+\ldots +i_u$; 

$\bullet $\ \ $\bar O=O_{\bar K}$, $\bar{\m }$ is the maximal ideal in $\bar O$;  
$\bar B=B\otimes _{O}\bar O$, $I_{\bar B}=I_{B}\otimes
_{O}\bar O$, $\Cal J_{\bar B}=\Cal J_B\otimes _O\bar O$;

$\bullet $\ \  
we use the generators $X_1,\dots ,X_u$ for $\bar B$ and the generators
$X_1\otimes 1, \dots , $
\linebreak 
$X_u\otimes 1, 1\otimes X_1, \dots , 1\otimes X_u$ 
for $\bar B\otimes \bar B$ 
to define for any $\alpha\in\bar O$, the ideals $I_{\bar B}(\alpha )$
and $I_{\bar B\otimes\bar B}(\alpha )$ in the same way as we defined
in 3.3 for $\alpha\in O$, the ideals $I_{B}(\alpha )$ and
$I_{B\otimes B}(\alpha )$. We set for any $\alpha\in\bar O$, 
$I_{B}(\alpha )=I_{\bar B}(\alpha )\cap B$ and similarly, 
$I_{B\otimes B}(\alpha )=I_{\bar B\otimes\bar B}(\alpha )\cap
(B\otimes B)$;   

$\bullet $ \ \  
any element of $\bar B$ can be uniquely written 
as an $\bar O$-linear combination of all  
$X^{\underline i}=X_1^{i_1}\dots X_u^{i_u}$  
with all multi-indices $\underline {i}$. Similarly, 
any element of $\bar B\otimes\bar B$ can be uniquely written as an $\bar
O$-linear combination of all $X^{\underline {i}_1}\otimes
X^{\underline {i}_2}$ with multi-indices
$\underline{i}_1,\underline{i}_2$. We shall use the following obvious
property 
for $B$ (and its analogue for $B\otimes B$):
\medskip 

6.1.1)  
{\it if $\alpha\in\bar{\m }$, all $C_{\underline i}, D_{\underline i}\in\bar O$ and 
$\dsize\sum\Sb \underline i\endSb C_{\underline i}X^{\underline
i}\equiv \sum\Sb \underline i\endSb D_{\underline
i}X^{\underline i}\operatorname{mod}I_{\bar B}(\alpha )$ 
then for all $\underline i$,  
$C_{\underline i}X^{\underline i}\equiv D_{\underline i}X^{\underline
i}\operatorname{mod}I_{\bar B}(\alpha )$.} 
\medskip

$\bullet $\ \ in this section we calculate  
 in $\m (I_{\bar B})$ and $\m
(I_{\bar B\otimes\bar B})$, which are  provided with the Lubin-Tate addition
cf. 4.5; one must bear in mind the following agreement: if say,  
$a\in I_{\bar B}$ and it 
appears in the form $[a]$ then $a$ must be always considered as an element
of $\m (I_{\bar B})$;

$\bullet $ \ \ as earlier, we introduce $O'=O[\pi ']$, where $\pi
^{\prime p}=\pi $; all appropriate extensions of scalars from $O$ to
$O'$ will have the dashed notation, eg. $B'=B\otimes _OO'$, 
$I_{B'}$, $\Cal J_{B'}$ and so on.   
\medskip

\proclaim{Main Lemma} 
Suppose $\tilde\lambda\in O$, $\tilde\lambda ^{p-1}|p$. Suppose 
$f\in\m (I_{\bar B})$ and
$\delta _{\LT }(f)\in\tilde\lambda I_{B\otimes B}$. 
Then there are 

--- $f_0\in \tilde\lambda I_{B}$; 

--- for all $i$ and $0\leqslant l<p$, 
$o'_{il}\in \pi ^{\prime l}O$;

--- for all $\underline{i}$, $D_{\underline{i}}\in\bar O$ 
\medskip 
such that all $o_{il}^{\prime p}\in\tilde\lambda O$, $D_iX_i\in
I_{\bar B}(p)$
and 

$$f=[f_0]+\sum\Sb 0\leqslant l<p\\
1\leqslant i\leqslant s \endSb [o'_{il}X_i]+\sum\Sb 1\leqslant
i\leqslant s\endSb [D_iX_i]+\left [p\dsize\sum\Sb r(\underline i)\geqslant 2\endSb 
D_{\underline i}X^{\underline i}\right ]\operatorname{mod}p^2I_{\bar B}.\tag{6.1.2}$$
\endproclaim 

The Main Lemma will be proved in subsections 6.2-6.6 below. 

\remark {Remark} We need this lemma to study the extensions 
$G\in\Ext _{\Gr _O}(H,G_{\tilde\lambda })$. By Proposition 5.5.4 such
extensions appear from elements $f\in\Cal H_{\LT }$ satisfying the
conditions of our lemma. On the other hand, we expect that the
$O$-algebra $A(G)$ can be obtained (at least over $O'$) via the
special construction from Section 3. The Main Lemma shows that we can
replace $f$ by $[f]-[f_0]$, which gives rise to the same
extension $G$. Then in Section 7 we show that these special elements
from $\Cal H_{\LT }$ give rise to special extensions constructed in
Section 3. 
\endremark 
\medskip 

\subsubhead{\rm 6.2.} Auxiliary lemmas  
\endsubsubhead 

\proclaim{Lemma 6.2.1} Suppose 
$C_1,\dots ,C_u\in\bar O$, $g\in\m (I_{\bar B})$, 
 and there is  
$\beta_0\in \bar{\m }$ such that 
$g\equiv \sum\Sb i\endSb [C_iX_i]
\operatorname{mod}I_{\bar B}(\beta _0)$.  
Then there are $C'_{\underline i}\in\bar O$ such that: 
\medskip 
 
{\rm a)}\ $g=\sum\Sb i\endSb [C_i'X_i]+\left [\sum_{r(\underline i)\geqslant
2} C_{\underline i}'X^{\underline i}\right ]$;
\medskip 
  
{\rm b)}\ for all $1\leqslant i\leqslant u$, it holds 
$C_iX_i\equiv C_i'X_i\operatorname{mod}I_{\bar B}(\beta
_0)$;
\medskip 

{\rm c)}\   for all multi-indices $\underline{i}$ with $r(\underline{i})\geqslant
2$, it holds 
$C'_{\underline i}X^{\underline i}\in I_{\bar B}(\beta _0)$. 
\endproclaim 

\demo{Proof} Because $g\in\m (I_{\bar B})$ there is  $\alpha
_0\in\bar{\m }$ such that 
$g\in I_{\bar B}(\alpha _0)$ and we can assume that 
$I_{\bar B}(\alpha _0)\supset I_{\bar B}(\beta _0)$. Then   
all $C_iX_i\in I_{\bar B}(\alpha _0)$. Suppose 
$\beta\in\bar m$, $\beta _0|\beta $ and the statement of our
lemma is proved modulo $I_{\bar B}(\beta )$. (We can start with  
$\beta =\beta _0$.) Then 
$g=\sum\Sb \underline{i}\endSb  [C'_{\underline i}
X^{\underline i}]+[a]$, where $a\in I_{\bar B}(\beta )$. 
It is easy to show that such $a$ can be written as    
$a=\sum \Sb \underline {i} \endSb D_{\underline i}
X^{\underline i}$, where all $D_{\underline i}\in\bar O$ and
$D_{\underline i}X^{\underline i}\in I_{\bar B}(\beta )$. 
This implies that  
$a\equiv \sum\Sb \underline i\endSb [D_{\underline i}
X^{\underline i}]\operatorname{mod}I_{\bar B}(\beta ^p)$ 
and for all multi-indices 
$\underline i$, 
$$[C'_{\underline i}X^{\underline i}]+
[D_{\underline i}X^{\underline i}]\equiv 
[(C'_{\underline i}+D_{\underline i})X^{\underline i}]
\operatorname{mod}I_{\bar B}(\beta\alpha _0^{p-1}).$$ 
This proves our lemma modulo $I_{\bar B}(\beta\alpha _0^{p-1})$ 
and the proof can be finished by repeating this procedure.
$\square $
\enddemo

\proclaim{Lemma 6.2.2} Suppose $0\leqslant l<p$, $o_1',o_2'
\in\pi ^{\prime l}O$, $o_1^{\prime p}, o_2^{\prime p}\in\tilde\lambda O$. 
Then for any $a\in I_{B}$, there is $b\in I_B$ such that 
$$[o_1'a]+[o_2'a]=[(o_1'+o_2')a]+[\tilde\lambda b].$$ 
\endproclaim 

\demo{Proof} Just apply 4.5.1) and use that for $n\geqslant 1$, $\phi
_n(o_1',o_2')\in O\cap\tilde\lambda O'=\tilde\lambda O$. $\square $
\enddemo

\proclaim{Lemma 6.2.3} Suppose $1\leqslant i\leqslant u$, 
$o'\in O'$ and $o^{\prime
p}\in\tilde\lambda O$. Then there is 
\linebreak 
$a_i\in I_{B\otimes
B}$ such that 
$\delta _{\LT }[o'X_i]=[o'j_i]+[\tilde\lambda a_i]$. 
\endproclaim 

\demo{Proof} 
$\delta _{\LT }(o'X_i)=[o'(X_i\otimes 1+1\otimes X_i+j_i)]-
[o'(X_i\otimes 1)]-[o'(1\otimes X_i)]$ 
$$=[o'(X_i\otimes 1+1\otimes X_i)]+[o'j_i]-
\sum\Sb n\geqslant 1\endSb [o^{\prime p^n}\phi _n(X_i\otimes
1+1\otimes X_i,j_i)]$$
$$-[o'(X_i\otimes 1+1\otimes X_i)]-\sum\Sb n\geqslant 1\endSb
[o^{\prime p^n}\phi _n(X_i\otimes 1,1\otimes X_i)]=
[o'j_i]+[\tilde\lambda a_i],$$
with some $a_i\in I_{B\otimes B}$ because $o^{\prime p^n}\in
\tilde\lambda O$ if $n\geqslant 1$. $\square $
\enddemo

\definition{Definition} If $\alpha ,\beta\in\bar O$ then  
$I_{\bar B}(\alpha , \beta )= 
I_{\bar B}(\alpha )$ if
$\alpha |\beta $ and  $I_{\bar B}(\alpha ,\beta )=I_{\bar B}(\beta )$ 
if $\beta |\alpha $. (So, $I_{\bar B}(\alpha ,\beta )=I_{\bar B}(\alpha
)+I_{\bar B}(\beta )$.)  
Similarly, define the ideals $I_{\bar B\otimes\bar B}(\alpha ,\beta
)$. 
\enddefinition 

\proclaim{Lemma 6.2.4} Suppose $CX_i\in I_{\bar B}(\alpha )$, where
$C\in\bar O$ and $\alpha\in\bar{\m }$. Then 
\medskip 

{\rm a)} $\delta _{\LT }(CX_i)\equiv 
-[C^p\phi (X_i)]\operatorname{mod}
I_{\bar B\otimes\bar B}(\alpha ^{p^2},p^p)$; 
\medskip 

{\rm b)} if $\alpha =p$ then 
$\delta _{\LT }(CX_i)\equiv $
$$[Cj_i]-[C^p\phi (X_i)]-[C^p\phi
(X_i\otimes 1+1\otimes X_i, j_i)]
\operatorname{mod}(pI_{\bar B\otimes\bar B}\Cal J_{\bar B}+
p^2I_{\bar B\otimes \bar B}),$$
in particular, 
 $\delta _{\LT }(CX_i)\in \Cal J_{\bar B}+pI_{\bar B\otimes\bar B}$. 
\endproclaim 

\demo{Proof} First, notice that 
$$\delta _{\LT }(CX_i)=[C(X_i\otimes 1+1\otimes X_i+j_i)]-[CX_i\otimes
1]-[1\otimes CX_i]$$
$$=[C(X_i\otimes 1+1\otimes X_i)]+[Cj_i]-[C^p\phi (X_i\otimes
1+1\otimes X_i, j_i)]$$
$$-\sum\Sb n\geqslant 2\endSb [C^{p^n}\phi _n(X_i\otimes 1+1\otimes
X_i, j_i)]-[C(X_i\otimes 1+1\otimes X_i)]$$
$$-[C^p\phi (X_i)]-\sum\Sb n\geqslant 2\endSb [C^{p^n}\phi
_n(X_i\otimes 1, 1\otimes X_i)].$$

In the case a), the condition $CX_i\in I_{\bar B}(\alpha )$ implies
that  
all terms from the above both sums belong to $I_{\bar B\otimes\bar B}(\alpha ^{p^2})$. 
It remains to note that 
$C^p\phi (X_i\otimes 1+1\otimes X_i, j_i)$ and $Cj_i$ belong to
$I_{\bar B\otimes\bar B}(p^p)$, because $j_i\in I_{\bar B\otimes\bar B}(p^p)$.   

In the case b), $CX_i\in I_{\bar B}(p)$ implies that  
for all $n\geqslant 2$, 
$$C^{p^n}\phi _n(X_i\otimes 1, 1\otimes X_i)\in p^{p-1}I_{\bar
B\otimes\bar B}.$$ 
Indeed, 
$\phi _n$ is homogeneous of degree $p^n$ and is a linear combination
of terms 
\linebreak 
$X_i^{s_1}\otimes X_i^{s_2}$ with $s_1+s_2=p^n$. If $n\geqslant 2$
then  we can
apply to any such term at least $p-1$ times the relation $C^pX_i^p\in
pI_{\bar B}$. 
\medskip 

For $n\geqslant 2$, we have also that 
$$C^{p^n}\phi _n(X_i\otimes 1+1\otimes X_i, j_i)\in
pI_{\bar B\otimes\bar B}\Cal J_{\bar B}.$$
Indeed, we can use that $j_i^p\in p\Cal J_{\bar B}$ and that for $s>p^n-p\geqslant
2p$, the elements $C^s(X_i\otimes 1+1\otimes X_i)^s$ belong to 
$pI_{\bar B\otimes\bar B}$.  
The lemma is proved. 
$\square $
\enddemo 

\proclaim{Lemma 6.2.5} If $\alpha\in\bar{\m }$ and  for all 
$\underline i$ with $r(\underline i)\geqslant 2$, it holds 
$C_{\underline i}X^{\underline i}\in I_{\bar B}(\alpha ^p)$ then 
$$\delta _{\LT }\left (\sum_{\underline i }
C_{\underline i}X^{\underline i}\right )\equiv \sum\Sb \underline i \endSb 
C_{\underline i}\sum\Sb \underline{\jmath}'+\underline{\jmath
}''=\underline i\endSb
A_{\underline{\jmath}'\underline{\jmath}''}X^{\underline{\jmath}'}\otimes
X^{\underline{\jmath}''}\operatorname{mod}\left (I_{\bar B\otimes\bar B}(\alpha
^{p^2})+I_{\bar B\otimes\bar B}\Cal J_{\bar B}\right ),$$
where $r(\underline\jmath '), r(\underline\jmath ^{\prime\prime })>0$
and all coefficients
$A_{\underline{\jmath}'\underline{\jmath}''}\in\Bbb Z_p^*$. 
\endproclaim 

\demo{Proof} Notice that  
 the
Lubin-Tate group law  on $I_{\bar B\otimes\bar B}(\alpha ^p)$ 
can be replaced 
modulo $I_{\bar B\otimes\bar B}(\alpha ^{p^2})$ 
by the usual addition. Then use
that  
for any multi-index $\underline i$ with $r(\underline i)\geqslant 2$,
it holds   
$\Delta (X^{\underline i})\equiv 
(X_1\otimes 1+1\otimes X_1)^{i_1}\dots (X_s\otimes 1+1\otimes
X_s)^{i_s}
\operatorname{mod}I_{\bar B\otimes\bar B}\Cal J_{\bar B}$, 
where all appropriate binomial coefficients are prime to $p$. 
This implies the statement of our lemma. 
$\square $
\enddemo 

\remark{Remark 6.2.6} Notice that we've just  proved that 
$\delta ^+\left (\sum\Sb \underline{i}\endSb 
C_{\underline i}X^{\underline i}\right )$ (where $r(\underline
i)\geqslant 2$) is congruent 
modulo the ideal 
$I_{\bar B\otimes\bar B}J_{\bar B}$ 
to the
right-hand side of the formula from above Lemma 6.2.5 . 
\endremark 
\medskip

\subsubhead{\rm 6.3.} Step 1 
\endsubsubhead 

Suppose $\alpha\in\bar{\m }$, $\alpha |p$ is such that the statement
of the Main Lemma  
holds modulo $I_{\bar B}(\alpha ^p)$ with all $D_iX_i\in I_{\bar
B}(\alpha )$, i.e.  
$$f\equiv [f_{\alpha }]+\sum\Sb i,l\endSb [o'_{il}X_i]+\sum\Sb i\endSb
[D_iX_i]\operatorname{mod}I_{\bar B}(\alpha ^p),\tag{6.3.1}$$
where all $D_iX_i\in I_{\bar B}(\alpha )$, $f_{\alpha }\in
\tilde\lambda I_{B}$ and  
all $o'_{il}\in \pi ^{\prime l}O$ are such that 
$o^{\prime p}_{il}\in \tilde\lambda O$. Such $\alpha $ always
exists, e.g. we can take $\alpha =\alpha _0^{1/p}$, where $\alpha
_0\in\bar{\m }$ is such that $f\in I_{\bar B}(\alpha _0)$. 
We are going to prove a similar congruence 
modulo the smaller ideal $I_{\bar B}(\alpha ^{p^2}, p^p)$. 

Apply Lemma 6.2.1 to 
$g=[f]-[f_{\alpha }]-\sum\Sb i,l \endSb [o'_{il}X_i]$
and 
$\beta _0=\alpha ^p$. Then 
\linebreak 
$f=[f_{\alpha }]+\sum\Sb i,l\endSb 
[o'_{il}X_i]+\sum\Sb i\endSb [C_iX_i]+
\left [\sum\Sb r(\underline i)\geqslant 2\endSb 
C_{\underline i}X^{\underline i}\right ]$, 
where all $C_{\underline i}\in\bar O$, 
for all $1\leqslant i\leqslant u$, it holds $C_iX_i\equiv
D_iX_i\operatorname{mod}I_{\bar B}(\alpha ^p)$ (and therefore 
all $C_iX_i\in I_{\bar B}(\alpha )$) and for 
all $\underline i$ with 
$r(\underline i)\geqslant 2$, it holds $C_{\underline
i}X^{\underline i}\in I_{\bar B}(\alpha ^p)$. 

Now apply Lemmas 6.2.3-6.2.5 and 
notice that $\Cal J_{\bar B}\subset I_{\bar B\otimes\bar B}
(\alpha ^{p^2}, p^p)$. Then 
the condition $\delta _{\LT }(f)\in \tilde\lambda I_{B\otimes B}$
implies that  

$$-\sum\Sb i \endSb C_i^p\phi (X_i)+\sum\Sb 
r(\underline i)\geqslant 2\endSb C_{\underline i}
\sum\Sb \underline{\jmath}'+\underline{\jmath}''=\underline
i\endSb A_{\underline{\jmath}'\underline{\jmath}''}
X^{\underline{\jmath}'}\otimes X^{\underline{\jmath}''}
\in \tilde\lambda I_{B\otimes B}\operatorname{mod}I_{\bar
B\otimes\bar B}(\alpha ^{p^2},
p^p).$$

Notice that all monomials in the both above sums are different 
and belong to $\bar O^{<p}[X_1\otimes 1,\dots ,X_u\otimes 1, 1\otimes
X_1,\dots ,1\otimes X_u]$.  
Due to Remark 6.1.1 this implies the following two facts:
\medskip 

1) for all $1\leqslant i\leqslant u$, $C_i^p\phi (X_i)\equiv o_i\phi
(X_i)\operatorname{mod}I_{\bar B\otimes\bar B}(\alpha ^{p^2}, p^p)$, where
$o_i\in\tilde\lambda O$; 
\medskip 

2) if $r(\underline i)\geqslant 2$ then 
$C_{\underline i}X^{\underline i}\equiv o_{\underline i}X^{\underline
i}
\operatorname{mod}I_{\bar B}(\alpha ^{p^2}, p^p)$, where 
$o_{\underline i}\in\tilde\lambda O$. 
\medskip 

The first fact implies that $(C_i^p-o_i)^p\eta _i^p\equiv
0\operatorname{mod}(\alpha ^{p^2},p^p)$ and, therefore, 
$$(C_i^p-o_i)\eta _i\equiv 0\operatorname{mod}(\alpha ^p,p).$$
  Decompose each $o_i$ in the
form 
$o_i\equiv\dsize\sum\Sb 0\leqslant l<p\endSb o^{\prime\prime p}_{il}
\operatorname{mod}p\tilde\lambda $, 
where for all $0\leqslant l<p$, 
$o^{\prime\prime }_{il}\in \pi ^{\prime l}O$. Notice that all  
$o^{\prime\prime p}_{il}\in\tilde\lambda
O$. Then 
$(C^p_i-\dsize\sum\Sb l \endSb o^{\prime\prime p}_{il})\eta _i\equiv
0\operatorname{mod}(\alpha ^p,p)$ 
or, equivalently, 
$$C_iX_i\equiv \dsize\sum\Sb l\endSb o^{\prime\prime }_{il}
X_i\operatorname{mod}I_{\bar B}(\alpha ^p,p).$$ 
Notice that also all $C_iX_i, o^{\prime\prime }_{il}X_i\in I_{\bar B}(\alpha )$. 
Let $\dsize\sum\Sb l\endSb o^{\prime\prime }_{il}=o_i'$. Then 
there are $C'_{ij}, C^{\prime\prime }_{ij}\in\bar O$ such that by 4.5
a), 

$$[C_iX_i]-[o'_{i}X_i]
\equiv [(C_i-o'_i)X_i]+[\phi (C_i, -o_i')X_i^p]
\equiv\sum \Sb j\endSb [C'_{ij}X_j]\operatorname{mod}I_{\bar B}(\alpha
^{p^2},p^p)$$  
(use that $X_i^p$ is a linear combination of $X_j$, $1\leqslant
j\leqslant u$, and that all terms in the middle and in the 
right hand side belong to $I_{\bar B}(\alpha ^p, p)$)
and for similar reasons 
$$[o_i'X_i]=\sum\Sb l \endSb [o^{\prime\prime }_{il}X_i]+
\sum\Sb j\endSb [C_{ij}^{\prime\prime }X_j]
\operatorname{mod}I_{\bar B}(\alpha
^{p^2},p^p),$$  
where all $C'_{ij}X_i, C^{\prime\prime }_{ij}X_j\in I_{\bar B}(\alpha ^{p},p)$. 
This gives, finally, that 
$$\sum\Sb i\endSb [C_iX_i]=\sum\Sb i,l\endSb
[o^{\prime\prime }_{il}X_i]+\sum\Sb i\endSb 
[D_i'X_i]\operatorname{mod}I_{\bar B}(\alpha ^{p^2},p^p),$$
where all $D_i'\in\bar O$ and $D_i'X_i\in I_{\bar B}(\alpha ^{p},p)$. 

The fact 2) means that with $\tilde\lambda f'=\sum\Sb r(\underline
i)\geqslant 2\endSb o_{\underline i}X^{\underline i}\in\tilde\lambda
I_B$, it holds 
$$\sum\Sb r(\underline i)\geqslant 2\endSb C_{\underline
i}X^{\underline i}\equiv \tilde\lambda f'\operatorname{mod}I_{\bar
B}(\alpha ^{p^2},p^p).$$
By Lemma 6.2.4 there is $f^{\prime\prime }\in \tilde\lambda I_B$ such that 
$$f\equiv [f^{\prime\prime }]+\sum\Sb i,l\endSb 
[(o'_{il}+o^{\prime\prime }_{il})X_i]+\sum\Sb i\endSb 
[D_i'X_i]\operatorname{mod}I_{\bar B}(\alpha ^{p^2},p^p)$$
and we obtained an analogue of (6.3.1) modulo $I_{\bar B}(\alpha
^{p^2},p^p)$ with all $D_iX_i\in I_{\bar B}(\alpha ^p,p)$.  

If $I_{\bar B}(\alpha ^{p^2},p^p)=I_{\bar B}(\alpha ^{p^2})$
repeat this step with $\alpha $ replaced by $\alpha ^p$. Then in
finitely many steps we obtain that 
$I_{\bar B}(\alpha ^{p^2},p^p)=I_{\bar B}(p^p)$. This means that we 
proved an 
analogue of formula (6.3.1) modulo $I_{\bar B}(p ^p)$ (with all
$D_iX_i\in I_{\bar B}(p)$). 
\medskip 

\subsubhead{\rm 6.4.} Step 2 
\endsubsubhead 

Suppose $\alpha\in\bar O$ is such that $p|\alpha $ and suppose  

$$f\equiv [f_{\alpha }]+\sum\Sb i,l\endSb [o'_{il }X_i]+\sum\Sb
i\endSb [D_iX_i]
\operatorname{mod}\left (I_{\bar B}(\alpha ^p)+pI_{\bar B}\right ),\tag{6.4.1}$$
where as earlier, $f_{\alpha }\in \tilde\lambda I_{B}$, 
for all $i$ and $0\leqslant l<p$, $o'_{il}\in \pi ^{\prime l}O$,  
$o^{\prime p}_{il}\in\tilde\lambda O$,   
$D_i\in\bar O$ and $D_iX_i\in I_{\bar B}(p)$. This congruence holds for $\alpha
=p$ by results of n.6.3.

Prove that there is a similar 
congruence modulo $I_{\bar B}(\alpha ^{p^2})+pI_{\bar B}$. 

Apply Lemma 6.2.1. Then 

$$f\equiv [f_{\alpha }]+
\sum\Sb i,l \endSb [o'_{il}X_i] 
+\sum\Sb i\endSb [C_iX_i]+\left [\sum\Sb r(\underline i)\geqslant 2
\endSb C_{\underline i}X^{\underline i}\right
]\operatorname{mod}pI_{\bar B},$$
where all $C_{\underline i}\in\bar O$, 
$C_iX_i\equiv D_iX_i\operatorname{mod}I_{\bar B}(p)$ and all terms from the
last sum belong to $I_{\bar B}(\alpha ^p)$. 
Apply Lemmas 6.2.3-6.2.5. Then  

$$
\sum\Sb r(\underline i)\geqslant 2\endSb C_{\underline i}
\sum\Sb \underline{\jmath}'+\underline{\jmath}''=
\underline i\\
r(\underline{\jmath }'),r(\underline{\jmath}^{\prime\prime })>0\endSb 
A_{\underline{\jmath}'\underline{\jmath}''}
X^{\underline{\jmath}'}\otimes X^{\underline{\jmath}''}
\in\tilde\lambda I_{B\otimes B}\operatorname{mod}\left (\Cal J_{\bar B}+pI_{\bar
B\otimes\bar B}+I_{\bar B\otimes\bar B}(\alpha ^{p^2})\right ).\tag{6.4.2}$$

\proclaim{Lemma 6.4.3} Suppose 
$\dsize\sum\Sb \underline i_1,\underline i_2\endSb 
a_{\underline i_1\underline i_2 }X^{\underline i_1}\otimes
X^{\underline i_2}\in \bar O^{<p}[X_1\otimes 1, \dots ,1\otimes X_u]$ 
belongs to ${\Cal J_{\bar B}}$. 
If multi-indices 
$\underline i^0_1,\underline i_2^0$ are such that 
the total degree of the monomial 
$X^{\underline i_1^0}\otimes X^{\underline i_2^0}$ is less than $p$
in each variable $X_1,\dots ,X_u$ then $a_{\underline i_1^0\underline
i_2^0}\in p\bar O$. 
\endproclaim 

\demo{Proof} The elements of $\Cal J_{\bar B}$ are $\bar O$-linear combinations
of the terms 
\linebreak 
$\tilde\eta _i(X_i^{r}\otimes X_i^{p-r})b$, where 
$1\leqslant i\leqslant u$, $1\leqslant r<p$ and $b$ is a monomial of the form 
$X^{\underline{\jmath }_1}\otimes X^{\underline {\jmath }_2}$. 
Such product makes a non-zero contribution
to $a_{\underline i_1^0\underline i_2^0}$ only if its total degree in
$X_i$ 
will become $<p$. This will be a chance only if $X_i^p$ appears in 
a left or right side 
of this tensor product. It remains to notice that
$\tilde\eta _iX_i^p\in pI_{B}$. 
$\square $
\enddemo 

The above lemma together with relation (6.4.2) implies that for any
$\underline i$ with 
\linebreak 
$r(\underline i)\geqslant 2$, there is an
$o_{\underline i}\in\tilde\lambda O_1$ such that 
$C_{\underline i}X^{\underline i}\equiv o_{\underline i}X^{\underline
i}
\operatorname{mod}\left (I_{\bar B}(\alpha ^{p^2})+pI_{\bar B}\right )$. 
Therefore, 

$$f\equiv [f_{\alpha ^p}]+
\sum\Sb i,l\endSb [o'_{il}X_i]+\sum\Sb i\endSb [C_iX_i]
\operatorname{mod}\left (I_{\bar B}(\alpha ^{p^2})+pI_{\bar B}\right ),$$
where $f_{\alpha ^p}\in \tilde\lambda I_{B}$ is such that 
$[f_{\alpha }]-\left [\sum\Sb \underline i\endSb 
o_{\underline i}X^{\underline i}\right ]=[f_{\alpha
^p}]$. 
\medskip 

\subsubhead{\rm 6.5.} Step 3
\endsubsubhead 

By repeating the procedure from n.6.4 sufficiently many times 
we shall obtain  

$$f=[f_0]+\sum\Sb i,l \endSb [o'_{il }X_i]+\sum\Sb
i \endSb [D_iX_i]+[pg],\tag{6.5.1}$$
where its ingredients $f_0$, $o'_{il}$ and $D_i$ satify 
the requirements of the Main Lemma.
 
Let $g=\sum\Sb \underline i\endSb A_{\underline i}X^{\underline i}$
with all $A_{\underline i}\in\bar O$. It remains to prove that we can
get rid of all linear terms $A_iX_i$ in $g\operatorname{mod}pI_{\bar
B}$. 

Notice that for all indices
$1\leqslant i\leqslant u$, 
$$[C_iX_i]+[pA_iX_i]\equiv 
[(C_i+pA_i)X_i]+[\phi (C_i,pA_i)X_i^p]\operatorname{mod}p^2I_{\bar B},$$
because for $n\geqslant 2$, $\phi _n(C_i,pA_i)X_i^{p^n}\equiv
0\operatorname{mod}p^2$ (use that $C_i^{p^n-1}X_i^{p^n}$ is divisible
by $C_i^pX_i^p\in pI_{\bar B}$). Notice also that 
$$\phi (C_i,pA_i)X_i^p\equiv -pC_i^{p-1}A_iX_i^p\equiv p\sum\Sb j\endSb
C_{ij}X_j\operatorname{mod}p^2I_{\bar B},$$
where all $C_{ij}\in\bar O$ are divisible by 
$C_i^{p-1}\eta _i\equiv 0\dsize\operatorname{mod}p^{1-1/p}$, 
because $C_i^p\eta _i\equiv 0\operatorname{mod}p$. This
implies that 

$$\sum\Sb i\endSb [C_iX_i]+\sum\Sb i\endSb [pA_iX_i]\equiv $$

$\sum\Sb i\endSb [(C_i+pA_i)X_i]+\sum\Sb i,j\endSb 
[pC_{ij}X_j]\equiv 
\sum\Sb i\endSb \left [(C_i+pA_i+p\sum \Sb j\endSb C_{ji})X_i\right ]
\operatorname{mod}p^{2-1/p}I_{\bar B}$. 

\ \ 
\newline 
This relation implies an analogue of formula (6.5.1) with  $D_i$
replaced by $C_i+pA_i+p\sum\Sb j\endSb C_{ji}$ and where  
$g\equiv \sum\Sb r(\underline i)\geqslant 2\endSb A_{\underline
i}X^{\underline i}\operatorname{mod}p^{1-1/p}I_{\bar B}$. 
Repeating this step one time more we shall get a similar congruence
for $g$ modulo $pI_{\bar B}$, i.e. that a new $g$ will not 
contain linear terms $A_iX_i$ modulo $pI_{\bar B}$. 

The Main Lemma is proved. 
$\square $
\medskip 
\medskip

\subsubhead{\rm 6.6.} Application of the Main Lemma  
\endsubsubhead 

\proclaim{Proposition 6.6.1} Suppose $f\in \m(I_{\bar B})$ 
satisfies the assumptions of
the Main Lemma and is given by the corresponding formula  (6.1.2). 
Then for $1\leqslant i\leqslant u$, there are $o_i\in\tilde\lambda O$
such that 

$$\sum\Sb r \endSb (D_r+\sum\Sb 0\leqslant l<p\endSb o'_{rl})
\tilde\eta _id_{ir}-D_i^p\equiv o_i\operatorname{mod}p\tilde\eta
_i.$$ 
\endproclaim 

\demo{Proof} Via Lemmas 6.2.3-6.2.5 the condition $\delta _{\LT
}(f)\in\tilde\lambda I_{B\otimes B}$ implies that 

$$\sum\Sb r\endSb \left (\sum\Sb l \endSb [o'_{rl}j_r]+[D_rj_r] 
-[D_r^p\phi (X_r)]-[D_r^p\phi
(X_r\otimes 1+1\otimes X_r, j_r)]\right )+
[p\delta ^+g_0]\tag{6.6.2}$$ 
belongs to $\tilde\lambda I_{B\otimes
B}$ modulo $pI_{\bar B\otimes\bar B}{\Cal J_{\bar B}}+p^2I_{\bar B\otimes\bar B}$. 
We are going to follow the coefficient for $\phi (X_i)$ in this
formula written as an element of $\bar O^{<p}[X_1\otimes 1,\dots
,1\otimes X_u]$.  

\proclaim{Lemma 6.6.3} Elements from the ideal $I_{\bar B\otimes\bar B}{\Cal
J_{\bar B}}$ contain  the  
monomials 
\linebreak 
$X_i^j\otimes X_i^{p-j}$, where $1\leqslant i\leqslant u$
and $1\leqslant j<p$, with coefficients divisible by $p$. 
\endproclaim 

\demo{Proof of lemma} It is quite similar to the proof of Lemma
6.4.2. 
$\square $
\enddemo 

Now the proof of our Proposition 6.6.1 can be finished as follows:  
\medskip  

---  working with 
formula (6.6.2) modulo the ideal $I_{\bar B\otimes\bar B}(p^{2p})$ 
we can find the coefficient
for $\phi (X_i)$ modulo $p\tilde\eta _i$. Indeed, if 
$C\phi (X_i)\in I_{\bar B\otimes\bar B}(p^{2p})$ then $C^p\eta _i^p\equiv
0\operatorname{mod}p^{2p}$ and $C\equiv 0\operatorname{mod}p\tilde\eta
_i$; 
\medskip 

---   
if we take relation (6.6.2) modulo $I_{\bar B\otimes\bar B}(p^{2p})$ we can
replace the Lubin-Tate group law by the usual addition, because all
terms belong to $I_{\bar B\otimes\bar B}(p^p)$; 
\medskip 

--- the term $p\delta ^+g_0$ gives the zero contribution
to the coefficient for $\phi (X_i)$ modulo $p^2$, because of 
Remark 6.2.6 and above Lemma 6.6.3;
\medskip 

--- the term $D_r^p\phi (X_r\otimes 1+1\otimes X_r,j_r)$
gives the zero contribution to the coefficient for $\phi (X_i)$
modulo $p\tilde\eta _i$ (apply similar arguments as in the proof of
Proposition 3.5.1); 
\medskip 

--- finally, use Proposition 3.5.1 to find out that the 
coefficient for $\phi (X_i)$ in (6.6.2) coincides 
modulo  $p\tilde\eta _i$
with the left-hand side of the formula from our proposition. 
\medskip 
The proposition is proved. 
$\square $
\enddemo

\subhead 7. Epimorphic property of $\Cal G^O_{O_0}$  
\endsubhead 
\medskip 

In this section we use the notation and assumptions 
about $O_0$ and $O$ from n.4. As in Section 6 we use $O'=O[\pi ']$
where $\pi ^{\prime p}=\pi $. 
We are going to prove that for any $G_0\in\Gr _{O_0}$ there is an $\Cal
M\in\MF ^e_{S}$ such that $G_0=\Cal G_{O_0}^O(\Cal M)$ or,
equivalently, 
$G:=G_0\otimes _{O_0}O=\Cal G_O(\Cal M)$. Notice
that by the Tate-Oort classification of group schemes of order $p$, 
 [TO],  
this is true for group schemes of order $p$. Therefore, we can
assume that the order $|G_0|$ of $G_0$ is bigger than $p$ and 
the above property holds  
for all $H_0\in\Gr _{O_0}$ such that $|H_0|<|G_0|$. 

7.1. By results of Section 4 we can replace $O_0$ by the valuation ring of
sufficiently large tamely ramified extension of $K_0$. Therefore, we can
assume the existence of $\tilde\lambda\in O_0$ such that
$\tilde\eta:=\tilde\lambda ^{p-1}$ divides $p$ and if $\eta
=-p/{\tilde\eta }$ then $G_0\in\Ext _{\Gr _{O_0}}(H_0, \G_{0\tilde\eta })$.  Here 
$H_0=\Spec B_0\in\Gr _{O_0}$ and $\G_{0\tilde\eta }=\Cal G^O_{O_0}(\Cal M_{\tilde s})$,
where $\Cal M_{\tilde s}\in\MF ^e_S$ and $\kappa _{SO}(\tilde s^p)=
\tilde\eta\operatorname{mod}p$. By inductive assumption there is $\Cal
N\in\MF ^e_S$ such that 
$H_0=\Cal G^O_{O^0}(\Cal N)$, where $\Cal N\in\MF _S^e$.

We assume that the structure of $\Cal N$ as an object of the category
$\MF _S^e$ as well as the structure of the coalgebra 
$B=B_0\otimes _{O_0}O$ are given in notation from n.3.2. By enlarging
(if necessary) the residue field $k$ we can assume that 
$\tilde\lambda $ and all $\tilde\eta _i$, $1\leqslant i\leqslant u$, 
are just powers of the uniformising element $\pi _0$ of $O_0$. 

Let $\bar f\in\m (I_{\bar B})$ be such that $G_0=\Theta _{\LT }(\bar
f)$, cf. 5.5.4.  
Notice that 
$\bar f$ is 
defined modulo $\tilde\lambda I_{B_0}\cap\m (I_{B_0})$
(with respect to the Lubin-Tate addition) and satisfies the
requirements 
$$\delta _{\LT }(\bar f)\in\tilde\lambda I_{B_0\otimes B_0},\ \ 
[p](\bar f)\in\tilde\lambda ^pI_{B_0}. \tag {7.1.1} $$

Apply the Main Lemma from Section 6 
to $\bar f$. Then there is $f_0\in\m (I_{B})\cap\tilde\lambda I_B$ such that 

$$[\bar f]\equiv [f_0]+\sum\Sb 1\leqslant i\leqslant u\\ 0\leqslant l<p\endSb 
[o'_{il}X_i] +\sum\Sb 1\leqslant i\leqslant u\endSb
[D_iX_i]\operatorname{mod}pI_{\bar B},$$
where all $o'_{il}\in\pi ^{\prime l}O\subset O'$, 
$o_{il}^{\prime p}\in\tilde\lambda
O$, $D_i\in\bar O$ and $D_iX_i\in I_{\bar B}(p)$. By 
Proposition 6.6.1 
for all $i$, 
$$\sum\Sb 1\leqslant r\leqslant u\endSb (D_r+\sum\Sb 0\leqslant
l<p\endSb o'_{rl})\tilde\eta _id_{ir}-D_i^p\equiv o_i
\operatorname{mod}p\tilde\eta _i,\tag{7.1.2}$$
where all $o_i\in\tilde\lambda O$. Notice that 
$[\bar f_1]:=[\bar f]-[f_0]$ corresponds to $G\in\Ext
_{\Gr _{O}}(H,G_{\eta })$ under the map $\Theta _{\LT }$ from 5.5.4. 
Also notice that congruences 
(7.1.2) imply that all $o_i\in\tilde\eta _iO$ (use that 
$D_i^p\equiv 0\operatorname{mod}\tilde\eta _i$ because $D_iX_i\in
I_{\bar B}(p)$). 

For $1\leqslant i\leqslant u$ and $0\leqslant l<p$, introduce
$o_{il}^{\prime\prime }\in\pi ^{\prime l}O\subset O'$ such that 
$$o_i\equiv \sum\Sb l\endSb o_{il}^{\prime\prime
p}\operatorname{mod}p\tilde\eta _i.$$
 Clearly, all $o_{il}^{\prime\prime p}\in\tilde\lambda O\cap
\tilde\eta _iO$. 
Set $o_{il}:=o_{il}'-o_{il}^{\prime\prime }$. 

\proclaim{Proposition 7.1.3} There is $h\in \tilde\lambda I_{B}\cap
I_{B}(p)$ such that for $[\bar f']=[\bar f_1]-[h]$, it holds 

$$l_{\LT }(\bar f')\equiv -\sum\Sb i,l\endSb o_{il}X_i+\sum\Sb i,l \endSb
l_{\LT }(o_{il}X_i)\operatorname{mod}pI_{\bar B}.$$
\endproclaim 

\demo{Proof} 
Proceed with the following computation modulo $pI_{\bar B}$.

$$l_{\LT }(\sum\Sb i\endSb [D_iX_i])=\sum\Sb i\endSb l_{\LT }(D_iX_i)
\equiv \sum\Sb i\endSb \left (D_iX_i+D_i^pX_i^p/p\right )$$

$$\equiv \sum\Sb i \endSb D_iX_i -\sum\Sb i\endSb
\frac{o_iX_i^p}{p}+
\frac{1}{p}\sum\Sb i,r \endSb (D_r+\sum\Sb l \endSb 
o'_{rl })\tilde\eta _id_{ir}X_i^p\equiv $$

$$-\sum\Sb i,l \endSb o'_{il}X_i-\frac{1}{p}\sum\Sb i\endSb o_i
X_i^p\equiv-\sum\Sb i,l\endSb o_{il}X_i-\sum\Sb i,l\endSb 
l_{\LT }(o^{\prime\prime }_{il}X_i).$$
(Use that $D_iX_i\in I_{\bar B}(p)$ and $\sum\Sb i\endSb\tilde\eta
_id_{ir}X_i^p=
\sum\Sb i\endSb \tilde\eta _id_{ir}\eta _i
\sum\Sb j\endSb X_jc_{ji}=-pX_r$.)

Therefore, 
$$l_{\LT }(\bar f_1)\equiv-\sum\Sb i,l\endSb o_{il}X_i+\sum\Sb i,l\endSb
l_{\LT }([o_{il}'X_i]-[o_{il}^{\prime\prime
}X_i])\operatorname{mod}pI_{\bar B}.$$
Now note that for all $i$ and $l$, 
$[o'_{il}X_i]-[o_{il}^{\prime\prime }X_i]-[o_{il}X_i]=[h_{il}]$, 
where 
$$h_{il}=\sum\Sb n\geqslant 1\endSb \left [X_i^{p^n}\phi
_n(o'_{il},o^{\prime\prime }_{il})\right ]\in\tilde\lambda I_{B}\cap
I_{B}(p).$$

Indeed by Lemma 6.2.2, all $h_{il}\in\tilde\lambda I_B$ and notice
that 
for all $n\geqslant 1$, 
 $\phi _n(o'_{il},o^{\prime\prime }_{il})\equiv
0\operatorname{mod}o^{\prime\prime }_{il}$ (because $\phi _n(X,0)=0$)
and $o^{\prime\prime }_{il}X_i\in I_{B}(p)$ (because 
$o^{\prime\prime p}_{il}\equiv 0\operatorname{mod}\tilde\eta _i$). 

So, if $[h]=\sum\Sb i,l\endSb [h_{il}]$ then 
$h\in \tilde\lambda I_{B}\cap I_{B}(p)$ and 
$$l_{\LT }(\bar f')=l_{\LT }([\bar f_1]-\sum\Sb i,l \endSb [h_{il}])
\equiv -\sum\Sb i,l\endSb o_{il}X_i+\sum\Sb
i,l\endSb l_{\LT }(o_{il}X_i)\operatorname{mod}pI_{\bar B}.$$

The proposition is proved. 
$\square $
\enddemo 

7.2. Let $[g]=[\bar f']-\sum\Sb i,l \endSb
[o_{il}X_i]$. Remind that $B'=B\otimes _OO'$ with the augmentation
ideal $I_{B'}=I_B\otimes _OO'$ and $\Cal J_{B'}=\Cal J_B\otimes _OO'$. 

\proclaim{Proposition 7.2.1} 
\newline 
{\rm a)} $[p](g)\equiv  -p\sum\Sb i,l\endSb 
o_{il}X_i\operatorname{mod}p^2I_{B'}$; 
\newline 
{\rm b)} $\delta _{\LT }(g)\in\tilde\lambda ^{1/p}\Cal
J_{B'}+pI_{B'\otimes B'}$. 
\endproclaim 

\demo{Proof} a) In the notation from the proof of proposition 7.1.3 it
holds   

$$[g]=[\bar f']-\sum\Sb i,l \endSb [o_{il}X_i]=[\bar f_1]-\sum\Sb i,l \endSb
[o'_{il}X_i]+\sum\Sb i,l\endSb [o^{\prime\prime }_{il}X_i]
=\sum\Sb i\endSb [D_iX_i]+\sum\Sb i,l\endSb [o^{\prime\prime
}_{il}X_i]\in I_{B'}(p).$$
(Use that $[p](\bar f_1)\in\tilde\lambda ^pI_B$.) 
Therefore, 
$[p](g)\equiv [g^p]+[pg]\equiv 0
\operatorname{mod}pI_{B'}$ and  
$l_{\LT }([p]g)\equiv [p]g\operatorname{mod}p^2I_{B'}$. So, by
Proposition 7.1.3 it holds 
$$l_{\LT }([p]g)=pl_{\LT }(g)\equiv -p\sum\Sb i,l\endSb
o_{il}X_i\operatorname{mod}p^2I_{B'}.$$

b) As earlier, let $\delta ^+=\Delta -\id\otimes 1-1\otimes\id $. Then 
$$l_{\LT }(\delta _{\LT }(g))=\delta ^+l_{\LT }(g)\equiv -\sum\Sb
i,l\endSb o_{il}\delta ^+X_i\equiv 0\operatorname{mod}\left (\tilde\lambda
^{1/p}\Cal J_{B'}+pI_{B'\otimes B'}\right )\tag{7.2.2}$$
because all $o_{il}\equiv 0\operatorname{mod}\tilde\lambda ^{1/p}$ and 
$\delta ^+X_i\in\Cal J_B$, cf. Proposition 3.2.3. 
On the other hand by Lemma 6.2.4 
$$\delta _{\LT }(g)\equiv \sum\Sb i\endSb [\delta _{\LT }(D_iX_i)]+\sum\Sb
i,l\endSb[\delta _{\LT }(o^{\prime\prime }_{il}X_i)]\equiv 0
\operatorname{mod}(\Cal
J_{\bar B}+pI_{\bar B\otimes\bar B}).$$
Now notice that $\Cal J_{\bar B}+pI_{\bar B\otimes\bar B}\subset I^{DP}_{\bar
B\otimes\bar B}$,  
$l_{\LT }$ induces a one-to-one 
transformtion of  
$I^{DP}_{\bar B\otimes\bar B}$, 
$\tilde\lambda ^{1/p}\Cal J_{B'}+pI_{B'\otimes B'}\subset
I^{DP}_{\bar B\otimes\bar B}$, 
$l_{\LT }(\tilde\lambda ^{1/p}\Cal J_{B'})=\tilde\lambda
^{1/p}\Cal J_{B'}$ 
and $l_{\LT }(pI_{B'\otimes B'})=pI_{B'\otimes B'}$. 
Therefore, (7.2.2) implies that 
$\delta _{\LT }(g)\in\tilde\lambda ^{1/p}\Cal
J_{B'}+pI_{B'\otimes B'}$.

The proposition is proved.
$\square $
\enddemo

7.3. By results of n.5, $G=\Spec A$ where 
the structure of the $O$-bialgebra $A$ is uniquely recovered from 
the following conditions

$$A=B[\theta ],\ \ [p](\tilde\lambda \theta )=[p](\bar
f')\in\tilde\lambda ^pI_B,\ \ 
\delta _{\LT }(\tilde\lambda\theta )=\delta _{\LT }(\bar
f')\in\tilde\lambda I_{B\otimes B}.\tag{7.3.1}$$

Let $A'=A\otimes _OO'$. Introduce $Y,Z\in I_{A'}$ such that 

$$[\tilde\lambda ^{1/p}Z]=[\tilde\lambda \theta ]-\sum\Sb i,l\endSb 
[o_{il}X_i],\ \ 
\tilde\lambda Y=\tilde\lambda ^{1/p}Z+\sum\Sb i,l\endSb o_{il}X_i.
\tag{7.3.2}$$
(The existence of $Y$ and $Z$ follows from the congruences 
$o_{il}\equiv 0\operatorname{mod}(\tilde\lambda ^{1/p})$ 
and 
\linebreak 
$\sum\Sb i,l\endSb [o_{il}X_i]\equiv 
\sum\Sb i,l\endSb o_{il}X_i\operatorname{mod}\tilde\lambda I_{B'}$.)

\proclaim{Proposition 7.3.3} 
\medskip 

{\rm a)} $\delta ^+Z\in I_{A'\otimes
A'}(p)^p+\dsize\frac{p}{\tilde\lambda ^{1/p}}I_{A'\otimes A'}$;
\medskip 

{\rm b)} $\dsize -\frac{1}{p}Z^p\equiv
Y\operatorname{mod}\frac{p}{\tilde\lambda }I_{A'}.$
\endproclaim 

\demo{Proof} By 7.2.1 and 7.3.1 
$\delta _{\LT }(\tilde\lambda ^{1/p}Z)=\delta _{\LT
}(g)\in\tilde\lambda ^{1/p}I_{B'\otimes B'}(p)^p+pI_{B'\otimes B'}$. 
Therefore, the part a) will be implied by the following congruence 
$$[\tilde\lambda ^{1/p}Z\otimes 1]+[1\otimes\tilde\lambda ^{1/p}Z]
\equiv [\tilde\lambda ^{1/p}(Z\otimes 1+1\otimes Z)]
\operatorname{mod}\tilde\lambda ^{1/p}I_{A'\otimes A'}(p)^p.$$
By 4.5.1) this will follow from the congruences 
$$\phi _n(\tilde\lambda ^{1/p}Z\otimes 1, 1\otimes \tilde\lambda
^{1/p}Z)
=\tilde\lambda ^{p^{n-1}}\phi _n(Z\otimes 1,1\otimes Z)
\equiv 0\operatorname{mod}\tilde\lambda ^{1/p}I_{A'\otimes A'}(p)^p,$$
where $n\geqslant 1$.
Because all $\phi _n$ are homogenious polynomials of degree $p^n$ 
it will be enough to prove that for all $1\leqslant r<p$, 
$\tilde\lambda ^{1-1/p}Z^r\otimes Z^{p-r}\in I_{A'\otimes
A'}(p)^p$. 
First, notice that 
the congruence 
$$[p](g)=[p](\tilde\lambda ^{1/p}Z)\equiv 
[\tilde\lambda Z^p]+[p\tilde\lambda
^{1/p}Z]\operatorname{mod}p^2I_A$$
implies by (7.2.1) that 
$\tilde\lambda Z^p\in p\tilde\lambda ^{1/p}I_{A'}$, or equivalently, 
$\tilde\lambda ^{1/p-1/{p^2}}Z\in I_{A'}(p)$. Therefore, 

$$\tilde\lambda ^{1-1/p}Z^r\otimes Z^{p-r}\in 
I_{A'\otimes A'}(p)^p$$
and the part a) is proved. 

Now we can apply 7.2.1 a) to obtain that 
$$-p\sum\Sb i,l\endSb o_{il}X_i\equiv [p](g)\equiv \tilde\lambda
Z^p+p\tilde\lambda ^{1/p}Z\operatorname{mod}p^2I_{A'}$$
and dividing it by 
$p\tilde\lambda $ we obtain the part b) of our proposition. 
$\square $
\enddemo

7.4. Remind that $\Spec B=\Cal G_O(\Cal N)$, where 
$\Cal N\in\MF ^e_S$ is given in the notation from n.3.2. 
Let $S'=S[t']$, where $t^{\prime p}=t$ and let $\Cal N'=\Cal N\otimes
_SS'$. Then $\Cal N'\in\MF _{S'}^{ep}$ and we can extend the
identification $\kappa _{SO}:S\operatorname{mod}t^{ep}\longrightarrow
O\operatorname{mod}p$ to the identification 
$\kappa _{S'O'}:S'/t^{ep}S'\longrightarrow
O'/pO'$. 

For all $i$ and $l$, let $\alpha _{il}\in S'$ be such that 
$\kappa _{S'O'}(\alpha _{il}\operatorname{mod}t^{ep})=
\tilde\lambda ^{-1/p}o_{il}\operatorname{mod}p$.

Then in notation from 1.3 and 3.2 we have the following 

\proclaim{Lemma 7.4.1} $\sum\Sb i,l\endSb \alpha _{il}n_i\in
Z_{\tilde s}(\Cal N')$.  
\endproclaim 

\demo{Proof} It will be sufficient to prove that for all $i$ and $l$, 
$t^e\tilde s^{-1}\alpha _{il}\equiv 0\operatorname{mod}\tilde s_i$. 
Via the identification $\kappa _{S'O'}$ these conditions can be
rewritten in the form 
\linebreak 
$\pi ^e\tilde\lambda ^{-1}o_{il}\equiv 0\operatorname{mod}\tilde\eta
_i^{1/p}$ or, equivalently, 
$$o_{il}\eta _i^{1/p}\equiv 0\operatorname{mod}\tilde\lambda .
\tag{7.4.2}$$
 (Use that $\kappa _{SO}(\tilde s_i\operatorname{mod}t^{ep})
=\tilde\eta _i^{1/p}\operatorname{mod}p$ and 
$\kappa _{SO}(\tilde s\operatorname{mod}t^{ep})
=\tilde\lambda ^{1-1/p}\operatorname{mod}p$.) 

These conditions can be verified as follows. By 7.2.1 it holds 
$$[p](\bar f')=[p](g)+\sum\Sb i,l\endSb [p](o_{il}X_i)\equiv \sum\Sb
i,l\endSb [o_{il}^pX_i^p]\operatorname{mod}p^2I_B.\tag{7.4.3}$$
Therefore, the relation $[p](\bar f')\in\tilde\lambda ^pI_B$ implies
that 
$$\sum\Sb i,l\endSb [o_{il}^pX_i^p]\equiv\sum\Sb i,l\endSb
o_{il}^pX_i^p\equiv 0\operatorname{mod}(\tilde\lambda ^pI_B).$$
(Use that all $o_{il}^p\equiv 0\operatorname{mod}\tilde\lambda $ and
the formulae 4.5.3). Finally, the explicit description 
of the structure of the $O$-algebra $B$ from n.3.2 implies that 
for all $i$ and $l$, 
$o_{il}^p\eta _i\equiv 0\operatorname{mod}(\tilde\lambda ^p)$ and
congruences (7.4.2) are proved. 

The lemma is proved. 
$\square $
\enddemo 

Suppose $\Cal M'=(M^{\prime 0}, M^{\prime 1}, \varphi _1)\in\Ext 
_{\MF _{S'}^{ep}}(\Cal N', \Cal M_{\tilde s}\otimes _SS')$ is given by 
the cocycle 
$\sum\Sb i,l\endSb \alpha _{il}n_i$ from above Lemma 7.4.1. 
Remind that $M^{\prime 0}=(N^0\otimes _SS')\oplus mS'$ and 
$M^{\prime 1}=(N^1\otimes _SS')+m^1S'$, where 
$m^1=\tilde sm+\sum\Sb i,l\endSb \alpha _{il}n_i$ and $\varphi
_1(m^1)=m$. 

Clearly, the correspondences $m^1\mapsto Z\operatorname{mod}
I_{A'}^{DP}$ and $m\mapsto Y\operatorname{mod}I_{A'}^{DP}$, where $Y$
and $Z$ were introduced in 7.3.2, define a unique 
$\Cal F\in\Hom _{\FM _S}(\Cal M',\iota (A'))$ such that 
$\Cal F|_{\Cal N'}$ coincides with the canonical morphism 
$\iota _{B'}:\Cal N'\longrightarrow\iota (B')$. 

Consider $A_1'\in\Cal A(\Cal M')$ such that 
$A_1'=B'[Y_1]=B'[Y_1,Z_1]$, where 
$$Z_1=\tilde \lambda ^{1-1/p}Y_1+\sum\Sb i,l \endSb (\tilde\lambda
^{-1/p}o_{il})X_i,\ \ Z_1^p=-pY_1.$$

Then $G_1'=\Spec A_1'=\Cal G_{O'}(\Cal M')\in\Gr _{O'}$. 

By Proposition 2.4.1, $\Cal F$ gives rise to a unique $O'$-algebra morphism 
$F:A_1'\longrightarrow A'$ such that 
$\iota _{A_1'}\circ \iota (F)=\Cal F$. 

\proclaim{Proposition 7.4.4} $F$ is an isomorphism of coalgebras. 
\endproclaim 

\demo{Proof} Let $\Delta '$ and $\Delta _1'$ be the comultiplications
on $A'$ and $A_1'$. Prove that 
$$F\circ\Delta '=\Delta _1'\circ
(F\otimes F).$$ 

It will be sufficient to prove that the elements 
$\iota _{A_1'}\circ \iota (F\circ\Delta ')$ and 
$\iota _{A_1'}\circ \iota (\Delta _1'\circ (F\circ F)$ of 
$\Hom _{\FM _{S'}}(\Cal M',\iota (A'\otimes A'))$ coincide. 
Clearly, their restrictions to $\Cal N'$ coincide because 
$F|_{B'}=\id $. It remains to note that the both maps 
send $m^1$ to $(Z\otimes 1+1\otimes Z)
\operatorname{mod}I_{A'\otimes A'}^{DP}$. This follows from 
$\delta ^+Z_1\in I^{DP}_{A_1'\otimes A_1'}$ (by the definition of the 
comultiplication on $A_1'$) and $\delta ^+Z\in I_{A'\otimes A'}^{DP}$
by proposition 7.3.3. 

It remains to prove that $F$ is an isomorphism of $O'$-algebras. 

Consider the induced homomorphism of geometric points 
$$F^*:G'(\bar K)\longrightarrow G_1'(\bar K).$$
Then $F^*$ induces the identity map on the common quotient $H'(\bar
K)$ and we have the following two cases:
\medskip 

a) $F^*$ is a group isomorphism;
\medskip 

b) $G_1'(\bar K)=\G_{\tilde\eta }(\bar K)\times\operatorname{Im}F^*$.
\medskip 

In the case a), $F\otimes _{O'}K':A_1'\otimes _{O'}K'\longrightarrow
A'\otimes _{O'}K'$ is an isomoprhism of $K'$-algebras and, therefore, 
$F$ is an embedding of $A_1'$ into $A'$. This embedding is the
identity 
map on the common subalgebra $B'$. Therefore, 
$F(A_1')=A'$ because the differentes 
$\Cal D(A'/B')$ and $\Cal D(A_1'/B')$ coincide. So, 
$F$ is an isomorphism of $O'$-algebras. 

In the case b), $F(A_1')$ is an $O'$-subalgebra in 
$A'$. It contains $B'$ and $\rk _{O'}F(A_1')=\rk _{O'}B'$. 
Therefore, $F(A_1')=B'$ because 
the quotient $A'/B'$ has no $O'$-tosion. 
In particular, the elements $Y$ and $Z$ from 7.3.2 belong to 
$I_{B'}+I_{A'}^{DP}$. This implies that 
$$\tilde\lambda\theta\in\tilde\lambda ^{1/p}I_{B'}+\tilde\lambda
^{1/p}I_{A'}^{DP}, \tag{7.4.5}$$
where $\theta\in A$ was introduced in 7.3.1.
 
Suppose $g\in \G_{\tilde\eta }(O)$, $g\ne 0$. Then $\theta (g)=\lambda v$, 
where $v \in O^*$, and 7.4.3 implies that 
$$\pi ^*=\lambda\tilde\lambda\in\tilde\lambda ^{1/p}I_{O'}^{DP}
\subset I_{O'}^{DP}=(\pi ^*\pi '),$$
where $\pi ^*\in O$ is such that $\pi ^{*p-1}=-p$. The contradiction. 

The proposition is completely proved. 
$\square $
\enddemo 

7.5. It remains to prove that there is an $\Cal M\in\MF ^e_S$ such that
$\Cal M\otimes _SS'=\Cal M'$. 

By Proposition 1.3.5 it will be sufficient to prove that for all $i$, 
$\left (\sum\Sb l\endSb \alpha _{il}\right )\in
S\operatorname{mod}\tilde s_i$, or equivalently, 
for all $i$ and $1\leqslant l<p$, 
$\alpha _{il}\equiv 0\operatorname{mod}\tilde s_i$. 
Applying the identification $\kappa _{S'O'}$ we can replace these
conditions by the following equivalent ones 
$$o_{il}^p\equiv 0\operatorname{mod}\tilde\lambda\tilde\eta _i,\tag{7.5.1}$$
where as earlier, $1\leqslant i\leqslant u$ and $1\leqslant l<p$. 

Remind that we started with $\bar f\in\m (I_{\bar B})$ such that 
$[p](\bar f)\in\tilde\lambda ^pI_{B_0}$ and for $\bar f'$ such that 
$[\bar f']=[\bar f]+[h_1]$, where 
$[h_1]=[f_0]+[h]$ cf. 7.1.3, we have 
$[p](\bar f')\equiv\sum\Sb i,l\endSb 
[o_{il}^pX_i^p]\operatorname{mod}p^2I_B$. Now notice that 
$h_1\in\tilde\lambda I_B$ and, therefore, $h_1^p\in\tilde\lambda
^pI_{B_0}\operatorname{mod}p\tilde\lambda\pi $. This implies that 
$$\sum\Sb i,l\endSb [o_{il}^pX_i^p]\in
I_{B^0}\operatorname{mod}p\tilde\lambda I_B.\tag{7.5.2}$$

\proclaim{Lemma 7.5.3} For any $i$ and $l$, 
$o_{il}^p\equiv\pi ^l\tilde\lambda
u^0_{il}\operatorname{mod}p\tilde\lambda\pi ^l$, 
where all $u^0_{il}\in O_0$. 
\endproclaim 

\demo{Proof} For all $i$ and $l$, $o_{il}=\pi ^{\prime l}u_{il}$ 
with $u_{il}\in O$. Then $o_{il}^p=\pi ^lu_{il}^p\in\tilde\lambda O$
implies that all $u_{il}^p\in\tilde\lambda O$ and , therefore, 
$u_{il}^p\equiv\tilde\lambda u^0_{il}\operatorname{mod}p\tilde\lambda
$, where all $u^0_{il}\in O_0$. 

The lemma is proved. 
$\square $
\enddemo 

As earlier, for all $i$, 
$X_i^p\equiv\eta _i f^0_i\operatorname{mod}p\pi I_B$, 
where all $f^0_i\in I_{B_0}$. Notice also that 
equations for $X_i$ from n.3.2 imply that the residues of 
$f^0_1\operatorname{mod}\pi,\dots ,f^0_u\operatorname{mod}\pi $ 
are linearly independent modulo $\pi I_B$. Now we can rewrite 
condition 
(7.5.2) in the following form 
$$\sum\Sb i,l\endSb [\pi ^lu^0_{il}\eta _if^0_i]\in
I_{B^0}\operatorname{mod}p\tilde\lambda I_B.\tag{7.5.4}$$
Clearly, the terms with $l=0$ already belong to
$I_{B_0}\operatorname{mod}p\tilde\lambda I_B$. Therefore, we can 
assume that in above relation (7.5.4) the index $l$ 
varies from 1 to $p-1$.  

Suppose the ideal in $O$, which is generated by all 
$\pi ^lu^0_{il}\eta _i$, where $1\leqslant l<p$ and 
$1\leqslant i\leqslant u$, equals $\pi ^{c}O$. Therefore,  
$c\in\Bbb N$ and $c\not\equiv 0\operatorname{mod}p$. 
Then the left-hand sum in (7.5.4) belongs to 
$\pi ^{c}I_B\setminus \pi ^{c+1}I_B$. For this reason, it  
belongs to 
$I_{B^0}\operatorname{mod}p\tilde\lambda I_B$ if and only if 
$\pi ^{c}\equiv 0\operatorname{mod}p\tilde\lambda $. 
In other words, all 
$o_{il}^p\eta _i\equiv 0\operatorname{mod}p\tilde\lambda $ 
if $l\ne 0$. This gives conditions (7.5.2) because $\eta _i\tilde\eta
_i =-p$. 

So, the existence of $\Cal M$ is proved and this implies that 
$G=\Cal G_O(\Cal M)$. 
\medskip 
\medskip

\subhead 8. Applications  
\endsubhead 
\medskip 

In this section we prove  
\medskip 

1) that under the same choice of the uniformising element $\pi _0\in
O_0$, our antiequivalence essentially
coincides with Breuil's antiequivalence restricted to the category 
of group schemes killed by $p$; 
\medskip 

2) a criterion for 
a finite $\Bbb F_p[\Gamma _{K_0}]$-module to be isomorphic to 
$G_0(\bar K)$, where $G_0\in\Gr _{O_0}$; 
\medskip 

3) establish via the Fontaine-Wintenberger field-of-norms functor 
a relation between the Galois modules coming from  Faltings's strict
modules and  the Galois modules of the form $G_0(\bar K)|_{\Gamma
_{K_{\infty }}}$, where $G_0\in\Gr _{O_0}$ and 
$$K_{\infty }=K(\{\pi _n\
|\ n\geqslant 0, \pi _{n+1}^p=\pi _n\});$$
\medskip 

4) that a natural duality in the category $\MF ^e_S$ is transformed to
the Cartier duality in $\Gr _{O_0}$ via the functor $\Cal G_{O_0}^O$. 
\medskip 

\subsubhead {\rm 8.1.} Relation to Breuil's antiequivalence  
\endsubsubhead  Denote by $\Br ^O_{O_0}:\MF ^e_S\longrightarrow\Gr _{O_0}$ the
restriction of Breuil's antiequivalence from [Br1] to our categories. 
Notice that by [Br2, Theorem 3.1.1] Breuil's category of filtered
modules over a suitable divided powers envelope of $S$ can be replaced
by $\MF ^e_S$. 
Let $\Br _O:\MF ^e_S\longrightarrow\Gr _O$ be the extension of scalars
of Breuil's functor, i.e. for any $\Cal
M\in\MF ^e_S$, $\Br _O(\Cal M)=\Br ^O_{O_0}(\Cal M)\otimes _{O_0}O$. 

For any $G\in\operatorname{Im}\Cal G_O(=\operatorname{Im}\Br _O)$
introduce 
$\Cal M(G),\Cal M_{\Br }(G)\in\MF ^e _S$ such that 
$\Cal G_O(\Cal M(G))=G$ and $\Br _O(\Cal M_{\Br }(G))=G$. 
Clearly, $\Cal M$ and $\Cal M_{\Br }$ can be considered 
as contravariant functors from $\operatorname{Im}\Cal G_O$ to 
$\MF ^e_S$. 

The essential  coincidence of $\Br ^O_{O_0}$ 
and $\Cal G^O_{O_0}$ will
be proved in the following form. 

\proclaim{Theorem A} For all $G\in\operatorname{Im}\Cal G_O$ 
there are isomorphisms 
$$f(G)\in\Hom _{\MF ^e _S}(\Cal M_{\Br }(G), \Cal M(G)),$$ 
which are functorial in $G$, in other words for all 
$\Pi\in\Hom _{\Gr _O}(G,G')$, it holds  
$\Cal M_{\Br }(\Pi )\circ f(G)=f(G')\circ\Cal M(\Pi )$.
\endproclaim 

\demo{Proof} 
Suppose $H\in\operatorname{Im}\Cal G_O$. Then its $O$-algebra $A(H)$
can be presented in the form $O[X_1,\dots ,X_u]/(f_1,\dots ,f_u)$
where the generators $f_i$, $1\leqslant i\leqslant u$, of the ideal 
$(f_1,\dots ,f_u)$  are the left hand sides of equations from the
beginning of n.3.2. This $O$-algebra is syntomic 
and following [Br1, Lemma 2.3.2] introduce for all $n\geqslant
0$, 
$$A(H)_n=O[\pi ^{p^{-n}}, X_1^{p^{-n}},\ldots ,X_u^{p^{-n}}]/(f_1,\dots
,f_u)$$
and $A(H)_{\infty }=\cup _{n\geqslant 0}A(H)_n$. 
For any flat $O$-algebra $B$ consider the ideal $B(p)=\{b\in B\ |\
b^p\in pB\}$ in $B$ and 
introduce 
$$\theta (B)=(B/B(p)^p, B(p)/B(p)^p,\varphi _1)\in\FM _S,$$
where $\varphi _1$ is induced 
for all $b\in B(p)$, by the correspondences $b\mapsto -b^p/p$ 
and the corresponding $S$-module structure comes from the
identification $\kappa _{SO}$. 
The following proposition is just an adjustment of Lemmas 
3.1.6 and 3.1.7 from [Br1] to our situation. 

\proclaim{Proposition 8.1.1} There is a functorial in
$G,H\in\operatorname{Im}\Cal G_O$ identification of abelian groups 
$G(A(H)_{\infty })=\Hom _{O-alg}(A(G),A(H)_{\infty })=
\Hom _{\FM _S}(\Cal M_{\Br }(G), \theta (A(H)_{\infty })). \square $
\endproclaim

\proclaim{Proposition 8.1.2} Suppose $G,H\in\operatorname{Im}\Cal G_O$
and $\Cal M\in\MF ^e_S$. Then the natural embedding $A(H)\subset
A(H)_{\infty }$ induces the following functorial in $G,H$ and $\Cal M$
identifications 
\medskip 

{\rm a)} $\Hom _{\FM _S}(\Cal M,\theta (A(H))=\Hom _{\FM _S}(\Cal
M,\theta (A(H)_{\infty }))$;
\medskip 

{\rm b)} $\Hom _{O-alg}(A(G), A(H))=\Hom _{O-alg}(A(G),A(H)_{\infty
}).$
\endproclaim 

\demo{Proof} By Lemma 2.4.1 it will be enough to prove a). 
Then we can use the description of $A(H)$ from n.3.1. 
Any element $a'\in A(H)_{\infty }
\operatorname{mod}A(H)_{\infty }(p)^p$ appears as
a $k$-linear combination  of monomials 
$P_{\underline\alpha }=\pi ^{\alpha _0}X_1^{\alpha _1}\ldots
X_u^{\alpha _u}$, where all $\alpha _i\in\Bbb N[1/p]\cup\{0\}$ and 
if the monomial $P_{\underline\alpha }$ appears with a non-zero coefficient then
$P_{\underline\alpha }\in A(H)(p)$, cf. n.3.3.  This implies that 
$\varphi _1(a')\operatorname{mod}A(H)_{\infty }(p)^p$ is a $k$-linear 
combination of $\varphi _1(P_{\underline\alpha })$ for such monomials
$P_{\underline\alpha }$. Notice that if for $n\geqslant 1$, 
$P_{\underline\alpha }\in A(H)_n$ then $\varphi _1(P_{\underline\alpha
})\in A(H)_{n-1}$. Now the proof can be finished by applying these 
arguments to the images of elements of $M^1$, where $\Cal
M=(M^0,M^1,\varphi _1)$. $\square $ 
\enddemo 

So, we have a functorial in $G,H\in\operatorname{Im}\Cal G_O$ 
identification 
$$\Hom _{O-alg}(A(G),A(H))=
\Hom _{\FM _S}(\Cal M_{\Br }(G),\theta (A(H))$$
and, therefore, the induced  functorial identification 
of abelian groups 
$$\Hom _{O-bialg}(A(G),A(H))=
\Hom _{\MF ^e_S}(\Cal M_{\Br }(G), \Cal M(H)).$$

Take $H=G$ and denote by $f(G)$ the morphism from 
$\Hom _{\MF ^e_S}(\Cal M_{\Br }, \Cal M(H))$ which corresponds under
this identification to the identity morphism $\id _{A(G)}$. 

Clearly, $f(G)$ is functorial in $G$. At the same time, 
$f(G)$ is an isomorphism in $\MF ^e_S$. Otherwise, 
there is a proper subgroup scheme $G_1\subset G$ such that the composition of
$f(G)$ with the natural projection $A(G)\longrightarrow A(G_1)$ is
zero, but the corresponding composition of $\id _{A(G)}$ and the
natural projection $A(G)\longrightarrow A(G_1)$ is not equal to the 
counit morphism $A(G)\longrightarrow O$. 

Theorem A is completely proved. $\square $
\enddemo 
\medskip 

\subsubhead {\rm 8.2.} Galois modules 
$G_0(\bar K)$ and $G(\bar K)$ with $G_0\in\Gr _{O_0}$ and 
 $G\in\operatorname{Im}\Cal G_O$
\endsubsubhead    

Suppose $V$ is a finite abelian group killed by $p$ and 
provided with a 
continuous action of $\Gamma _K$. Introduce 
$T(V)=(T(V)^0, T(V)^1, \varphi _1)\in\FM _S$ such that 
$T(V)^0=\Hom ^{\Gamma _K}(V,\bar O/p\bar O)$, 
$T(V)^1=\Hom ^{\Gamma _K}(V, (\pi ^e\bar O)/p\bar O)$ 
and $\varphi _1$ is induced by the map $a\mapsto -a^p/p$, $a\in
\pi ^e\bar O$. As earlier, the corresponding $S$-module structures 
appear via the identification $\kappa _{SO}: 
S/t^{ep}S\longrightarrow O/pO$.

Let $A(V)=\Map ^{\Gamma _K}(V,\bar O)\in\Aug _O$ with the
augmentation ideal 
$I_{A(V)}=$
\linebreak 
$\{a\in A(V)\ |\ a(0)=0\}$. 
Notice that if $A(V)_K:=A(V)\otimes _OK=\Map ^{\Gamma _K}(V,\bar K)$
then $\Spec A(V)_K$ has a natural structure of a finite group scheme
over $K$ and $V$ is the $\Gamma _K$-module of its $\bar K$-points.

 Introduce the functor $\iota ^{(p)}:\Aug _O\longrightarrow \FM
_S$ such that for $A\in\Aug _O$, 
$\iota ^{(p)}(A)=(I_A/I_A(p)^p, I_A(p)/I_A(p)^p, \varphi _1)$, where 
as usually we use the identification $\kappa _{SO}$ to provide 
$I_A/I_A(p)^p$ with an $S$-module structure (notice that 
$I_A(p)^p\supset pI_A$) and $\varphi _1$ is induced by 
the correspondence $a\mapsto -a^p/p$, $a\in I_A(p)$. 
Notice that the embedding $I_A(p)^p\subset I_A^{DP}$ induces 
a strict epimorphism 
$\iota _p^{DP}:\iota ^{(p)}(A)\longrightarrow \iota ^{DP}(A)$ in the
category $\FM _S$. Suppose $\Cal M=(M^0,M^1,\varphi _1)
\in\MF ^e_S$ and $A\in\Cal A(\Cal
M)$. 
By Proposition 2.3.2, 
$\theta _A^{DP}:\Cal M\longrightarrow\iota ^{DP}(A)$ is $\varphi
_1$-nilpotent. Then there is a unique morphism 
$\theta  _A^{(p)}\in\Hom _{\FM _S}(\Cal M,\iota ^{(p)}(A))$ such that 
$\theta _A^{(p)}\circ\iota _p^{DP}=\theta _A^{DP}$ and $\theta _A^{(p)}$
is $\varphi _1$-nilpotent. 

Let $G=\Cal G_O(\Cal M)=\Spec A$. Suppose 
$$\alpha :V\longrightarrow G(\bar K)\in\Hom _{\Bbb F_p[\Gamma
_K]}(V,G(\bar K)).$$
Consider the morphism of augmented $O$-algebras $\alpha
^*:A\longrightarrow A(V)$ given by 
the correspondence $a\mapsto\{a(\alpha (v))\ |\ v\in V\}$. Then 
we obtain the morphism 
\linebreak 
$\alpha _A:=\iota ^{(p)}(\alpha ^*):\iota ^{(p)}(A)\longrightarrow\iota
^{(p)}(A(V))$ in the category $\FM _S$. 
Let $\gamma ^{(p)}=\iota ^{(p)}\circ\alpha _A:\Cal
M\longrightarrow\iota ^{(p)}(A(V))$. Then  
$\gamma ^{(p)}(M^0)$ is contained in 
$$\{a\in I_{A(V)}\operatorname{mod}pI_{A(V)}
\ |\ a(v+v')\equiv a(v)+a(v')\operatorname{mod}p\bar O, \forall
v,v'\in V\}\subset T^0(V)$$
and $\gamma ^{(p)}(M^1)$ is contained in 
$T^1(V)=\pi ^eT^0(V)$. Therefore, 
$\gamma ^{(p)}$ induces a morphism 
$\gamma :\Cal M\longrightarrow T(V)$ in the category $\FM _S$, 
and the correspondence $\alpha\mapsto\gamma $ gives a map 
$$\Cal B:\Hom _{\Bbb F_P[\Gamma _K]}(V,G(\bar K))\longrightarrow \Hom
_{\Cal FM _S}(\Cal M,T(V)).$$

\proclaim{Proposition 8.2.1} Suppose $|V|=|G(\bar K)|$. 
Then $\Cal B$ induces a bijective map 
from the subset of isomorphisms 
in $\Hom _{\Bbb F_p[\Gamma _K]}(V,G(\bar K))$ to the subset of  
$\varphi _1$-nilpotent morphisms in $\Hom _{\FM _S}(\Cal M,T(V))$. 
\endproclaim

\demo{Proof} 
Suppose $\alpha :V\longrightarrow G(\bar K)$ is an 
isomorphism 
of $\Bbb F_p[\Gamma _K]$-modules. Prove that 
$\gamma =\Cal B(\alpha )$ is $\varphi _1$-nilpotent. 

Because $\gamma $ factors through $\theta _A^{(p)}$ and $\theta
_A^{(p)}$ is $\varphi _1$-nilpotent it will be sufficient to prove
that $\alpha _A$ induces a $\varphi _1$-nilpotent morphism on 
$\theta ^{(p)}_A(\Cal M)=(N^0,N^1,\varphi _1)\subset\iota ^{(p)}(A)$. 

$\bullet $\ \ Prove that $T:=\Ker \alpha _A|_{\theta _A^{(p)}(\Cal
M)}\subset N^1$. 

Suppose $n_0\in N^0\setminus N^1$ and $\alpha _A(n_0)=0$. 

Take $m_0\in M^0$ such that $\theta _A^{(p)}(m_0)=n_0$; clearly, 
$n_0\notin N^1$ implies that 
$m_0\notin M^1$. 

Because $\Cal M\in\MF ^e_S$, we can choose for all $i\geqslant 1$,
$s_i\in S$, $m_i^1\in M^0\setminus tM^1$ and $m_i\in M^0\setminus
tM^0$ such that $m_i^1=s_im_{i-1}$ and $m_i=\varphi _1(m_i^1)$. 
Set $n_i^1=s_in_{i-1}$ and $n_i=\varphi _1(n_i^1)$; clearly, all $n_i$
and $n_i^1$ belong to $\Ker\alpha _A$. 

Notice that all $n_i^1\notin tN^1$. (Otherwise, there is $m'\in M^1$
such that $m_i^1-tm'\in\Ker \theta ^{(p)}_A$, but $\Ker\theta
_A^{(p)}\subset tM^1$ because $\theta _A^{(p)}$ is 
$\varphi _1$-nilpotent, cf. 1.2.3, and, therefore, $m_i'\in tM^1$.)  
In particular, 
all $n_i^1$ and $n_i$ are not equal to zero. 

Let $f_0\in I_A$ be such that $f_0\operatorname{mod}I_A(p)^p=n_0$. For
$i\geqslant 1$, choose  
$o_i\in O$ such that 
$\kappa _{SO}(s_i\operatorname{mod}t^{ep})=o_i\operatorname{mod}p$, 
and then by induction on $i$ choose $f_i\in I_A, f_i^1\in I_A(p)$ such that 
$f_i^1=o_if_{i-1}$ and $f_i=-(f_{i-1}^1)^p/p$. Then all 
$f_i\notin pI_A\subset I_A(p)^p$, because 
$f_i\operatorname{mod}I_A(p)^p=n_i\ne 0$. On the other hand, 
$\alpha _A(n_0)=0$ implies that $f_0\in pA(V)$. 
Therefore, all $f_i\in p^{c_i}A(V)$, where 
$c_i=p^i-(1+p+\dots +p^{i-1})\to +\infty $ if $i\to\infty $. This
implies the existence of $i_0\in\Bbb N$ such that 
$f_{i_0}\in pI_A$ and, therefore, $n_{i_0}=0$. The contradiction. 
\medskip

$\bullet $\ \ Prove that $\varphi _1|_T$ is nilpotent.

First, $\varphi _1(T)\subset T$ because $\alpha _A$ commutes with
$\varphi _1$. 

Now suppose $n_0\in T$. Then there is $f_0\in I_A(p)\cap pA(V)$ such that 
$f_0\operatorname{mod}I_A(p)^p=n_0$. Define by induction on
$i\geqslant 1$, $n_i\in T$ and $f_i\in I_A(p)$ such that 
$n_i=\varphi _1(n_{i-1})$ and $f_i=-f_{i-1}^p/p$. 
(Notice that $n_i\in T$ because $\varphi _1(T)\subset T$, and 
$f_i\in I_A(p)$ because $f_i\operatorname{mod}I_A(p)^p=n_i$.) 
As earlier, $f_0\in pA(V)$ implies that $f_i\in p^{c_i}A(V)$ where
$c_i\to +\infty $ if $i\to\infty $. Therefore, 
there is an $i_0\in\Bbb N$ such that $f_{i_0}\in pI_A$ and
$n_{i_0}=0$. 
So, $\gamma =\Cal B(\alpha )$ is $\varphi _1$-nilpotent. 
\medskip  

Now suppose there is a
$\varphi _1$-nilpotent morphism 
$\gamma :\Cal M\longrightarrow T(V)$ 
in the category $\FM _S$. 

Let $\bar A(V)=\Map (V,\bar O)\in\Aug _{\bar O}$ 
with $I_{\bar A(V)}=\{a\in\bar A(V)\ |\ a(0)=0\}$. 
Then the natural embedding of $T(V)$ into $\iota ^{(p)}(\bar A(V))
=(I_{\bar A(V)}/pI_{\bar A(V)}, \pi ^eI_{\bar
A(V)}/pI_{\bar A(V)},\varphi _1)$
allows us to consider $\gamma $ as a $\varphi _1$-nilpotent 
morphism from  
$\Hom _{\FM _S}(\Cal M,\iota ^{(p)}(\bar A(V)))$. 
If $A\in\Cal A(\Cal M)$ then by Proposition 2.4.1
there is a unique morphism of
augmented $\bar O$-algebras $\bar\Cal F:A\otimes _O\bar
O\longrightarrow\bar A(V)$, which corresponds to 
the composition $\gamma ^{DP}$ of $\gamma $ and 
$\iota ^{DP}_p:\iota ^{(p)}(\bar A(V))\longrightarrow\iota ^{DP}(\bar
A(V))$. 

Then:

a) because $T(V)$ is $\Gamma _K$-invariant in $\iota ^{(p)}(\bar
A(V))$, $\bar\Cal F=\Cal F\otimes _O\bar O$, where $\Cal F\in\Hom
_{\Aug _O}(A,A(V))$;

b) because $T^0(V)\subset\Hom (V,\bar O\operatorname{mod}p)$, 
$\bar\Cal F_K=\bar\Cal F\otimes _{\bar O}\bar K$ is a morphism of
$\bar K$-coalgebras $A\otimes _O\bar K\longrightarrow\Map (V,\bar
K)$. 

The above properties a) and b) imply that $\Cal F_K=\Cal F\otimes _OK$
is a morphism of $K$-bialgebras 
$A_K=A\otimes _OK\longrightarrow A(V)\otimes _OK=\Map ^{\Gamma
_K}(V,\bar K)$. 

Let $B:=\Cal F(A)\subset A(V)$ with the induced structure of augmented
$O$-algebra. Then $B$ is a flat $O$-algebra and, if $\Delta
_V:A(V)\longrightarrow A(V\times V)\supset A(V)\otimes _OA(V)$ 
is induced by the addition 
$V\times V\longrightarrow V$, then $\Delta _V(B)\subset B\otimes
_OB$. This implies that $H=\Spec B$ has the induced structure of an
object of the category $\Gr _O$ and $\Spec\Cal F:H\longrightarrow G$
is a closed embedding in $\Gr _O$. Therefore, 
there is an $\Cal N\in\MF ^e_S$ such that 
$\Cal G_O(\Cal N)=H$ and a strictly epimorphic 
$f\in\Hom _{\MF ^e_S}(\Cal M,\Cal N)$ such that $\Cal G_O(f)=\Spec\Cal
F$. From the definition of $f$ it is clear that 
$\gamma ^{DP}:\Cal M\longrightarrow\iota ^{DP}(A(V))$ factors through
$f$ and the corresponding morphism $\Cal N\longrightarrow\iota
^{DP}(A(V))$ is still $\varphi _1$-nilpotent. 
Therefore, $f$ is an isomorphism in $\MF ^e_S$ and $\Cal F_K$ is an
isomorphism of the $K$-bialgebra $A_K$ and a $K$-sub-bialgebra 
of $\Map ^{\Gamma _K}(V,\bar K)$. 

It remains to notice that  
$\rk _KA_K=p^{\rk _SM^0}=p^{\rk _{\Bbb F_p}V}=\rk _KA(V)_K$ implies that 
$\Cal F_K(A_K)=A(V)_K$. So,  
$G(\bar K)\simeq V$ as $\Bbb F_p[\Gamma _K]$-modules. 
Clearly, if $\Cal F_K$ is induced by the $\Bbb F_p[\Gamma
_K]$-isomorphism $\alpha :V\longrightarrow G(\bar
K)$ then $\Cal B(\alpha )=\lambda $. 

Proposition 8.2 is completely proved. 
$\square $
\enddemo 

We can reformulate it as the following criterion (use results of
section 4 for part b)). 

\proclaim{Theorem B} {\rm a)} A finite $\Bbb F_p[\Gamma _K]$-module 
$V$ is isomorphic to $G(\bar K)$, where 
$G\in\operatorname{Im}(\Cal G_O)$, if and only there is an 
$\Cal M=(M^0,M^1,\varphi _1)\in\MF ^e_S$ such that
$\operatorname{dim}_{\Bbb F_p}V=\operatorname{rk}_S(M^0)$, and a
$\varphi _1$-nilpotent morphism from $\Cal M$ to $T(V)$ in the
category $\FM _S$.
\medskip 
{\rm b)} A finite $\Bbb F_p[\Gamma _{K_0}]$-module $V_0$ 
is isomorphic to $G_0(\bar K)$, where $G_0\in\Gr _{O_0}$, if and only if 
$V|_{\Gamma _K}$ appears in the form $G(\bar K)$ with
$G\in\operatorname{Im}\Cal G_O$ and the ramification subgroup  
$\Gamma _{K_0}^{(e^*)}$ with $e^*=ep/(p-1)$ act trivially on $V$.$\square $
\endproclaim 
\medskip 

\remark{Remark} 1) Clearly, in part a) of the above Theorem, $G=\Cal
G_0(\Cal M)$;
\medskip 
2) proposition 8.2.1 is also interesting if $K_0$ is big enough, e.g. 
all points of $G\in\Gr _{O_0}$ are defined over $K_0$; in particular,
it allows to retrieve the main results of the paper [Ab1] 
\endremark 
\medskip 
\medskip 

\subsubhead {\rm 8.3.} Group schemes from $\Gr _{O_0}$ and Faltings's
strict modules  
\endsubsubhead 

Remind that (cf. basic notation) $S=k[[t]]$ and $S_0=k[[t_0]]$ with
$t_0=t^p$. 

\subsubhead {\rm 8.3.1.} Characteristic $p$ analogues of 
$\Cal G_{O}$ and $\Cal G^O_{O_0}$ 
\endsubsubhead 

 Suppose $S_{00}=\Bbb F_p[\tau _{00}]$ 
where $\tau _{00}\notin S_0^*$. 
Then the completion $\hat S_{00}$ is a closed subring in $S_0$ with the residue field $\Bbb
F_p$ and a uniformising element $\tau _{00}$. Consider
the categories $\Gr (S_{00})_{S_0}$ and $\Gr (S_{00})_{S}$ 
of finite flat commutative 
group schemes over $S_0$ and, resp., $S$, which are provided 
with strict action of
$S_{00}$ and are killed by the corresponding action of $\tau _{00}$. 
The general concept of such strict modules was introduced in [Fa] and
was studied in details in [Ab4]. 

As earlier, suppose $\Cal K_0=\Frac S_0$, 
$\Cal K=\Frac S$ and $\Cal K_{00}=\Frac \hat S_{00}$. Then 
the ramification index of $\Cal K_0$ over $\Cal K_{00}$ is $e$. 
The objects of the category 
$\Gr (S_{00})_{S_0}$ (it was denoted by 
$\operatorname{DGr}^{\prime *}_1(S_{00})_{S_0}$ in [Ab4]) can be
described via the antiequivalence 
$\Cal G^S_{S_0}:\MF ^e_S\longrightarrow\Gr (S_{00})_{S_0}$  
as follows (cf. [Ab, 4.5.3], where 
$\MF ^e_S$ was denoted by
$\operatorname{BR}_1(S_{00})_{S_0}$).

Suppose $\Cal M=(M^0,M^1,\varphi _1)\in\MF ^e_S$ is given via 
an $S$-basis $\bar m^1=(m_1^1,\dots ,m_u^1)$ of $M^1$,  
an $S$-basis $\bar m=(m_1,\dots ,m_u)$ of $M^0$   
and $U\in M_u(S)$, such that 
$\varphi _1(\bar m^1)=\bar m$, $\bar m^1=\bar mU$ and $U$ divides the
scalar matrix $(t^e\delta _{ij})\in M_u(S)$. Then we can define the
functor $\Cal G_S$ from $\MF ^e_S$ to 
$\Gr (S_{00})_S$.  By definition, 
$\Cal G_{S}(\Cal M)=\Spec A=\Cal H$, where 
$A=S[\bar X]$, $\bar X=(X_1,\dots ,X_u)$ and 
$$\bar X^{(p)}+\tau _{00}{U^{(p)}}^{-1}\bar X=0.$$
Notice that these equations come from the relation 
$(\bar XU)^{(p)}+\tau _{00}\bar X=0$, which is a complete analogue 
of the corresponding relation from 2.2. 
The coalgera structure on $A$ is given via the counit 
$e:A\longrightarrow S$ such that 
$e(\bar X)=0$ and the comultiplication 
$\Delta :A\longrightarrow A\otimes _{S}A$ 
such that $\Delta (\bar X)=\bar X\otimes 1+1\otimes \bar X$. 
Finally, the action $[r]:A\longrightarrow A$ of any 
$r\in S_{00}$ is uniquely determined by the following conditions 
$[\tau _{00}](a)=0$ and $[\alpha ](a)=\alpha a$ if 
$a\in I_{A}:=\Ker e$ and $\alpha\in\Bbb F_p$. It remains to notice 
that $\Cal H$ appears already as extension of scalars of $S_0$-scheme
$\Cal H_0$, so in this case 
the problem of descent to $S_0$ has a trivial solution and 
the functor $\Cal G_S$ induces the required functor $\Cal G_{S_0}^S$. 

Notice that $A\otimes _{S}\Cal K$ is an etale 
$\Cal K$-algebra and it makes sense to introduce 
the $\Gamma _{\Cal K}=\Aut _{\Cal K}(\bar\Cal K)$-module 
$\Cal H(\bar\Cal K)$ of $\bar\Cal K$ points of $\Cal H$. 

We can see now that the functor $\Cal G_{S}$, resp. $\Cal G_{S_0}^S$, is 
just a simplified characteristic $p$ 
version of the functor $\Cal G_O$, resp. $\Cal G^O_{O_0}$. Nevertheless, 
the functors $\Cal G_S$ and $\Cal G_{S_0}^S$ are not still very far
from the functors $\Cal G_O$ and $\Cal G_{O_0}^O$. 
In Theorem C below we prove that 
for any $\Cal M\in\MF ^e_S$,  
the Galois modules $\Cal G^S_{S_0}(\Cal M)(\bar\Cal K)$ and 
$\Cal G^O_{O_0}(\Cal M)(\bar K)$ can be identified via the 
Fontaine-Wintenberger construction of the field-of-norms functor 
and even more, they can be uniquely recovered one from another.

Suppose  
$\bar S$ is the valuation ring of $\bar\Cal K$. Let $\Cal V$ be a
finite $\Bbb F_p[\Gamma _{\Cal K}]$-module and 
$T(\Cal V)=(T(\Cal V)^0,T(\Cal V)^1,\varphi _1)\in\FM _S$, where 
$T(\Cal V)^0=\Hom ^{\Gamma _{\Cal K}}(\Cal V,\bar
S/\tau _{00}\bar S)$, 
$T(\Cal V)^1=\{a\in T(\Cal V)^0\ |\ a^p=0\}$ and 
$\varphi _1:T(\Cal V)^1\longrightarrow T(V)^0$ is induced by the map 
$s\mapsto -s^p/\tau _{00}$, where $s\in t^e\bar S$.

The following property can be obtained in the same way as 
above Theorem B.

\proclaim{Theorem B$^\prime $} With the above notation 
suppose $\Cal M\in\MF ^e_S$, $\Cal H=\Cal G_{S}(\Cal M)$ and 
$|\Cal V|=|\Cal H(\bar \Cal K)|$. 
Then $\Cal V\simeq\Cal H(\bar \Cal K)$ 
as $\Bbb F_p[\Gamma _{\bar K}]$-modules if and
only if there is a $\varphi _1$-nilpotent morphism 
in $\Hom _{\FM _S}(\Cal M,T(\Cal V))$. $\square $
\endproclaim 
\medskip

\subsubhead{\rm 8.3.2.} Galois modules $G_0(\bar K)$ and the
field-of-norms functor 
\endsubsubhead 

Consider Fontaine's ring $R=\varprojlim \left (\bar O/p\bar
O\right )_n$, where for $n\geqslant 1$, 
the connecting morphisms $(\bar O/p\bar O)_{n+1}
\longrightarrow (\bar O/p\bar O)_n$ are induced by the 
$p$-th power map on $\bar O$. 
Let $R_0=\Frac R$ be the fraction field of $R$. Then $R_0$ is a
complete algebraicly closed valuation field of characteristic $p$, 
the embedding 
of $k$ into $\bar O/p\bar O$ induces a canionical embedding of $k$
into $R_0$. We extend it to the closed embedding of $S$ into $R$ by 
identifying the uniformising element $t$ of $S$ with 
$(\pi _n\operatorname{mod}p)_{n\geqslant 1}\in R$ such that 
$\pi _1=\pi $ and for all $n\geqslant 2$, $\pi _n\in\bar O$ are such
that 
$\pi _n^p=\pi _{n-1}$. Therefore, $\Cal K=\Frac S$ is identified with
a closed subfield in $R_0$ and by 
the Fontaine-Wintenberger theory of the field-of-norms functor,  
$R_0$ coincides with the completion of the algebraic closure 
of $\Cal K$ in $R_0$. Notice also that $\Cal K$ is an inseparable
extension of $\Cal K_0=\Frac S_0$ of degree $p$. 

On the other hand, the absolute Galois group 
$\Gamma _K=\Gal (\bar K/K)$ acts on $R_0$ and this allows to identify
its subgroup 
$\Gamma _{K_{\infty }}=\Gal (\bar K/K_{\infty })$ with the absolute
Galois group $\Gamma _{\Cal K}=\Aut _{\Cal K}(\bar K)$ of $\Cal K$. 
Here $K_{\infty }=\cup _{n\geqslant 0}K_n$ and $K_n=K(\pi _n)$ for all
$n\geqslant 0$. 

Now notice that:
\medskip

a) the above embedding $S\subset R$ induces an embedding $\bar
S\subset R$ (where $\bar S$ is the valuation ring of the algebraic
closure $\bar\Cal K$ of $\Cal K$ in $R_0$) and the identification 
$$\kappa _{\bar S\bar O}:\bar S/t^{ep}\bar S
\longrightarrow R/t^{ep}R=\bar O/p\bar O$$
(use the projection of $R$ to $(\bar O/p\bar O)_1$), which extends 
our original identification $\kappa _{SO}$;
\medskip  

b) with respect to the above identification $\Gamma _{\Cal K}=
\Gamma _{K_{\infty }}$, the identification $\kappa _{\bar S\bar O}$ is
compatible with the action of $\Gamma _{\Cal K}$;
\medskip 

c) suppose $F(T)\in W(k)[T]$ is the minimal monic polynomial 
for $\pi _0\in K_0$ over $K_{00}=\Frac W(k)$. Then 
$F(T)=T^e+p(b_1T^{e-1}+\dots +b_{e-1}T+b_e)$, 
where all $b_i\in W(k)$ and $b_e\in W(k)^*$. Let 
$$\tau _{00}=-t_0^e(\bar b_e+\bar b_{e-1}t_0+\dots +\bar
b_1t_0^{e-1})^{-1}\in S_0,$$
where all $\bar b_i:=b_i\operatorname{mod}p\in k$. Then 
$\Cal K_{00}=k[[\tau _{00}]]$ is a closed subfield 
in $\Cal K_0$ and $\Cal K_0$ is a totally ramified extension 
of $\Cal K_{00}$ of degree $e$; 
\medskip 

d) suppose $s\in\bar S$ and $o\in\bar O$ are such that 
$\kappa _{\bar S\bar
O}(s\operatorname{mod}t^{ep})=o\operatorname{mod}p$; 
then $s\in t^e\bar S$ implies that $o\in\pi ^e\bar O$ and 
$$\kappa _{\bar S\bar O}\left ((-s^p/\tau
_{00})\operatorname{mod}t^{ep}\right )=
\left (-o^p/p\right )\operatorname{mod}p.$$
Indeed, it will be sufficient to verify this formula for $s_0=t^e$ and
$o_0=\pi ^e$; then $-t^{pe}/\tau _{00}=
\bar b_e+\bar b_{e-1}t_0+\dots +\bar b_1t_0^{e-1}$ and 
$$\kappa _{\bar S\bar O}(-t^{pe}/\tau _{00}\operatorname{mod}t^{pe})
=(b_e+b_{e-1}\pi _0+\dots +b_1\pi _0^{e-1})\operatorname{mod}p=
(-\pi ^{ep}/p)\operatorname{mod}p.$$
\medskip 

\proclaim{Theorem C} Suppose $\Cal M\in\MF^e_S$, 
$H_0=\Cal G^O_{O_0}(\Cal M)$ and 
$\Cal H_0=\Cal G^S_{S_0}(\Cal M)\in\Gr (S_{00})_{S_0}$. Then 
\medskip 

{\rm a)} with respect to the field-of-norms identification 
$\Gamma _{\Cal K_0}=\Gamma _{K_{\infty }}\subset\Gamma _{K_0}$, 
the $\Gamma _{\Cal K_0}$-modules 
$H_0(\bar K)|_{\Gamma _{K_{\infty }}}$ and $\Cal
H_0(\bar\Cal K)$ are isomorphic; 
\medskip 

{\rm b)} the $\Gamma _{K_0}$-module $V_0=H_0(\bar K)$ can be uniquely recovered
from the $\Gamma _{\Cal K_0}$-module $\Cal H_0(\bar\Cal K)$.
\endproclaim 

\demo{Proof} Suppose $V$ is a finite $\Bbb F_p[\Gamma
_K]$-module and $\Cal V=V|_{\Gamma _{K_{\infty }}}$ is the $\Bbb F_p[\Gamma
_{\Cal K}]$-module with respect to our identification $\Gamma _{\Cal
K}=\Gamma _{K_{\infty }}$. Then the embedding 
$$\Hom ^{\Gamma _K}(V,\bar O/p\bar O)\longrightarrow\Hom
^{\Gamma _{K_{\infty }}}(V, \bar O/p\bar O)$$
together with the identification $\kappa _{\Bar S\bar O}$ induce 
the embedding $\omega :T(V)\longrightarrow T(\Cal V)$ in the category
$\FM _S$. 

Now notice that 
if $V=H(\bar K)$ then there is a $\varphi _1$-nilpotent morphism 
$\gamma\in\Hom _{\FM _S}(\Cal M,T(V))$. Therefore, 
for $\Cal V=\Cal H(\bar\Cal K)$, 
$\gamma\circ\omega _*\in\Hom _{\FM _S}(\Cal M, T(\Cal V))$ is also
$\varphi _1$-nilpotent and, therefore, 
$\Cal V\simeq\Cal H(\bar \Cal K)$. This proves the part a), because 
$\Gamma _{\Cal K_0}=\Gamma _{\Cal K}=\Aut _{\Cal K_0}(\bar \Cal K)$. 
($\Cal K$ is purely inseparable over $\Cal K_0$.) 

In order to prove b) let $e^*=ep/(p-1)$ and 
notice that by Fontaine's estimates, the ramification subgroup 
$\Gamma _{K_0}^{(e^*)}$ acts trivially on $V_0$ and by estimates from
[Ab4], $\Gamma _{\Cal K_0}^{(e^*)}$ acts trivially on $\Cal V_0$. 

Therefore, it will be sufficient to prove that the field-of-norms
embedding 
$\Gamma _{\Cal K_0}=\Gamma _{K_{\infty }}\subset\Gamma _{K_0}$ induces 
group isomorphism $\Gamma _{\Cal K_0}/\Gamma _{\Cal K_0}^{(e^*)}\simeq 
\Gamma _{K_0}/\Gamma _{K_0}^{(e^*)}$ or, equivalently, we have the
following two properties:
\medskip 

1) $\Gamma _{\Cal K_0}^{(e^*)}=\Gamma _{\Cal K_0}\cap\Gamma
_{K_0}^{(e^*)}$; 
\medskip 

2) $\Gamma _{K_0}=\Gamma _{\Cal K_0}\Gamma _{K_0}^{(e^*)}$. 
\medskip 

Now remind that the ramification theory attaches to any finite
separable extension of complete discrete valuation fields  with
perfect residue fields $L\supset E$, 
the Herbrand function $\varphi _{L/E}(x)$,
$x\geqslant 0$. One can extend the well-known definition of Herbrand's function
for Galois extensions from [Se] by the use of the composition property 
$\varphi _{L'/E}(x)=\varphi _{L/E}(\varphi _{L'/L}(x))$, 
where $L'$ contains $L$ and is Galois over $E$. Alternatively, 
the Appendix to [De] contains a direct construction of the ramification
theory for arbitrary (not necessarily Galois) separable finite extensions $L/E$ 
by using $E$-embeddings of $L$ into a fixed algebraic closure $\bar E$
of $E$. Anyway, such Herbrand's function satisfies the following two
propertis:
\medskip

$\bullet $ \ \ if $L_1\supset L\supset E$ are finite separable 
field extensions
then 
for all $x\geqslant 0$, 
$$\varphi _{L_1/L}(x)=\varphi _{L/E}(\varphi _{L_1/L}(x));$$
\medskip 

$\bullet $\ \ the ramification subgroup $\Gamma _E^{(v)}$ acts
trivially on $L$ if and only if $v>v(L/E)$, where $v(L/E)$ is the
value of $\varphi _{L/E}$ at its last edge point. (By definition, 
$(0,0)$ is always an edge point of Herbrand's function of any ramified
extension.)
\medskip 

Now notice that for any $n\geqslant 0$, 

$$\varphi _{K_{n+1}/K_n}(x)=\cases x, &\text{if } 0\leqslant
x\leqslant e^*p^n \\ 
e^*p^n+(x-e^*p^n)/p, &\text{if }x\geqslant e^*p^n
\endcases 
$$ 
Therefore, if 
$$\varphi _{K_{\infty }/K_0}(x)=\lim _{n\to\infty }
(\varphi _{K_n/K_{n-1}}\circ\ldots \circ\varphi _{K_1/K_0})(x)$$
then $\varphi _{K_{\infty }/K_0}(e^*)=e^*$ and by the corresponding
property of the field-of-norms functor [Wi, Cor. 3.3.6], it holds 
$\Gamma _{\Cal K_0}^{(e^*)}=\Gamma _{K_{\infty }}\cap\Gamma
_{K_0}^{(\varphi _{K_{\infty }/K_0}(e^*))}=
\Gamma _{K_{\infty }}\cap\Gamma _{K_0}^{(e^*)}$. 
This proves the property 1). 

Prove the property 2). Suppose $L$ is the subfield of $\bar K$ 
fixed by $\Gamma _{\Cal K_0}\Gamma _{K_0}^{(e^*)}=
\Gamma _{K_{\infty }}\Gamma _{K_0}^{(e^*)}$. Then $L$ is a 
finite extension of $K_0$ in $K_{\infty }$ and $v(L/K_0)<e^*$. 
If $L\ne K_0$ then there is an $s\geqslant 0$ such that 
$LK_s=K_{s+1}$ (use that for all $n\geqslant 0$, $[K_{n+1}:K_n]=p$). 

Notice that for all $n\geqslant 1$, $v(K_n/K_0)=e^*+e(n-1)$ 
(use that $\varphi _{K_n/K_0}=\varphi _{K_n/K_{n-1}}\circ\ldots
\circ\varphi _{K_1/K_0}$). Therefore, $s\ne 0$ (otherwise, 
$LK_0=L=K_1$ but $v(L/K_0)<v(K_1/K_0)=e^*$). But if $s\geqslant 1$ then 
$$e^*+es=v(K_{s+1}/K_0)=\max (v(L/K_0), v(K_s/K_0))=e^*+e(s-1).$$
The contradiction. So, $L=K_0$ and the property 2) is proved. 

Theorem C is completely proved. $\square $
\enddemo 

\subsubhead{\rm 8.3.3.} Full faithfulness of the restriction from
$\Gamma _{K_0}$ to $\Gamma _{K_{\infty }}$ 
\endsubsubhead 

The above methods can be applied to study a more general situation. 

Suppose $\Cal C_{K_0}$ is a full 
 subcategory of the category of finite $p$-torsion modules with
 continuous action of $\Gamma _K$. Let $\MG _{K_{\infty }}$ be the
 category of $\Gamma _{K_{\infty }}$-modules. Then we have the functor 
$\Cal F:\Cal C_{K_0}\longrightarrow\MG _{K_{\infty }}$ of restriction 
of the action of $\Gamma _{K_0}$ to the action of $\Gamma _{K_{\infty }}\subset
 \Gamma _{K_0}$. For $n\in\Bbb N$, let $\Cal C^{(n)}_{K_0}$ be the full subcategory in
 $\Cal C_{K_0}$ consisting of modules killed by $p^n$.  

\proclaim{Theorem C$^{\prime }$} Suppose for any $H\in\Cal C^{(1)}_{K_0}$ the
ramification subgroups $\Gamma _{K_0}^{(e^*)}$ act trivially on
$H$. Then $\Cal F$ is fully faithful.  
\endproclaim

\demo{Proof} For $n\in\Bbb N$, let $\Cal F^{(n)}$ be the restriction of $\Cal F$ to 
$\Cal C^{(n)}_{K_0}$. Then we can proceed as in 8.3.2 to deduce from 
$\Gamma _{K_{\infty }}\Gamma
_{K_0}^{(e^*)}=\Gamma _{K_0}$  
that $\Cal F^{(1)}$ is fully faithful. Now notice that for any
$H_1,H_2\in\Cal C_{K_0}$ there is a short exact sequence 
$0\longrightarrow \Ker (p\id _{H_2})\longrightarrow
H_2\longrightarrow \operatorname{Im}(p\id _{H_2})\longrightarrow 0$ 
and we can use a devissage procedure (based on the standard 6-terms
$\Hom -\Ext $ exact sequence) to deduce by induction that 
$\Cal F^{(n)}$ is also fully faithful. $\square $ 
\enddemo 

This theorem can be applied in the following cases: 

$\bullet $ \ \ if $\Cal C _{K_0}$ is the category of 
Galois modules $G_0(\bar K)$ where $G_0$ is an arbitrary finite flat
commutative $p$-group scheme over $O_0$ we 
retrieve Breuil's result [Br3, Theorem 3.4.3]; 

$\bullet $  if  $K_0$ is unramified over $\Bbb Q_p$  then we can apply
Theorem $C'$ to the category  
of all finite subquotients 
of crystalline $\Bbb Z_p[\Gamma _{K_0}]$-modules 
with Hodge-Tate weights of length $<p$ because of  
the ramification estimates from [Ab5] (if the above length is
$\leqslant p-2$ we retrieve the main result of [Br3] where it is
sufficient to use Fontaine's ramification estimates from [Fo5]).
\medskip

8.4. {\it Cartier duality.} In this subsection we prove that if $\Cal
N\in\MF ^e_S$ and $\widetilde{\Cal N}\in MF ^e_S$ is its dual (cf. the definition below)
then $\Cal G_{O_0}^O(\Cal N)$ and $\Cal G_{O_0}^O(\widetilde{\Cal N})$ are
Cartier dual group schemes. Clearly, it will be 
sufficient to verify this over $O$, i.e. that $G=\Cal G_O(\Cal N)$ and $\widetilde{G}=\Cal
G_o(\widetilde{\Cal N})$ are Cartier dual. We are going to prove this 
by constructing a non-degenrate
bilinear pairing of group functors $G\times
\widetilde{G}\longrightarrow\mu _{p,O}$, where as usually $\mu _{p,O}$
is the constant multiplicative group scheme of order $p$ over $O$. 

Let $-p=\pi _0\varepsilon _0$, where $\varepsilon _0\i O_0^*$. 
Let $\omega _0\in S_0$ be such that 
$\kappa _{SO}(\omega _0\operatorname{mod}t^{ep})=
\varepsilon _0\operatorname{mod}p$. Define the $\sigma $-linear
morphism $\varphi _1:t^eS\longrightarrow S$ by the relation 
$\varphi _1(t^es)=\omega _0\sigma (s)$, $s\in S$. Notice that 
$\varphi _1(t^e\sigma ^{-1}(\omega _0))=1$ and 
$\Cal S=(S,t^eS,\varphi _1)\in\MF ^e_S$. 

Suppose $\Cal N=(N^0, N^1,\varphi _1)\in\MF ^e_S$. 

\definition{Definition} Let $\widetilde{\Cal N}=(\widetilde{N}^0,
\widetilde{N}^1,\varphi _1)\in \FM _S$ be such that 
\medskip 

a) $\widetilde{N}^0=\Hom _S(N^0,S)$;
\medskip 

b) $\widetilde{N}^1=\{f\in \widetilde{N}^0\ |\ f(N^1)\subset t^eS\}$;
\medskip 

c) for any $f\in \widetilde{N}^1$, $\varphi _1(f)\in \widetilde{N}^0$ is such that 
for any $n\in N^1$, $\varphi _1(f)(\varphi _1(n))=\varphi _1(f(n))$
(cf. the above definition of $\varphi _1|_{S}$).
\enddefinition

\remark{Remarks} 1) The condition c) determines $\varphi _1$ uniquely
because $\varphi _1 (N^1)S=N^0$;

2) one can verify easily that $\widetilde{\Cal N}\in\MF ^e_S$; 

3) In the above definition $\widetilde{\Cal N}\operatorname{mod}t^{ep}$ 
does not depend on a choice of $\omega _0$; therefore, 
for different choices of $\omega _0$, the corresponding objects 
$\widetilde{\Cal N}=\widetilde{\Cal N}(\omega _0)$ are related via unique isomorphisms in the
category $\MF ^e_S$ as different $\varphi _1$-nilpotent lifts of $\widetilde{\Cal
N}\operatorname{mod}t^{ep}$.
\endremark 

\proclaim{Theorem D} With the above notation, 
$\widetilde{H}=\Cal G_O(\widetilde{\Cal N})$ is the Cartier dual to $H=\Cal
G_O(\Cal N)$. 
\endproclaim  

\demo{Proof} Suppose $\Cal N$ is given in notation similar to those
from n.3.1. Then we have 
\medskip  

--- the vector $n=(n_1,\dots ,n_u)$ consisting of elements of an 
$S$-basis of $M^0$;
\medskip 

--- for $1\leqslant i\leqslant u$, 
there are $\tilde s_i'\in S$ such that all $\tilde s_i'|t^e$ and 
the vector $n^1=(n_1^1,\dots ,n_u^1)=(\tilde s_1'n_1,\dots ,\tilde s_u'n_u)$ 
consists of elements of an $S$-basis of $N^1$;
\medskip 

--- there is a matrix $U\in \GL _u(S)$ such that $\varphi _1(n^1)=nU$. 
\medskip 

Then one can verify that $\widetilde{\Cal N}=(\widetilde{N}^0,\widetilde{N}^1,\varphi _1)$ can be
described via the following data: 
\medskip 

--- the vector $\tilde n=(\tilde n_1,\dots ,\tilde n_u)$ consisting of elements of 
the $S$-basis of $\widetilde{N}^0$ which is dual to the basis $n_1,\dots ,n_u$ of
$N^0$;
\medskip 

--- the vector $\tilde n^1=(\tilde n_1^1,\dots ,\tilde n_u^1):=(s_1'\tilde
n_1,\dots ,s_u'\tilde n_u)$
consisting of elements of an $S$-basis of $\widetilde{N}^1$, where for
$1\leqslant i\leqslant u$, $s_i'=t^e\sigma ^{-1}(\omega _0)(\tilde
s_i')^{-1}$; 
\medskip 

--- the relation 
$\varphi _1(\tilde n^1)=\tilde n\widetilde{U}$, 
where $\widetilde{U}^t=U^{-1}$ (here $\widetilde{U}^t$ is the
transposed to $\widetilde{U}$ matrix from $\GL _u(S)$).
\medskip 

For $1\leqslant i\leqslant u$, let $\tilde\eta '_i, \eta _i'\in O$ be
such that $\kappa _{SO}(\tilde s_i'\operatorname{mod}t^{ep})=
\tilde\eta '_i\operatorname{mod}p$ and 
$\kappa _{SO}(s_i'\operatorname{mod}t^{ep})=
\eta '_i\operatorname{mod}p$. Set 
$\eta _i=-p/\tilde\eta ^{\prime p}_i$ and $\tilde\eta _i=-p/\eta
_i^{\prime p}$. 

Then $A=A(H)=O[X_1,\dots ,X_u]$ with the equations 
$X_i^p=\eta _i\sum \Sb j\endSb X_jc_{ji}$, 
where $1\leqslant i\leqslant u$ and $C=(c_{ji})\in\GL _u(O)$ is such
that 
$\kappa _{SO}(U\operatorname{mod}t^{ep})=C\operatorname{mod}p$. 

Similarly, $\widetilde{A}=A(\widetilde{H})=O[\widetilde{X}_1,\dots
,\widetilde{X}_u]$ with the equations 
$\widetilde{X}_i^p=\sum\Sb j\endSb \widetilde{X}_j\tilde c_{ji}$, 
where $1\leqslant i\leqslant u$ and $\widetilde{C}=(\tilde
c_{ji})\in\GL _u (O)$ is such that 
$\kappa _{SO}(\widetilde{U}\operatorname{mod}t^{ep})=
\widetilde{C}\operatorname{mod}p$.

\proclaim{Lemma 8.4.1} Let $Z=\sum\Sb i\endSb 
X_i\otimes\widetilde{X}_i\in I_{A\otimes\widetilde{A}}$. Then 
$$Z^p+pZ\equiv 0\operatorname{mod}
(pI_{A\otimes\widetilde{A}}(p)^p+p^2I_{A\otimes\widetilde{A}}).$$
\endproclaim 

\demo{Proof} Use that for all $1\leqslant i\leqslant u$, 
$\eta _i\tilde\eta _i\equiv
-p\operatorname{mod}p^2$, 
$X_i\otimes \widetilde{X}_i\in I_{A\otimes\widetilde{A}}(p)$ and 
$C\widetilde{C}^t\equiv E\operatorname{mod}p$, 
where $E$ is the unit matrix of order $u$. $\square $
\enddemo 

Let $\widetilde{\Cal S}=(Sm, Sm,\varphi _1)\in\MF ^e_S$  be such that $\varphi
_1(m)=m$. Then $\Cal G_O(\widetilde{\Cal S})=\mu _{p,O}$ is the constant
multiplicative group scheme of order $p$ over $O$ with the algebra 
$A(\mu _{p,O})=O[X]$, $X^p+pX=0$. 

Now notice that the correspondence $m\mapsto
Z\operatorname{mod}I_{A\otimes\widetilde{A}}^{DP}$ determines the
morphism 
$q\in\Hom _{\FM _S}(\Cal S, \iota ^{DP}(A\otimes\widetilde{A}))$.
Therefore, by Lemma 2.4.1, $q$ is induced by 
a unique morphism of $O$-algebras 

$$e^*:A(\mu _{p,O})\longrightarrow A(G)\otimes A(\widetilde{G})$$

Clearly, the definitions of the coalgebra structures on $A(H)$ and
$A(\widetilde{H})$ from n.2.4, immediately imply that $e^*$ is 
co-bilinear.   
It remains to verify that $e^*$ gives a non-degenerate
pairing of group functors. 

We can assume that $K$ is so large that all $\bar K$-points of 
group schemes $H$, $\widetilde{H}$ and $\mu _{p,O}$ are defined over
$K$. Then it will be sufficient to verify that  
if $\tilde h_0\in \widetilde{H}(O)$ is such that 
for any $h\in H(O)$, 

$$e^*(X)(h, \tilde h_0)=e_{\mu _{p,O}}(X)=0\tag{8.4.2}$$ 
then $\tilde h_0=0$. (Here $e_{\mu _p,O}:A(\mu _{p,O}\longrightarrow O$
is the counit map.) 

For all $i$, denote by $\bar X_i\in\Hom (H(O), O/pO)$ the images of 
\linebreak 
$X_i\in A(H)\subset \Map (H(O),O)$ with respect to the natural maps 
$$\Map (H(O),O)\longrightarrow \Map (H(O),O/pO)\supset \Hom
(H(O),O/pO)$$
(use that $\delta ^+X_i\in I_{A\otimes A}(p)^p$ implies that 
$\bar X_i\in\Hom (H(O),O/pO)$). By the results of n.3.4, the
generated by $\bar X_i$, $1\leqslant i\leqslant u$, 
$O$-submodule $\Cal H(H)$ in $\Hom (H(O),O/pO)$ can be defined in an
invariant way just in terms of the image of the corresponding
$S$-module $N^0$ in $I_A/I_A(p)^p$. This module can't be too small, 
for $\pi _0^*\in O$ such that $v_p(\pi _0^*)=1/(p-1)$, it holds 
$$\Cal H(H)\supset \pi _0^*\Hom (H(O),O/pO). \tag {8.4.3 } $$
(Use the embedding of $O$-bialgebras $A(H)\supset A(\mu
_{p,O})^{\otimes u}$ and that $\Cal H(\mu _{p,O}^{\times u})=
\pi _0^*\Hom (H(O),O/pO)$.) 

\proclaim{Lemma 8.4.4} 
\medskip 

{\rm a)} $e^*(X)\equiv
Z\operatorname{mod}(I_{A\otimes\widetilde{A}}(p)^p)$;
\medskip 

{\rm b)} $e^*(X)\equiv
Z+Z'\operatorname{mod}(I_{A\otimes\widetilde{A}}(p)^{2p-1})$, where 
$Z'\in I_{A\otimes\widetilde{A}}(p)$ is an $O$-linear combination 
of the terms $X_{i_1}\dots X_{i_p}\otimes 
\widetilde{X}_{i_1}\dots \widetilde{X}_{i_p}$ for all $1\leqslant
i_1,\dots ,i_p\leqslant u$.  
\endproclaim 

\demo{Proof} The part a) follows from Lemma 8.4.1 and $b)$ is obtained
from a) and the relation $-e^*(X)^p/p=e^*(X)$. $\square $
\enddemo 

Now 8.4.2 and 8.4.4 a) imply that in $\Hom (H(O),O/pO)$ it holds 
$$\sum\Sb i\endSb \widetilde{X}_i(\tilde h_0)\bar X_i=0.$$
Then (8.4.3) implies that all $\widetilde{X}_i(\tilde h_0)\equiv
0\operatorname{mod}(p/\pi _0^*)$. 

If $p\geqslant 5$ then $v_p(\widetilde{X}_i(\tilde h_0))
\geqslant (p-2)/(p-1)>1/(p-1)$ and by Lemma 2.4.1 we can conclude that 
$\tilde h_0=0$. In order to finish the proof in general case, just use
that 
for all $i$, $\widetilde{X}_i(\tilde h_0)\in\pi _0^*O$. This imples
that, cf. Lemma 8.4.4 b), 
 $Z'(h,\tilde h_0)\in p\pi _0^*O$ and, therefore, 
$\sum\Sb i\endSb \widetilde{X}_i(\tilde h_0)X_i\in p\pi _0^*\Map
(H(O),O)$. Therefore, 
in $\Hom (H(O),O/pO)$ it holds 
$$\sum\Sb i\endSb (\widetilde{X}_i(\tilde h_0)/\pi _0^*)\bar X_i=0.$$
As earlier, this implies that all $\widetilde{X}_i(\tilde h_0)/\pi
_0^*
\in \pi _0^*O$, therefore, 
\linebreak 
$v_p(\widetilde{X}_i(\tilde h_0))\geqslant
2/(p-1)>1/(p-1)$ and $\tilde h_0=0$. 

Theorem $D$ is completely proved. $\square $

\enddemo

\enddocument